\documentclass[12pt,oneside,english,a4paper]{amsart}

\usepackage[T1]{fontenc}
\usepackage[utf8]{inputenc}
\usepackage[english]{babel}
\usepackage[dvipsnames]{xcolor}
\usepackage[yyyymmdd,hhmmss]{datetime}
\usepackage[numbers,sort,compress]{natbib}
\usepackage{verbatim,todonotes}
\usepackage{bm}
\usepackage{bbm}
\usepackage{amstext,amsthm,amssymb}
\usepackage{xparse,mathrsfs,mathtools,float}
\usepackage[pdftex,unicode=true,%
pdfusetitle,pdfdisplaydoctitle=true,%
pdfpagemode=UseOutlines,%
bookmarks=true,bookmarksnumbered=false,bookmarksopen=true,bookmarksopenlevel=2,%
breaklinks=true,%
pdfencoding=unicode,psdextra,
pdfcreator={},
pdfborder={0 0 0}
]{hyperref}
\usepackage{lmodern,microtype}
\usepackage{geometry}
\usepackage[capitalise,nameinlink]{cleveref}

\newtheorem{thm}{Theorem}[section]
\newtheorem{lem}[thm]{Lemma}
\newtheorem{cor}[thm]{Corollary}
\newtheorem{prop}[thm]{Proposition}
\theoremstyle{remark}
\newtheorem{rem}[thm]{Remark}

\numberwithin{equation}{section}

\newcommand{\sumstar}{\sideset{}{^*}\sum}

\RequirePackage{silence}
\title[Low moments of extended Rademacher sums]{Low moments of sums of extended Rademacher\\random multiplicative functions}
\author{Christopher Atherfold, Maxim Gerspach and Jad Hamdan}

\date{\today}

\begin{document}

\begin{abstract}
    Let $X^e \colon \mathbb{N} \to \{\pm 1\}$ be the random, completely multiplicative function whose values at the primes $(X^e(p))_{p \text{ prime}}$ are independent and uniformly distributed on $\{ \pm 1 \}$. Let $X=X^e\mu^2$. We determine the order of magnitude of 
    \[
    \mathbb{E} \Bigg[ \bigg| \frac{1}{\sqrt{x}}\sum_{n \le x} X^e(n) \bigg|^{2 q} \Bigg]\quad  \text{ and }\quad \mathbb{E} \Bigg[ \bigg|\sum_{n \le x} \frac{X(n)}{\sqrt{n}} \bigg|^{2 q} \Bigg] 
    \]
    for $q\in [0,1]$, showing that a transition from better-than-square-root cancellation to growth like a power of $\log x$ occurs when $q$ crosses $1/2$. In particular, the first absolute moments are of order $(\log\log x)^{3/4}$.

    The proof proceeds by comparing these sums to integrals of random Euler products, which are then studied by analogy with Gaussian multiplicative chaos integrals. The main difficulty lies in proving the lower bound at $q=1/2$, for which we also adapt arguments from the study of random walks.
\end{abstract}

\maketitle

\section{Introduction}

Let $X^e \colon \mathbb{N} \to \{\pm 1\}$ be the (random) completely multiplicative function whose values at the primes $(X^e(p))_{p \text{ prime}}$ are independent and uniformly distributed on $\{ \pm 1 \}$. This random function is sometimes referred to as an extended Rademacher random multiplicative function. We further denote by $X \colon \mathbb{N} \to \{ \pm 1, 0 \}$ the Rademacher random multiplicative function $ X := \mu^2 X^e$, which we have chosen here to be defined by the same realisation of the values of $X^e$ at the primes.

This paper is primarily concerned with the study of low moments (that is, the $2q$--th moment for $q\in [0,1]$) of the partial sums
\[ \sum_{n \le x} X^e(n), \]
which are meant to model the statistical properties of $\sum_{n\leq x} \chi(n)$ for a random quadratic Dirichlet character $\chi$. To see how, choose a large parameter $D$, and let $d$ be chosen uniformly at random among fundamental discriminants in $[D,2D]$. Then the (random) quadratic character $\chi=\chi_d$ associated to $d$ is completely multiplicative, and it is straightforward to verify that the distribution of $\chi(p)$ for a given prime $p$ is approximately uniform on $\{\pm 1\}$. One can also show that $\chi(p)$ and $\chi(p')$ are approximately independent for different primes $p$ and $p'$, so long as they are not too large compared to $D$.

The study of random character sums is usually much more tractable than that of $\sum_{n \le x} \chi_d(n)$
for a fixed $\chi_d$, as it gives access to various tools which are not a priori applicable to the latter. Indeed, we are essentially only able to detect non-trivial cancellation in $\sum_{n\le x} \chi_d(n)$ when $x \ge d^{1/4+o(1)}$, using Burgess' bound \cite{Burgess1963}. Even then, the resulting bound is much weaker than the ``square-root cancellation'' one might expect (meaning a bound of the shape $\ll x^{1/2+\varepsilon}$), with such a bound only being known for $x\geq d^{1-o(1)}$, as a consequence of the P\'olya--Vinogradov inequality (cf.\ \cite[Chapter 9.4]{MV07}). By contrast, square-root cancellation is known to hold in the random setting: uniformly in all $x$ and $D$, one has 
\begin{equation}\label{eq:armon}
   \frac{1}{\#\{ D < d \le 2D : d \text{ fund. disc.}\}} \sumstar_{D < d \le 2D} \bigg| \frac{1}{\sqrt{x}} \sum_{n \le x} \chi_d(n) \bigg|^2 \ll \log D,
\end{equation}
where the sum ranges over fundamental discriminants \cite{Armon}. 

Averaging over the character also allows for comparison with random models such as $X$ and $X^e$, an approach pioneered in \cite{gs01}. Furthering a line of investigation initiated by Davenport and Erd\H{o}s \cite{DavenportErdos}, recent works of Harper \cite{HarperShort} and Nosal \cite{n25} for instance 
used the extended model $X^e$ to study the distribution of quadratic character sums with a random starting point. Both works rely on establishing an estimate of the form
\[
    \frac{\log P}{P^\eta} \sum_{\substack{P\leq p\leq P+P^\eta\\p \text{ prime}}}\Big|\sum_{n\leq N}\alpha_n\Big(\frac{n}{p}\Big)\Big|^2 \ll \mathbb{E}\Big[\Big|\sum_{n\leq N} \alpha_n X^e(n)\Big|^2\Big]+\frac{N}{P^{\eta/2}}\Big(\sum_{n\leq N} |\alpha_n|\Big)^2,
\]
for $\eta\in [0.01,1]$ and $(\alpha_n)_n$ a sequence of complex numbers (see Number Theory Result/Fact 3 in \cite{HarperShort,n25}).

There are several challenges related to the study of the extended Rademacher model compared to the other two main models (namely Steinhaus and Rademacher), chief among them  the lack of orthogonality: whereas $\mathbb{E}[X(n) X(m)] = \mathbf{1}(n=m)$ on squarefree numbers $n, m$, 
\[ \mathbb{E}[X^e(n) X^e(m)] = \mathbf{1}(nm = \Box). \]
By factoring out the largest square divisor of $n$ and $m$, this for instance implies the second moment estimate
\begin{align*}
    \mathbb{E}\bigg[ \bigg| \frac{1}{\sqrt{x}} \sum_{n \le x} X^e(n) \bigg|^2 \bigg] &= \frac{1}{x} \sum_{n_1, n_2 \le x} \mathbf{1}(n_1 n_2 = \Box) = \frac{1}{x} \sum_{a\leq x} \mu(a)^2 \sum_{n_1, n_2 \le \sqrt{x/a}} 1 
    \asymp \sum_{a\leq x} \frac{\mu(a)^2}{a}\asymp \log x, 
\end{align*} 
(see, e.g., \cite[Theorem 2.2]{MV07}). This suggests that the upper bound in \eqref{eq:armon} is optimal on a wide range of $D$.
Moreover, by H\"older's inequality, this also implies that
\begin{equation}\label{eq:Holdertrivial}
    \mathbb{E}\bigg[ \bigg| \frac{1}{\sqrt{x}} \sum_{n \le x} X^e(n) \bigg|^{2q} \bigg] \ll  (\log x)^q,\quad q\in [0,1].
\end{equation} 
Comparing this to the statement of Theorem \ref{MainThm1_statement} shows that unlike for partial sums of $X$, where the analogous ``trivial'' bound only loses $(\log\log x)^{O(1)}$, the bound in \eqref{eq:Holdertrivial} overshoots by a power of $\log x$ when $q\in [0.01,0.99]$, say. 

\subsection{Main results}

Our primary goal is to establish the following theorem, precisely determining the order of magnitude of the low moments of $\sum_{n\leq x}X^e(n)$.

\begin{thm}\label{MainThm1_statement}
    Let $0 \le q \le 1$ be fixed. Uniformly for $x \ge 3$, we have
    \begin{equation*}
        \mathbb{E} \bigg[ \bigg| \frac{1}{\sqrt{x}}\sum_{n \le x} X^e(n) \bigg|^{2 q} \bigg] 
        \asymp
        \begin{cases}
            (\log \log x)^{-q/2}, &\mbox{ if } 0 \le q < \frac{1}{2}, \\
            (\log \log x)^{3/4}, &\mbox{ if } q = \frac{1}{2}, \\
            (\log x)^{q(2q-1)}, &\mbox{ if } \frac{1}{2} < q \le 1.
        \end{cases}
    \end{equation*}
\end{thm}
In doing so, we also establish the analogous claim for a closely related weighted sum of $X$.
\begin{thm}\label{MainThm2_statement}
    Let $0 \le q \le 1$ be fixed. Uniformly for $x \ge 3$, we have
    \begin{equation*}
        \mathbb{E} \bigg[ \bigg| \sum_{n \le x} \frac{X(n)}{\sqrt{n}} \bigg|^{2 q} \bigg] 
        \asymp
        \begin{cases}
            (\log \log x)^{-q/2}, &\mbox{ if } 0 \le q < \frac{1}{2}, \\
            (\log \log x)^{3/4}, &\mbox{ if } q = \frac{1}{2}, \\
            (\log x)^{q(2q-1)}, &\mbox{ if } \frac{1}{2} < q \le 1.
        \end{cases}
    \end{equation*}
\end{thm}
\noindent Both of these results display a transition at $q=1/2$, where better-than-square-root cancellation gives way to growth proportional to that of a positive power of $\log x$. We also obtain an upper bound with uniformity in $q$. Setting $Y = Y(x,q) = \min\{|1-2q|^{-1},\log\log x\}$ we show that 
\begin{align}\label{q_uniformity}
    \mathbb{E}\bigg[\Big| \frac{1}{\sqrt{x}}\sum_{n \le x}X^{e}(n)\Big|^{2q}\bigg] \ll\left\{
	\begin{array}{ll}
		Y (\log \log x)^{-q/2},  & \mbox{if } 0\leq q \le \tfrac{1}{2},\\
		Y^{3/4} (\log x)^{q(2q-1)}, & \mbox{if } \tfrac{1}{2}\le q \le 1.
	\end{array}
\right.
\end{align}
Moreover, we show that the same statement holds for sums of $X(n)/\sqrt{n}$ (see Equation \eqref{q_uniformity_weighted}).

While a similar transition has been shown to occur for moments of partial sums of $X$ \cite{Harper1, Harper2}, as well as for sums of Steinhaus random multiplicative functions, these transitions both occur at a value $q$ which is greater than or equal to $1$. The smaller threshold in our setting can be anticipated using results from the theory of Gaussian multiplicative chaos measures, more specifically their ``multifractality'' (cf.\ Section \ref{sec:GMCdiscussion}), and suggests that the distributions of our sums have much heavier tails. Indeed, combining (the uniform versions of) Theorems \ref{MainThm1_statement} and \ref{MainThm2_statement} with Markov's inequality yields the following.
\begin{cor}\label{Markov_upper}
    Suppose $x$ is sufficiently large and let $\lambda \ge 3$. Then
    $$\mathbb{P}\Big(\Big|\frac{1}{\sqrt{x}}\sum_{n \le x}X^e(n)\Big| > \frac{\lambda}{(\log\log x)^{1/4}}\Big) \ll \frac{\min\{\log \lambda, \log\log x\}}{\lambda}.$$
    Similarly we have
    \begin{equation*}\begin{split}
        \mathbb{P}&\Big(\Big|\sum_{n \le x}\frac{X(n)}{\sqrt{n}}\Big| > \frac{\lambda}{(\log\log x)^{1/4}}\Big) \ll \frac{\min\{\log \lambda, \log\log x\}}{\lambda}.
    \end{split}
    \end{equation*}
\end{cor}
\begin{proof}
    One applies Markov's inequality after raising both sides to the power $2q$. For $\lambda  \ge \log x$, picking $q=1/2$ and applying the bound in \eqref{q_uniformity} (resp. \eqref{q_uniformity_weighted}) yields the desired result. For smaller $\lambda$, we instead take $q=1/2-1/(2\log\lambda)$.
\end{proof}
\noindent By contrast, the tails of sums of Rademacher and Steinhaus random multiplicative functions were shown by Harper \cite{Harper1} to be $\asymp \lambda^{-2+o(1)}$, in an appropriate range of $\lambda$. It would be of interest to obtain a corresponding lower bound in Corollary \ref{Markov_upper}, which we conjecture should be $\asymp \lambda^{-1+o(1)}$. \par

Noting that \[
\mathbb{E}\Big[\frac{1}{\sqrt{x}}\sum_{n\leq x} X^e(n)\Big] =\frac{1}{\sqrt{x}}\sum_{\substack{n\leq x, \\n=\square}} 1 =1+O\Big(\frac{1}{\sqrt{x}}\Big),\quad \text{ and } \quad \mathbb{E}\Big[\sum_{n\leq x} \frac{X(n)}{\sqrt{n}}\Big]=1,
\]
Markov's inequality also yields the following convergence results upon restricting these partial sums to integers $n$ which are not perfect squares. One sees a clear bias towards $-1$ despite these restricted sums having mean zero, in the following sense.
\begin{cor}\label{randomised_sum_convergence}
    As $x$ tends to infinity, the following limits hold in probability:
    \begin{align*}
        \frac{1}{\sqrt{x}}\sum_{\substack{n \leq x \\ n \neq \square}}X^e(n) \longrightarrow-1,\quad \mathrm{and} \quad \sum_{\substack{1<n \leq x}}\frac{X(n)}{\sqrt{n}} \longrightarrow-1.
    \end{align*}
\end{cor}

Lastly we also prove results relating to a certain weighted sum of $X^e$, motivated by a heuristic connection to moments of $L$--functions $L(1/2,\chi_d)$.
\begin{prop}\label{weighted_extended_Rademacher}
    Let $0 \le q \le \frac{1}{2}$ be fixed. Uniformly for $x \ge 3$, we have
    \begin{equation}\label{extended_weighted_moments}
        \mathbb{E}\bigg[\bigg|\sum_{n \le x}\frac{X^e(n)}{\sqrt{n}}\bigg|^{2q}\bigg] \ll (\log x)^{q(2q+1)},
    \end{equation}
    and as $x$ tends to infinity, 
    \begin{equation}
                \frac{1}{\log x}\sum_{\substack{n \leq x \\ n \ne \square}}\frac{X^e(n)}{\sqrt{n}} {\longrightarrow}-\frac{1}{2}\quad \mathrm{in\,\, probability.}
    \end{equation}
\end{prop}
\begin{proof}
    For the first claim, see Section \ref{sec:proof_weighted_extended_Rademacher} in the appendix. For the second, note that 
    \[
\bigg|\frac{1}{\log x}\sum_{\substack{n \leq x \\ n \ne \square}}\frac{X^e(n)}{\sqrt{n}} +\frac{1}{2}\bigg| = \bigg|\frac{1}{\log x}\sum_{\substack{n \leq x}}\frac{X^e(n)}{\sqrt{n}}- \frac{1}{\log x}\sum_{\substack{n\leq x\\n=\square}}\frac{1}{\sqrt{n}}+\frac{1}{2}\bigg|=\bigg|\frac{1}{\log x}\sum_{\substack{n \leq x}}\frac{X^e(n)}{\sqrt{n}}+O\Big(\frac{1}{\log x}\Big)\bigg|.
    \]
    Convergence in probability then follows by Markov's inequality and the first claim, picking $q=1/4$.
\end{proof}
\noindent Sums of the form $\sum_{n\leq x}\chi_d(n)/\sqrt{n}$ arise naturally when computing moments of $L(1/2,\chi_d)$ via an approximate functional equation, and their averages are expected to behave similarly to the left-hand side in \eqref{extended_weighted_moments} whenever $x$ is small relative to $|d|$. One also finds the same exponent $q(2q+1)$ in \eqref{extended_weighted_moments} as in the Keating-Snaith conjectures \cite{KeatingSnaithL} for averages of $L(1/2,\chi_d)$ (over $|d|\leq x$), where it is a signature of averages over families of $L$-functions with symplectic symmetry.

\subsection{Connections to Gaussian multiplicative chaos with ``insertions''}\label{sec:GMCdiscussion} A now-common strategy to handle sums such as those in Theorem \ref{MainThm1_statement} and \ref{MainThm2_statement} is to relate them to certain random Euler product integrals, following pioneering work of Harper \cite{Harper1}. Focusing on $X^e$ for the purpose of this discussion, one roughly argues that
\begin{equation}\label{eq:discussion_extended}
    \mathbb{E}\Big[\Big|\frac{1}{\sqrt{x}}\sum_{n\leq x}X^e(n)\Big|^{2q}\Big]\approx  \mathbb{E}\Big[\Big(\frac{1}{\log x}\int_{\mathbb{R}}\frac{|F^e(1/2+ih)|^2}{|1/2+ih|^2}\mathrm{d}h\Big)^q\Big], \text{where } F^e(s)=\prod_{\substack{p\leq x\\ p \mathrm{\,prime}} } \Big(1-\frac{X^e(p)}{p^{s}}\Big)^{-1}.
\end{equation}
Were this a sum of Steinhaus random multiplicative functions, one would then recognise $(\log |F^{\mathrm{St}}(1/2+ih)|)_h$ (for $F^{\mathrm{St}}$ the corresponding truncated Euler product) to be well-approximated by a Gaussian process $(g_x(h))_h$ satisfying
\begin{equation}\label{eq:covariance}
        \mathbb{E}[g_x(h)]=0,\quad \mathbb{E}[g_x(h)g_x(h')]=\frac{1}{2}\log \min\{|h-h'|^{-1}, \log x\}+O(1),
\end{equation}
uniformly in $x\geq 3$ and $|h-h'| \le 1$.
Such \textit{logarithmically-correlated} (or \textit{log-correlated}) processes are well-studied in the probability literature, the theory of Gaussian multiplicative chaos (GMC) being specifically concerned with the asymptotic behaviour of integrals of the form $\int f(h) e^{\gamma g_x(h)}\mathrm{d}h$ as $x\to \infty$ (for $\gamma\in \mathbb{R}$ a parameter and $f$ a nice enough test function).

The main difficulty in studying such integrals stems from the covariance in \eqref{eq:covariance} acquiring a logarithmic singularity in this limit. As any process $g$ to which $g_x$ could converge to (in a suitable topology) would formally have covariance 
\begin{equation}\label{eq:logcor}
        \mathbb{E}[g(h)g(h')]=\frac{1}{2}\log\frac{1}{|h-h'|}+O(1),
\end{equation}
it cannot exist pointwise, and therefore cannot, a priori, be exponentiated. GMC theory nevertheless shows that a limiting measure $``e^{\gamma g(h)}\mathrm{d}h"$ can be obtained from $e^{\gamma g_x(h)}\mathrm{d}h$ under a suitable renormalisation procedure, and studies its properties. For $\gamma=2$ in particular, one must multiply $e^{\gamma g_x(h)}\mathrm{d}h$ by $$\frac{\sqrt{\mathbb{E}[g_x(h)^2]}}{\mathbb{E}[e^{2g_x(h)}]}\asymp \frac{\sqrt{\log\log x}}{\log x}$$ to recover a non-trivial limit, and this is known as the \textit{Seneta-Heyde} normalisation of \textit{critical} GMC (see \cite{PowellSurvey} for a survey of this topic).

\bigskip 

 Returning to the integral in \eqref{eq:discussion_extended}, one must make a few preliminary reductions in order to relate it to GMC. For simplicity, we will restrict our attention to the range $h\in [0,1]$ (in practice, one can do so using the decay of $|1/2+ih|^{-2}$, although this is far from immediate). By expanding the logarithm, we can write
 \begin{align*}
         \log|F^e(1/2+ih)|&=\sum_{p\leq x} \sum_{k\geq 1} \frac{X^e(p)^k \cos(kh\log p)}{kp^{k/2}} \\ &= \sum_{p\leq x}\frac{X^e(p)\cos(h \log p)}{p^{1/2}}+\frac{1}{2}\sum_{p\leq x}\frac{\cos(2h \log p)}{p}+O(1)\\
         &:= G_x(h)+\frac{1}{2}\sum_{p\leq x}\frac{\cos(2h \log p)}{p}+O(1)
 \end{align*}
 uniformly in $h$, and standard estimates on sums over primes (cf.\ Lemma \ref{lem:CosineSum}) imply that the second sum is $\log(\min\{\log x, |h|^{-1}\})+O(1)$. It follows that
 \begin{equation}\label{eq:passtoG_xintro}
          \mathbb{E}\Big[\Big(\frac{1}{\log x}\int_0^1\frac{|F^e(1/2+ih)|^2}{|1/2+ih|^2}\mathrm{d}h\Big)^q\Big] \asymp \mathbb{E}\Big[\Big(\frac{1}{\log x}\int_0^1 \min\{|h|^{-1}, \log x\}{e^{2G_x(h)}}\mathrm{d}h\Big)^q\Big],
 \end{equation}
 and one would then want to compare $G_x(h)$ to a Gaussian log-correlated process. It is reasonable to expect Gaussianity here, $G_x(h)$ being a sum of independent random variables for each $h$. However, it is not difficult to show that
 \begin{equation} \label{eq:RademacherCovariance}
          \mathbb{E}\big[G_x(h)G_x(h')\big]=\frac{1}{2}\log {\min\{|h-h'|^{-1},\log x\}}+\frac{1}{2}\log \min\{|h+h'|^{-1}, \log x\}+O(1).
 \end{equation}
Genuine log-correlation can nonetheless be recovered by the following simple observation: for very large $x$, this covariance kernel equals $-\frac{1}{2}\log|h^2-h'^2|+O(1)$ for most $h,h'$, making $G_x$ asymptotically log-correlated under a \textit{different metric}, namely $\mathrm{d}(h,h')=|h^2-h'^2|$. 

This suggests performing the change of variables $u=h^2$ to get
\begin{align}
    \mathbb{E}&\Big[\Big(\frac{1}{\log x}\int_0^1 \min\{|h|^{-1}, \log x\}{e^{2G_x(h)}}\mathrm{d}h\Big)^q\Big]\nonumber  \\
    &\asymp \mathbb{E}\Big[\Big(\frac{1}{\log x}\int_0^1 \min\{|u|^{-1/2}, \log x\}{e^{2G_x(\sqrt{u})}}\frac{\mathrm{d}u}{\sqrt{u}}\Big)^q\Big]\nonumber  \\
    &\asymp \mathbb{E}\Big[\Big(\int_0^1 \min\{|u|^{-1}, (\log x)^2\}\frac{{e^{2G_x(\sqrt{u})}}}{\mathbb{E}[e^{2G_x(\sqrt{u})}]}\frac{\mathrm{d}u}{\sqrt{u}}\Big)^q\Big],\label{eq:rademacher_finalintegral_intro}
\end{align}
and likening this quantity to the one obtained upon substituting $G_x$ by a (centred) Gaussian process $\tilde{g}_x$ with the same covariance structure \eqref{eq:RademacherCovariance}, so that the covariances of $u\mapsto \tilde{g}_x(\sqrt{u})$ now tend to those of $g$ in \eqref{eq:logcor}. Since
\begin{align}
         \mathbb{E}&\Big[\Big(\int_0^1 \min\{|u|^{-1}, (\log x)^2\}\frac{{\sqrt{\mathbb{E}[\tilde{g}_x(\sqrt{u})^2]}e^{2\tilde{g}_x(\sqrt{u})}}}{\mathbb{E}[e^{2\tilde{g}_x(\sqrt{u})}]}\frac{\mathrm{d}u}{\sqrt{u}}\Big)^q\Big] \label{eq:prelimit_rademacher}\\
         &\asymp  (\log\log x)^{q/2}\cdot \mathbb{E}\Big[\Big(\int_0^1 \min\{|u|^{-1}, (\log x)^2\}\frac{e^{2\tilde{g}_x(\sqrt{u})}}{\mathbb{E}[e^{2\tilde{g}_x(\sqrt{u})}]}\frac{\mathrm{d}u}{\sqrt{u}}\Big)^q\Big],\nonumber
\end{align}
($\sqrt{\mathbb{E}[\tilde{g}_x({\sqrt{u}})^2]}$ being the Seneta-Heyde normalisation here) we are led to conjecture that
\[
    \mathbb{E}\Big[\Big|\frac{1}{\sqrt{x}}\sum_{n\leq x}X^e(n)\Big|^{2q}\Big] \asymp (\log\log x)^{-q/2}
\]
whenever \eqref{eq:prelimit_rademacher} tends to a non-trivial limit, or at the very least, whenever 
\[
    \mathbb{E}\Big[\Big(\int_0^1 u^{-3/2} e^{2g(u)}\mathrm{d}u\Big)^q\Big] <\infty.
\]

Quantities of this kind have appeared in the GMC literature under the name of GMC with ``insertion points'' \cite{Remy, RemyZhu} (in this case, an insertion at $0$ of weight $3/2$), borrowing terminology from Liouville conformal field theory where similar integrals arise in correlation functions of so-called vertex operators \cite{DOZZ, DKRV, HRV}. It is known that for $\gamma< 2$, 
\begin{equation}\label{eq:seiberg}
        0<\mathbb{E}\Big[\Big(\int_0^1 u^{a} e^{\gamma g(u)}\mathrm{d}u\Big)^q\Big] <\infty\,\, \iff \,\, a>-\Big(1+\frac{\gamma^2}{4}\Big) \text{ and } q<\min\Big\{\frac{4}{\gamma^2},1+\frac{4(1+a)}{\gamma^2}\Big\},
\end{equation}
with the first condition\footnote{Such conditions sometimes appear under the name of a \textit{Seiberg}, or \textit{extended Seiberg} bound.} being necessary for the measure $e^{\gamma g(u)}\mathrm{d}u$ to be almost surely integrable against $u^a$ (see \cite{HRV}). In the specific case where $g$ is the Gaussian field on $[0,1]$ whose covariance \textit{equals} $-\tfrac{1}{2}\log|u-v|$,\footnote{See \cite{BacryMuzy} for the construction of such a field.} one can go much further: if $\gamma,a$ are such that $q$ is allowed to be a positive integer, these moments can be reduced to the famous Selberg integral \cite{FyodorovBouchaud, RemyZhu}, and expressed as a ratio of Gamma functions. It was then shown by Remy and Zhu \cite{RemyZhu} that this expression persists for fractional $q$, thus determining the exact distribution of these integrals.
 
 The ``if'' direction in \eqref{eq:seiberg} can essentially be seen as a consequence of the \textit{multifractality} of GMC measures, more specifically the fact that for any bounded open set $A\subseteq [0,1]$, 
\begin{equation}\label{eq:multifractality}
        \mathbb{E}\Big[\Big(\int_{rA} e^{\gamma g(u)}\mathrm{d}u\Big)^q\Big]\asymp r^{\xi_\gamma(q)}\,\, \text{as } r\to 0,\quad \text{ where } \quad \xi_\gamma(q)=\big(1+\frac{\gamma^2}{4}\big)q-\frac{\gamma^2}{4}q^2
\end{equation}
for $\gamma\leq 2$ (see \cite[Theorem 2.17]{RhodesVargasSurvey} and \cite[Corollary 7]{DRSVb}). The same argument extends \eqref{eq:seiberg} to $\gamma=2$, noting that \eqref{eq:multifractality} holds for such $\gamma$: subdivide $[0,1]$ into $e$-adic intervals and use the subadditivity of $x\mapsto x^q$ to write
\[
    \mathbb{E}\Big[\Big(\int_0^1 u^{a} e^{2 g(u)}\mathrm{d}u\Big)^q\Big]\ll_{a} \sum_{k} e^{-kaq}\mathbb{E}\Big[\Big(\int_{e^{-(k+1)}}^{e^{-k}}e^{2 g(u)}\mathrm{d}u\Big)^q\Big] \asymp \sum_k e^{-k(aq+\xi_2(q)) }=\sum_k e^{-kq(a+2-q) }.
\]
The resulting exponent is negative when the conditions in \eqref{eq:seiberg} are met. For $a=-3/2$, this translates to $q<1/2$, recovering the transition in Theorem \ref{MainThm1_statement}.

This discussion also suggests that this transition is due to \eqref{eq:passtoG_xintro} becoming dominated by the singularity forming at the origin as $x\to\infty$ (the integral in \eqref{eq:seiberg}, for instance, being finite all the way up to $q<1$ when $a=0$). Indeed, \eqref{eq:passtoG_xintro} is bounded below by
\[
     \mathbb{E}\Big[\Big(\int_0^{1/\log x} {{e^{2G_x(h)}}}{\mathrm{d}h}\Big)^q\Big],
\]
and if we assume that $G_x(h)= G_x(0)+O(1)$ on this range (due to ${G_x(h)}$ varying at a scale of $1/\log x$), we are led to predict that
\[
    \mathbb{E}\Big[\Big|\frac{1}{\sqrt{x}}\sum_{n\leq x}X^e(n)\Big|^{2q}\Big] \gg (\log x)^{-q}\mathbb{E}\big[e^{2q G_x(0)}\big].
\]
As $\mathbb{E}[e^{2qG_x(0)}]\asymp (\log x)^{2q^2}$ by a straightforward calculation, this is consistent with Theorem \ref{MainThm1_statement} for $q>1/2$.

\subsection{Heuristic discussion and proof strategy}

We first sketch the proof of Theorem \ref{MainThm2_statement}, with that of Theorem \ref{MainThm1_statement} essentially reducing to it by the arguments in Section \ref{sec:reduction}.

\subsubsection{Bounds for $q$ away from $1/2$}

For $x \ge 3$, set $t =\log \log x$ and define
\[ F_t(s) := \prod_{\substack{p \le \exp(e^t) \\ p \text{ prime}}} \left( 1 + \frac{X(p)}{p^{1/2+s}} \right) \]
to be the random Euler product associated to $X$, shifted by $\frac{1}{2}$ for notational convenience. (The reason behind the double-exponential indexing in $F_t$ will become apparent in the random walk heuristics used further below.)
We begin by showing that, roughly speaking,
\begin{equation}\label{eq:integral_intro}
    \mathbb{E} \bigg[ \bigg| \sum_{n \le x} \frac{X(n)}{\sqrt{n}} \bigg|^{2 q} \bigg] \approx \mathbb{E}\bigg[ \bigg( \frac{1}{\log x} \int_{\mathbb{R}} \frac{\left| F_t(ih) \right|^2}{\left| {1}/{\log x} + ih \right|^2} \mathrm{d}h \bigg)^q \bigg].
\end{equation} 
This follows the arguments in \cite[Section 2]{G1} fairly closely, which are themselves an adaptation of a similar reduction in \cite{Harper1}. We refer to both of these works for more detail.

Next, we divide the range of integration into $e$-adic intervals, starting at $|h| = \frac{1}{\log x}$, so that the denominator has the same order of magnitude on each interval
\begin{equation}\label{eq:dyadicdecompintro}
        \ll (\log x)^q\mathbb{E}\Big[ \Big( \int_{0}^{1/\log x} \left| F_t(ih) \right|^2 \mathrm{d}h \Big)^q \Big]+\sum_{\substack{H\geq 1/\log x\\ H\text{ $e$-adic}}} \frac{H^{-2q}}{(\log x)^q}\mathbb{E}\Big[ \Big( \int_{H}^{eH} \left| F_t(ih) \right|^2 \mathrm{d}h \Big)^q \Big],
\end{equation}
which essentially reduces the problem to bounding expressions of the shape
\[ \mathbb{E}\Big[ \Big( \int_{H}^{eH} \left| F_t(ih) \right|^2 \mathrm{d}h \Big)^q \Big] \]
for various $H \ge \frac{1}{\log x}$. The range $H \ge 1$ can be dealt with directly using bounds of Harper \cite{Harper4} in combination with an elementary bound for very large $H$. For the ``mesoscopic range'' $\frac{1}{\log x} \le H \le 1$, as well as the first term in \eqref{eq:dyadicdecompintro}, we make use of the fact that since $F_{\log 1/H}$ typically varies at scale $H$, $|F_{\log 1/H}(h)|^2$ does not change substantially over $h\in [0,eH]$. One then morally has
\begin{equation}\label{eq:meso_quantity}
    \mathbb{E}\Big[ \Big( \int_{H}^{eH} \left| F_t(ih) \right|^2 \mathrm{d}h \Big)^q \Big] \approx \mathbb{E}[|F_{\log 1/H}(0)|^{2q}] \mathbb{E}\Big[ \Big( \int_{H}^{eH} \left| F_{\log 1/H,t}(ih) \right|^2 \mathrm{d}h \Big)^q \Big],
\end{equation}  
where $F_{\log 1/H, t}$ denotes the Euler product ranging over primes $p\in [e^{1/H}, \exp(e^t)]$, and we make such a statement rigorous using Hölder's inequality.
It is then not difficult to show that
\begin{equation}\label{eq:F_0_intro}
    \mathbb{E}[|F_{\log 1/H}(0)|^{2q}] \asymp H^{-q(2q-1)},
\end{equation} 
(cf.\ Lemma \ref{lem:MomentsEulerProducts}).
For the second factor in \eqref{eq:meso_quantity}, we require a mesoscopic analogue of the better-than-squareroot cancellation result in \cite{Harper1}, showing in this case that
\begin{equation}\label{eq:meso_harper_intro}
    \mathbb{E}\Big[ \Big( \frac{1}{H} \int_{H}^{eH} \left| F_{\log 1/H,t}(ih) \right|^2 \mathrm{d}h \Big)^q \Big] \ll \left( \frac{H \log x}{\sqrt{\log (H \log x)}} \right)^q.
\end{equation}  
Combining \eqref{eq:F_0_intro}, \eqref{eq:meso_harper_intro}, and summing over $\frac{1}{\log x} \le H \le 1$ in \eqref{eq:dyadicdecompintro} then delivers the upper bound. 

As one would expect from the heuristics of the previous section, the main contribution for fixed $q>1/2$ turns out to come from the smallest values of $H$, namely those of size $\asymp \frac{1}{\log x}$. For $q<1/2$, the dominant contribution is instead from $H\asymp 1$. This makes the task of deriving a corresponding lower bound particularly simple, as one can simply restrict the integral in \eqref{eq:integral_intro} to such $H$ by positivity of the integrand. The corresponding lower bound in \eqref{eq:meso_harper_intro} for $H\asymp 1/\log x$ is then straightforward to show, and it is already known for $H\asymp 1$ \cite{Harper1}.

\subsubsection{The case $q=1/2$, and a branching random walk model}

In the transition as $q$ approaches $\frac{1}{2}$ sufficiently fast, \textit{all} values of $\frac{1}{\log x}\leq H\leq 1$ begin to contribute significantly to the final bound: there are $\asymp \log \log x$ such terms, each typically contributing $\asymp (\log (H \log x))^{-1/4}$, leading to the upper bound of  $\ll(\log\log x)^{3/4}$ in Theorem \ref{MainThm2_statement}.

This makes a corresponding lower bound for $q\approx 1/2$ significantly more challenging to obtain, since restricting to a single ({deterministic}) value of $H$ would now necessarily be lossy. However, it will turn out to be sufficient to use the approximation 
\begin{align*}
    \mathbb{E}\bigg[ \bigg( \frac{1}{\log x} \int_{\frac{1}{\log x}}^1 \frac{\left| F_t(ih) \right|^2}{\left| {1}/{\log x} + ih \right|^2} \mathrm{d}h \bigg)^{1/2} \bigg] &\asymp \mathbb{E}\bigg[ \bigg( \frac{1}{\log x} \sum_{\substack{\frac{1}{\log x} \le H \le 1 \\ H \, e\text{-adic}}} \frac{1}{H^2} \int_H^{eH} \left| F_t(ih) \right|^2 \mathrm{d}h \bigg)^{1/2} \bigg] \\
    &\approx \mathbb{E}\bigg[ \bigg( \frac{1}{\log x} \max_{\substack{\frac{1}{\log x} \le H \le 1 \\ H \, e\text{-adic}}} \frac{1}{H^2} \int_H^{eH} \left| F_t(ih) \right|^2 \mathrm{d}h \bigg)^{1/2} \bigg],
\end{align*}
stressing that this maximum is attained at a \textit{random} value of $H$. To study this quantity, we first factor out the contribution of the Euler product over $p \le e^{1/H}$ as was done earlier, which we recall should be well-approximated by $|F_{\log 1/H}(0)|^2$. This yields
\[ \approx \mathbb{E}\bigg[ \max_H {|F_{\log 1/H}(0)|} \Big( \frac{1}{H^2 \log x}\int_H^{eH} \left| F_{\log 1/H,t}(ih) \right|^2 \mathrm{d}h \Big)^{1/2} \bigg], \]
the maximum still ranging over the same set of $H$. For brevity, let us further denote $t_H=\log1/H$, and $W_{t_H} := (H^{2}\log x)^{-1}\int_H^{eH} | F_{t_H,t}(ih) |^2 \mathrm{d}h$ so that this is simply
\begin{equation}\label{eq:main_intro_lower}
    = \mathbb{E}\bigg[ \max_H|F_{t_H}(0)|\cdot W_{t_H}^{1/2} \bigg], 
\end{equation} 
noting that  $F_{t_H}(0)$ and $W_{t_H}$ are independent for every $H$.

\bigskip 

Prior to discussing how we bound this quantity, let us introduce a branching random walk model for $(\log |F_t(ih)|)_h$ which will be helpful to keep in mind. We follow \cite{ArguinBeliusHarper}, where such a model was first introduced in the setting of Steinhaus random multiplicative functions. Begin by noting that for any $h$,
\[
    \log|F_{t}(ih)|=\sum_{p\leq \exp(e^t)}\bigg(\frac{X(p)\cos(h\log p)}{p^{1/2}}-\frac{\cos(2h\log p)}{2p}\bigg)+O(1) 
\]
by expanding the logarithm. Ignoring the $O(1)$ error term, we want to view these sums as a random walk, and do so by subdividing them further into ``increments''
\[
    Y_k(h)=\sum_{\exp(e^{k-1})<p\leq \exp(e^k)}\bigg(\frac{X(p)\cos(h\log p)}{p^{1/2}}-\frac{\cos(2h\log p)}{2p}\bigg),\quad\quad 0\leq k \leq t.
\]
This ensures that each $Y_k(h)$ is well-approximated by a Gaussian random variable of the same mean and variance, being a sum of (many) independent random variables, and that these means and variances are comparable for different $k$. For instance, at $h=0$, Mertens' theorem implies that $\mathbb{E}[Y_k(0)]\approx -1/2$ and $\text{Var}{(Y_k(0))}\approx 1$ up to a small error term, uniformly over $k$. One can then (heuristically) view $\log|F_t(0)|$ as a Gaussian random walk with $t$ increments, each with distribution $\mathcal{N}(-1/2,1)$. Furthermore, it is not difficult to show that for any $h,h'$, very roughly,
\[
\text{Cov}\big[Y_k(h),Y_k(h')\big]\approx
\left\{
	\begin{array}{ll}
		\text{Var}(Y_k(h))  & \mbox{if } k\leq -\log|h-h'| \\
		0 & \mbox{if } k> -\log|h-h'|,
	\end{array}
\right.
\]
meaning that the sums $\sum_{k\leq t}Y_k(h)$ and $\sum_{k\leq t}Y_k(h')$ have almost the same increments up to $k=-\log|h-h'|$, after which they decorrelate very quickly.

These observations suggest modelling the process $(\log|F_t(ih)|)_{h\in [0,1]}$ by a \textit{branching random walk}, which is a tree-indexed Gaussian process where such a correlation structure is exact. To be precise, consider a full, rooted binary tree of depth $t$, with leaves indexed by $h_j=j/2^t$ for $0\leq j< 2^t$. To each edge $e$, attach an independent Gaussian random variable $N_e$ as follows: if $e$ has depth\footnote{On a path from the root to a leaf, the first edge has depth one, the second depth two, etc.} $k$, let $N_e$ have the same mean and variance as $Y_k(h_{j(e)})$, where $$j(e)=\inf\{j:\text{there is a path from the root to $h_j$ using $e$}\}.$$
Lastly, for each leaf $h_j$ and integer $1\leq k\leq t$, let $S_k(h_j)=\sum_{\ell\leq k} \mathcal{Y}_\ell(h_j)$, where $\mathcal{Y}_\ell(h_j)$ denotes the Gaussian random variable associated to the depth-$\ell$ edge on the path from the root to $h_j$. Then $(S_k(h_j))_k$ is a random walk for each $h_j$, and for most pairs $j,j'$, $S_k(h_j)$ and $S_k(h_{j'})$ share their first $[ -\log_2|h_j-h_{j'}|]$  increments (with $\log_2$ here denoting the base-2 logarithm), after which they evolve independently. The process $(S_t(h_j))_{j}$ can thus be seen as a toy model for $(\log|F_t(ih_j)|)_{j}$, which will be used as a guide in what follows.

\bigskip

\begin{figure}
    \centering
    \input{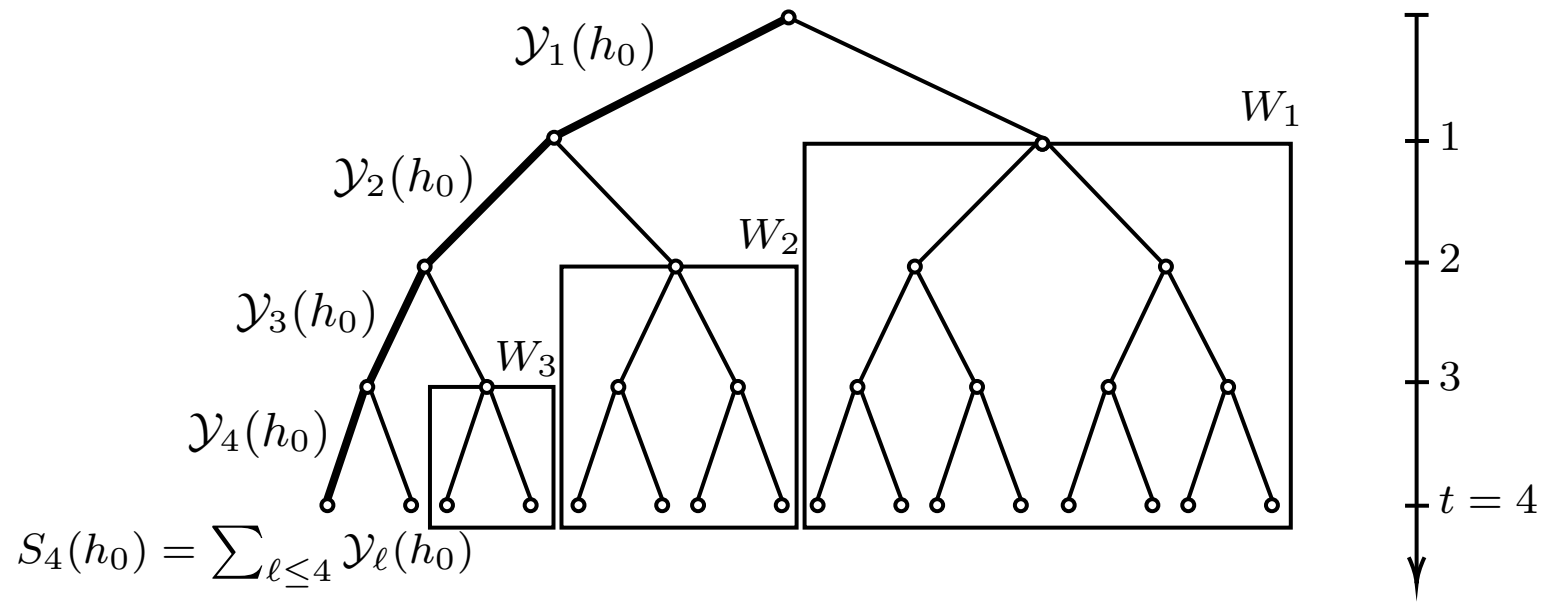}
    \caption{{A branching random walk toy model for $\log|F_t(ih)|$, $t=4$. The random walk $(S_k(h_0))_{k\leq t}$ corresponding to $(\log|F_{k}(0)|)_{k\leq t}$ is given by the leftmost branch, depicted in bold, with the leftmost leaf having value $S_t(h_0)$. For each $k$, $W_{k}$ is modelled by summing $e^{2(S_t(h_j)-S_k(h_j))}$ over the leaves $h_j$ of the corresponding framed subtree.}}
    \label{fig:BRW}
\end{figure}

Let us now return to our main task of estimating \eqref{eq:main_intro_lower}, in order to prove the lower bound in Theorem \ref{MainThm2_statement} at $q=1/2$. Our strategy begins by writing 
\[
         \mathbb{E}\bigg[ \max_H|F_{t_H}(0)|\cdot W_{t_H}^{1/2} \bigg] \geq \mathbb{E}\bigg[ |F_{ t_{H^*}}(0)|\cdot W_{{t_{H^*}}}^{1/2} \bigg],
\]
where $H^*$ is the value of $H$ at which the maximum \textit{of the Euler product} $|F_{t_{H^*}}(0)|$ is attained, then decomposing this quantity according to the value of $H^*$:
\begin{align}\label{eq:decompH}
        \sum_{H}\mathbb{E}\bigg[ |F_{t_{H^*}}(0)|\cdot W_{t_{H^*}}^{1/2}\cdot \mathbf{1}(H^*=H)\bigg].
\end{align}
The branching random walk model then suggests viewing $|F_{t_H}(0)|$ as the exponential of a random walk $(S_k(h_0))_{k\leq t_H}=(S_k)_{k\leq t_H}$ with $\approx \mathcal{N}(-1/2,1)$-distributed increments. Note also that in this model, the random variables $(W_{t_H})_{H}$ are not only independent from $S_{t_H}$ for each $H$, but also from the $\sigma$-algebra generated by the entire random walk path $(S_k)_{k\leq t}$ (cf.\ Figure \ref{fig:BRW}), and thus independent from the event $\{H^*=H\}$. One should then roughly have
\begin{align*}
        \sum_H \mathbb{E}\big[ |F_{t_{H^*}}(0)|\cdot W_{t_{H^*}}^{1/2}\cdot \mathbf{1}(H^*=H)\big] \approx \sum_H \mathbb{E}\big[ e^{S_{t_H}} \mathbf{1}(H^*=H)\big] \mathbb{E}\big[W_{t_H}^{1/2}\big],
\end{align*}
which by the estimate in \eqref{eq:meso_harper_intro}, is
\[
\asymp \sum_H \big(\log(H\log x)\big)^{-1/4}\mathbb{E}\big[ e^{S_{t_H}} \mathbf{1}(H^*=H)\big],
\]
thus reducing the task to a problem on Gaussian random walks.

Letting $\mathbb{Q}$ be the measure defined through
\[
        \frac{\mathrm{d}\mathbb{Q}}{\mathrm{d}\mathbb{P}}=\frac{e^{S_{t_H}}}{\mathbb{E}e^{S_{t_H}}}=e^{S_{t_H}}
\]
we can write
\[
\mathbb{E}\big[ e^{S_{t_H}} \mathbf{1}(H^*=H)\big] = \mathbb{Q}(H=H^*),
\]
and decompose $\{H=H^*\}$ into two independent events, each pertaining to its own random walk. Indeed, it is clear that
\begin{align*}
        \big\{H=H^*\big\}= \big\{S_{t_H}-S_{t_H-k}\geq 0\,\, \forall 1\leq k<t_H\big\}\cap \big\{S_{t_H+k}-S_{t_H}\leq 0\,\, \forall 1\leq k\leq t-t_H\big\},
\end{align*}
and not difficult to show that under $\mathbb{Q}$, $(S_{t_H}-S_{t_H-k})_{k<t_H}$ and $(S_k-S_{t_H})_{k>t_H}$ remain independent, Gaussian random walks, only with the former now having a (per-increment) drift of $+1/2$. The probability that a walk with positive (resp. negative) drift remains positive (resp. negative) is bounded below by a constant, uniformly in time, and we therefore conclude the desired bound:
\[
    \sum_H \big(\log(H\log x)\big)^{-1/4}\mathbb{E}\big[ e^{S_{t_H}} \mathbf{1}(H^*=H)\big]\gg  \sum_H \big(\log(H\log x)\big)^{-1/4} \asymp (\log\log x)^{3/4}.
\]
\bigskip

In order to make these heuristics rigorous, we proceed slightly differently and decompose according to both the value of $H^*$ (as in \eqref{eq:decompH}) as well as the approximate location of $\log|F_{k}(0)|$, using events of the shape
\[
    \mathcal{E}(k^*,V)=\big\{ \max_{k\leq t}\log |F_k(0)|=\log|F_{k^*}(0)|\in[V,V+1)\big\}.
\]
(The true definition of $\mathcal{E}$ differs slightly from this for technical reasons, cf.\ Section \ref{sec:lower}, Equation \eqref{eq:E(k,V)def}.) Noting that $t_{H^*}=k^*$ and $|F_{t_{H^*}}(0)|\asymp e^V$ on such an event, this allows us to bound $\mathbb{E}[|F_{t_{H^*}}(0)|\cdot W_{t_{H^*}}^{1/2}]$ below by
\begin{equation}\label{eq:lbsketch}
    \sum_{k^*}\sum_V\mathbb{E}\Big[|F_{t_{H^*}}(0)| W_{t_{H^*}}^{1/2}\cdot\mathbf{1}
    \big(\mathcal{E}(k^*,V)\big)\Big] \asymp  \sum_{k^*}\sum_V e^{V} \mathbb{E}\Big[W_{k^*}^{1/2}\cdot\mathbf{1}\big(\mathcal{E}(k^*,V)\big)\Big],
\end{equation}
with the sums running over $1\leq k^*\leq t$ and $V\in \mathbb{Z}$,  respectively. In practice, we will actually restrict to $C<V\leq t/C$ for a large constant $C>0$, and $|k^*-2V|\leq \sqrt{V}$, noting that this is where one expects the dominant contribution to come from according to the branching random walk model; we are free to do so for a lower bound, and this simplifies the proofs in Section \ref{sec:E_results}.

As $(|F_k(0)|)_k$ and $(W_k)_k$ do not enjoy the genuine independence present between their counterparts in said model, we instead rewrite  each summand in terms of a conditional expectation
\[
    e^V\cdot \mathbb{E}\Big[W_{k^*}^{1/2}\big|\mathcal{E}(k^*,V)\Big] \mathbb{P}\big(\mathcal{E}(k^*,V)\big),
\]
and argue that this conditioning does not affect the order of magnitude, in the sense that
\[
     \mathbb{E}\Big[W_{k^*}^{1/2}\big|\mathcal{E}(k^*,V)\Big] \gg \big(\log (e^{-k^*}\log x)\big)^{-1/4}.
\]
(see Section \ref{sec:lower}).
This relies on rather elaborate tilted measure arguments, the details of which can be found in Section \ref{sec:E_results}. The same section also shows that
\[
    \mathbb{P}\big(\mathcal{E}(k^*,V)\big)\asymp \frac{1}{\sqrt{V}} e^{-V},
\]
(cf.\ Proposition \ref{ProbLargeMax}) by comparing $(\log |F_k(0)|)_k$ to a Gaussian random walk, and inserting both of these estimates into \eqref{eq:lbsketch} yields the desired claim upon summing over $k^*$ and $V$.

\subsubsection{Proof of Theorem \ref{MainThm1_statement}}\label{sec:reduction} The  proof of Theorem \ref{MainThm1_statement} essentially reduces to that of Theorem \ref{MainThm2_statement} by the following observations. Any positive integer $n$ can be written as $n = s^2 m$, where $s^2$ is the largest square that divides $n$. Since $X^e(s^2) = X^e(s)^2 = 1$, this implies that 
\[ \sum_{n \le x} X^e(n) = \sum_{m \le x} X(m) \sum_{s^2 \le x/m} 1 \approx \sqrt{x} \sum_{m \le x} \frac{X(m)}{\sqrt{m}}. \]
There are some subtleties involved in making this rigorous, with our argument in particular requiring the estimate
\[
    \mathbb{E} \Bigg[ \bigg| \frac{1}{\sqrt{x}} \sum_{n \le x} X(n) \left\{ \sqrt{ \frac{x}{n} } \right\} \bigg|^{2q} \Bigg]\ll (\log\log x)^{-q/2},\quad 0\leq q\leq 1/2,
\]
with $\{\cdot\}$ denoting the fractional part  (see Proposition \ref{BoundFractPart}). We prove this by decomposing this into blocks where
\[
    \frac{x}{(m+1)^2}<n\leq \frac{x}{m^2},
\]
so as to fix the value of the integer part $[\sqrt{x/n}]$ of $\sqrt{x/n}$, and bounding each block's contribution by a result of Caich \cite{Caich1} on sums of $X$ in short intervals (see Section \ref{sec:ExtendedWeightedReduction} for more detail).

This reduction is used to get the most of the bounds in Theorem \ref{MainThm1_statement}, with the exception of the lower bound at $q<1/2$ where a different argument is required. In that case, we essentially compare the $2q$-th moment of $\sum_{n} X^e(n)$ to a quantity similar to the right-hand side in \eqref{eq:integral_intro}. The same argument used to bound this quantity then applies.
\subsection{Organisation} Section \ref{sec:prelim} records standard estimates on moments of Euler products and on prime sums. Section \ref{sec:ExtendedWeightedReduction} shows how to relate partial sums of $X^e(n)$ to those of $X(n)/\sqrt{n}$, and 
Section \ref{sec:pass_to_euler} relates the latter to  integrals involving random Euler products. In Sections \ref{sec:upper} and \ref{sec:lower} we prove the bounds obtained in Theorem \ref{MainThm2_statement}, with the auxiliary results needed for the case when $q=1/2$ being proved in Section \ref{sec:E_results} and the appendix.

\section*{Acknowledgements}

A large part of this research was conducted at the Centre de Recherches Mathématiques, as part of the ``Universal Statistics in Number Theory'' thematic program, and we would like to thank the organisers for support. We thank Nathan Creighton, Adam Harper, Christian Webb, and Max Xu for helpful conversations that contributed to the work, as well as Seth Hardy and Paweł Nosal for their feedback. C.A.\ thanks the Heilbronn Institute for Mathematical Research for support. J.H.\ is supported by the Engineering and Physical Sciences Research Council (EPSRC) Grant EP/Z535990/1. Claude was used as an auxiliary editorial tool for proofreading, and to generate the formalisation in Lean 4 (cf. Section \ref{lean}). 

\section{Preliminary Estimates}\label{sec:prelim}
We record a few estimates which will be used repeatedly throughout the paper.
\subsection{Moment bounds for random Euler products}

\begin{lem}[Euler Product Result 2, \cite{Harper2}]\label{lem:MomentsEulerProducts}
     Suppose $400 \leq y < z$ for sufficiently large $z$ with $\alpha,\beta,h_1,h_2 \in \mathbb{R}$ where $\alpha,\beta$ are fixed and $\sigma \geqslant -10/\log z$, then 
    \begin{equation*}\begin{split}
        \mathbb{E} & \bigg[\prod_{y < p \leq z}\Big|1+\frac{X(p)}{p^{1/2+\sigma+ih_1}}\Big|^{2\alpha}\Big|1+\frac{X(p)}{p^{1/2+\sigma+ih_2}}\Big|^{2\beta}\bigg] \\
        = \exp\bigg(&\sum_{y < p \leq z}\frac{\alpha^2 + \beta^2 + (\alpha^2 - \alpha)\cos(2h_1\log p) + (\beta^2-\beta)\cos(2h_2\log p)}{p^{1+2\sigma}} \\
        &+ \sum_{y < p \leq z} \frac{2\alpha\beta(\cos((h_1-h_2)\log p) + \cos((h_1+h_2)\log p)}{p^{1+2\sigma}} + O\Big(\frac{1}{\sqrt{y}\log y}\Big)\bigg). \end{split}
    \end{equation*} 
\end{lem}

\subsection{Prime number theorem estimates} The following estimates are a straightforward consequence of the prime number theorem and partial summation. 

\begin{lem}\label{lem:PNT}
There exists a $c>0$ such that, uniformly in $b>a>2$ and $s \ge 1/2-10/\log b$, we have
    \begin{equation}
        \sum_{a<p\leq b} \frac{1}{p^{2s}}= \int_a^b \frac{\mathrm{d}t}{t^{2s}\log t}+O\big(e^{-c\sqrt{\log a}}\big).
    \end{equation}
In particular, 
    \begin{equation*}
 \sum_{a<p\leq b} \frac{1}{p^{2s}} = 
\begin{cases}
    \log\log b-\log\log a+O\big(e^{-c\sqrt{\log a}}\big) & \text{if } s = 1/2 \\
    \mathrm{Ei}\big((1-2s)\log b\big)-\mathrm{Ei}\big((1-2s)\log a\big)+O\big(e^{-c\sqrt{\log a}}\big)  & \text{if } s \ne 1/2.
\end{cases}
\end{equation*}
\end{lem}

\begin{lem}\label{lem:log_PNT}
    Uniformly in $1\leq a\leq b$, and $m\geq 1$,
    \[
\sum_{a\leq p\leq b}\frac{(\log p)^m}{p}=O\big((\log b)^m\big).
    \]
\end{lem}

Using standard facts about the Riemann zeta function $\zeta$, we also derive the following.

\begin{lem}\label{lem:CosineSum}
    Let $y \ge 3$ and $1/\log y < \sigma \le 1$. Uniformly for $|h| \le \frac{1}{2}$, we have
    \begin{equation}\begin{split}\label{eq:CosineSum1}
        \sum_{p \le y} \frac{\cos(2h \log p)}{p^{1+2\sigma}} = \log |\sigma+ih|^{-1} + O(1).
    \end{split}\end{equation}
    Uniformly for $|h| \ge \frac{1}{2}$, we have
    \begin{equation}\label{eq:CosineSum2}
        \sum_{p \le y} \frac{\cos(2h \log p)}{p^{1+2\sigma}} \le \log \log (|h|+3) + O(1).
    \end{equation}
\end{lem}

\begin{proof}
    First, note that
    \[ \left| \sum_{p > y} \frac{\cos(2h\log p)}{p^{1+2\sigma}} \right| \le \sum_{p > y} \frac{1}{p^{1+2\sigma}} = O(1) \]
    uniformly for $\sigma \ge 1/\log y$. Hence
    \[ \sum_{p \le y} \frac{\cos(2h \log p)}{p^{1+2\sigma}} = \sum_{p } \frac{\cos(2h \log p)}{p^{1+2\sigma}} + O(1) = \log |\zeta(1+2 \sigma+2ih)| + O(1). \]
    Both claims then follow from the Laurent expansion of $\zeta$ around $1$ in combination with elementary upper bounds for $\zeta$ in vertical strips (cf. \cite[Corollary 1.17]{MV07}).
\end{proof}

\section{From Extended Rademacher Sums to Weighted Rademacher Sums}\label{sec:ExtendedWeightedReduction}

We begin by reducing the upper bound in Theorem \ref{MainThm1_statement} to the upper bounds in Theorem \ref{MainThm2_statement}, as well as the lower bound on the range $\frac{1}{2} \le q \le 1$. For the lower bound on the remaining range, a different strategy will be required (see Proposition \ref{lower_extended_reduction}). 

This reduction is a direct consequence of the following two propositions.

\begin{prop}\label{ExtToWeighted}
    Let $x \ge 3$ and $0 \le q \le 1$. Then we have
    \begin{equation}
        \mathbb{E} \bigg[ \bigg| \frac{1}{\sqrt{x}}\sum_{n \le x} X^e(n) \bigg|^{2 q} \bigg] \le 2 \mathbb{E} \bigg[ \bigg| \sum_{n \le x} \frac{X(n)}{\sqrt{n}} \bigg|^{2 q} \bigg] + 2 \mathbb{E} \bigg[ \bigg| \frac{1}{\sqrt{x}} \sum_{n \le x} X(n) \left\{ \sqrt{ \frac{x}{n} } \right\} \bigg|^{2q} \bigg]
    \end{equation}
    and
    \begin{equation}
        \mathbb{E} \bigg[ \bigg| \frac{1}{\sqrt{x}}\sum_{n \le x} X^e(n) \bigg|^{2 q} \bigg] \ge \frac{1}{2} \mathbb{E} \bigg[ \bigg| \sum_{n \le x} \frac{X(n)}{\sqrt{n}} \bigg|^{2 q} \bigg] - \mathbb{E} \bigg[ \bigg| \frac{1}{\sqrt{x}} \sum_{n \le x} X(n) \left\{ \sqrt{ \frac{x}{n} } \right\} \bigg|^{2q} \bigg].
    \end{equation}
\end{prop}

\begin{proof}
    We group integers $n$ by dividing out their largest square divisor, writing $n = s^2 m$ with $m$ squarefree. Since $X^e(n) = X^e(s)^2 X^e(m) = X(m)$, this gives
    \begin{align*}
        \mathbb{E} \bigg[ \bigg| \frac{1}{\sqrt{x}}\sum_{n \le x} X^e(n) \bigg|^{2 q} \bigg] &= \mathbb{E} \bigg[ \bigg| \frac{1}{\sqrt{x}}\sum_{m \le x} X(m) \sum_{s \le \sqrt{x/m}} 1 \bigg|^{2 q} \bigg] \\
        &= \mathbb{E} \bigg[ \bigg| \frac{1}{\sqrt{x}}\sum_{m \le x} X(m) \left( \sqrt{ \frac{x}{m}} - \left \{ \sqrt{ \frac{x}{m}} \right \} \right)\bigg|^{2 q} \bigg].
    \end{align*}
    The claim then follows from the elementary inequalities
    \[ \frac{1}{2}|a|^{2q} - |b|^{2q} \le |a+b|^{2q} \le 2(|a|^{2q} +|b|^{2q}), \]
    valid on the range $0 \le q \le 1$.
\end{proof}

\begin{prop}\label{BoundFractPart}
    Uniformly in $x \ge 3$ and $0 \le q \le 1$, we have
    \begin{equation}
        \mathbb{E} \bigg[ \bigg| \frac{1}{\sqrt{x}} \sum_{n \le x} X(n) \left\{ \sqrt{ \frac{x}{n} } \right\} \bigg|^{2q} \bigg] \ll
        \begin{cases}
            (\log \log x)^{-q/2}, &\mbox{for } 0 \le q \le \frac{1}{2}, \\
            1, &\mbox{for } \frac{1}{2} \le q \le 1.
        \end{cases}
    \end{equation}
\end{prop}

We note that one could likely improve this upper bound on the range $q > \frac{1}{2}$, but that is not needed for the final claims.

\begin{proof}
    The second bound follows from a simple mean square calculation by noting that
    \begin{align*} 
    \mathbb{E} \bigg[ \bigg| \frac{1}{\sqrt{x}} \sum_{n \le x} X(n) \left\{ \sqrt{ \frac{x}{n} } \right\} \bigg|^{2q} \bigg] &\le 
    \mathbb{E} \bigg[ \bigg| \frac{1}{\sqrt{x}} \sum_{n \le x} X(n) \left\{ \sqrt{ \frac{x}{n} } \right\} \bigg|^{2} \bigg]^q \\
    &= \bigg( \frac{1}{x} \sum_{\substack{n \le x \\ n \text{ squarefree}}} \left\{ \sqrt{ \frac{x}{n} } \right\}^2 \bigg)^q \le 1.
    \end{align*}
    Similarly, through an application of H\"older's inequality, we can reduce the remaining range to the case $q = \frac{1}{2}$. We split the sum into intervals according to the value of the integer part by looking at
    \[ \frac{x}{(m+1)^2} < n \le \frac{x}{m^2} \]
    for $m \in \mathbb{N}$, where
    $\left[ \sqrt{ \frac{x}{n} } \right] = m$. Through triangle inequality, this implies that
    \begin{align*}
        \mathbb{E} \bigg[ \bigg| \sum_{n \le x} X(n) \left\{ \sqrt{ \frac{x}{n} } \right\} \bigg| \bigg] &\le \mathbb{E} \bigg[ \bigg| \sum_{n \le \frac{x}{(M+1)^2}} X(n) \left\{ \sqrt{ \frac{x}{n} } \right\} \bigg| \bigg] + \sum_{m = 1}^M \mathbb{E} \bigg[ \bigg| \sum_{\frac{x}{(m+1)^2} <  n \le \frac{x}{m^2}} X(n) \left( \sqrt{ \frac{x}{n} } - m \right) \bigg| \bigg] 
    \end{align*}
    for some parameter $M$. For the first part, an application of the Cauchy-Schwarz inequality yields that
    \[ \mathbb{E} \bigg[ \bigg| \sum_{n \le \frac{x}{(M+1)^2}} X(n) \left\{ \sqrt{ \frac{x}{n} } \right\} \bigg| \bigg] \ll \frac{\sqrt{x}}{M+1}, \]
    which is sufficient with the choice $M = \left[ (\log \log x)^{1/4}\right]$.

    We will now proceed to show that 
    \begin{equation}\label{FracPartShortInt}
        \mathbb{E} \bigg[ \bigg| \sum_{\frac{x}{(m+1)^2} <  n \le \frac{x}{m^2}} X(n) \left( \sqrt{ \frac{x}{n} } - m \right) \bigg| \bigg] \ll \frac{\sqrt{x}}{(\log \log x)^{1/4}} \frac{\sqrt{\log (m+1)}}{m^{3/2}}
    \end{equation}
    uniformly for $1 \le m \le M$. Since the dependence on $m$ in this bound is summable, this implies the claim.

    To this end, set
    \[ g(t) = g(t; x, m) := \sqrt{ \frac{x}{t} } - m  \]
    and note that
    \[ |g'(t)| = \frac{\sqrt{x}}{2 t^{3/2}} \asymp \frac{m^3}{x} \]
    for $\frac{x}{(m+1)^2} < t \le \frac{x}{m^2}$, which is itself an interval of length $\asymp \frac{x}{m^3}$. Moreover, we have
    \[ g \left( \frac{x}{(m+1)^2} \right) = 1 \qquad \text{and} \qquad g \left( \frac{x}{m^2} \right) = 0.  \]
    Setting
    \[ S(x,x+y) := \sum_{x < n \le x+y} X(n), \]
    summation by parts thus gives that
    \begin{align*}
        \sum_{\frac{x}{(m+1)^2} <  n \le \frac{x}{m^2}} X(n) \left( \sqrt{ \frac{x}{n} } - m \right) 
        &= \sum_{\frac{x}{(m+1)^2} <  n \le \frac{x}{m^2}} X(n) g(n) 
        = - \int_{\frac{x}{(m+1)^2}}^{\frac{x}{m^2}} S \left( \frac{x}{(m+1)^2}, t \right) g'(t) \mathrm{d}t,
    \end{align*}
    hence
    \[ \mathbb{E} \bigg[ \bigg| \sum_{\frac{x}{(m+1)^2} <  n \le \frac{x}{m^2}} X(n) \left( \sqrt{ \frac{x}{n} } - m \right) \bigg| \bigg] \ll \frac{m^3}{x} \int_{\frac{x}{(m+1)^2}}^{\frac{x}{m^2}} \mathbb{E} \bigg[ \left| S \left( \frac{x}{(m+1)^2}, t \right) \right| \bigg] \mathrm{d}t. \]
    We can bound this inner expectation by applying work of Caich \cite{Caich1} for low moments of random multiplicative functions on short intervals. His results give that
    \begin{equation}
        \mathbb{E} \bigg[ \left| S \left( \frac{x}{(m+1)^2}, t \right) \right| \bigg] \ll \left( \left( t - \frac{x}{(m+1)^2}+1 \right) \min \left \{ 1, \frac{\log \left( 2 + \frac{x}{t(m+1)^2 - x} \right) }{\sqrt{\log \log x}} \right\} \right)^{1/2}
    \end{equation}
    after noting that
    \[ \log \log \frac{x}{(m+1)^2} \asymp \log \log x \]
    uniformly on the range $1 \le m \le M = \left[ (\log \log x)^{1/4}\right]$. Substituting $s = t - \frac{x}{(m+1)^2}$, we deduce
    \begin{align*} &\mathbb{E} \bigg[ \bigg| \sum_{\frac{x}{(m+1)^2} <  n \le \frac{x}{m^2}} X(n) \left( \sqrt{ \frac{x}{n} } - m \right) \bigg| \bigg] 
    \\&\ll \frac{m^3}{x} \int_0^{\frac{x}{m^2} - \frac{x}{(m+1)^2}} \sqrt{s+1} \min \left \{ 1, \frac{\log \left( 2 + \frac{x}{(m+1)^2 s} \right) }{\sqrt{\log \log x}} \right\}^{1/2} \mathrm{d}s \\
    &\leq \frac{m^3}{x} \int_0^{\frac{x}{(m+1)^2 \exp (\sqrt{\log \log x})}} \sqrt{s+1} \mathrm{d}s
    \\&+ \frac{m^3}{x (\log \log x)^{1/4}} \int_{\frac{x}{(m+1)^2 \exp (\sqrt{\log \log x})}}^{\frac{x}{m^2} - \frac{x}{(m+1)^2}} \sqrt{s+1} \log \left( 2 + \frac{x}{(m+1)^2 s} \right)^{1/2} \mathrm{d}s,
    \end{align*}
    implicitly using that
    \[ \frac{x}{(m+1)^2 \exp (\sqrt{\log \log x})} < \frac{x}{m^2} - \frac{x}{(m+1)^2} \]
    on the range of $m$ in question.
    The first integral here is
    \[ \ll \frac{m^3}{x} \left( \frac{x}{(m+1)^2 \exp (\sqrt{\log \log x})} \right)^{3/2} \ll \frac{\sqrt{x}}{\exp \left( \sqrt{\log \log x} \right)}, \]
    which is more than sufficient on the required range of $m \le \left[ (\log \log x)^{1/4}\right]$. For the second integral, we substitute
    \[ u = \frac{x}{(m+1)^2 s}, \qquad \mathrm{d}s = - \frac{x}{(m+1)^2 u^2} \mathrm{d}u. \]
    This yields
    \begin{align*}
        &\mathbb{E} \bigg[ \bigg| \sum_{\frac{x}{(m+1)^2} <  n \le \frac{x}{m^2}} X(n) \left( \sqrt{ \frac{x}{n} } - m \right) \bigg| \bigg] \\
        &\ll \frac{\sqrt{x}}{\exp \left( \sqrt{\log \log x} \right)} + \frac{\sqrt{x}}{(\log \log x)^{1/4}} \int_{\frac{m}{3}}^{\exp(\sqrt{\log \log x})} \sqrt{\log (u+2)} \frac{\mathrm{d}u}{u^{5/2}} \\
        &\ll \frac{\sqrt{x}}{(\log \log x)^{1/4}} \frac{\sqrt{\log (m+1)}}{m^{3/2}},
    \end{align*}
    implying \eqref{FracPartShortInt} and thus the claim.
\end{proof}

\section{Passing to Random Euler Products}\label{sec:pass_to_euler}

For any $y, z, \sigma,h\in \mathbb{R}$, let
\[
    F_{y, z}(\sigma+ih)=\prod_{\exp(e^y) < p \leq \exp(e^z)}\bigg(1+\frac{X(p)}{p^{1/2+\sigma+ih}}\bigg)
\]
and adopt the convention $F_z:=F_{-1,z}$. The following three propositions will allow us to reduce the analysis of moments of $\sum_{n\leq x} X(n)/\sqrt{n}$ and $\sum_{n\leq x} X^e(n)$ to those of (weighted) integrals of $|F|^2$. We then show how to bound the latter in the next section.

Note that despite $F$ not being the Euler product associated to $X^e$, Proposition \ref{lower_extended_reduction} nonetheless passes to a weighted integral of $|F|^2$, due to the ratio of 
\[ F_{y,z}^e(\sigma+ih) := \prod_{\exp(e^y) < p\leq \exp(e^z)} \bigg(1-\frac{X(p)}{p^{1/2+\sigma+ih}}\bigg)^{-1} \]
and $F_{y,z}(\sigma+ih)$ being, up to a small error, a deterministic sum over primes that can be evaluated asymptotically.

\begin{prop}\label{upper_reduction}
    Let $x$ be sufficiently large and set $\mathcal{K} = \lceil\log\log\log x \rceil$ and $\sigma_k = -\frac{2(k+1)}{\log x}$ as well as $t = \log \log x$. Uniformly for $\frac{1}{3} \le q \le \frac{1}{2}$ (say) and large enough $x$, we have
    \begin{equation}\label{upper_reduction_small}\begin{split}\mathbb{E}\bigg[\bigg|\sum_{n \leq x}\frac{X(n)}{\sqrt{n}}\bigg|^{2q}\bigg] 
    &\ll \sum_{0 \leq k \leq \mathcal{K}}\mathbb{E}\left[\left(\frac{1}{\log x}\int_{-\infty}^{\infty}\left|\frac{F_{t-k}(\sigma_k+ih)}{\sigma_k+ih}\right|^2\mathrm{d}h\right)^q\right] + (\log x)^{q (2q-1)}. 
    \end{split}\end{equation}
    Uniformly for $\frac{1}{2} \le q \le 1$, we have
    \begin{equation}\label{upper_reduction_large}\begin{split}\mathbb{E}\bigg[\bigg|\sum_{n \leq x}\frac{X(n)}{\sqrt{n}}\bigg|^{2q}\bigg]^{1/2q}
    &\ll \sum_{0 \leq k \leq \mathcal{K}}\mathbb{E}\left[\left(\frac{1}{\log x}\int_{-\infty}^{\infty}\left|\frac{F_{t-k}(\sigma_k+ih)}{\sigma_k+ih}\right|^2\mathrm{d}h\right)^q\right]^{1/2q} \! \! \! + (\log x)^{q-1/2}.
    \end{split}\end{equation}
\end{prop}
\begin{prop} \label{lower_extended_reduction}
    Let $x$ be sufficiently large and $0 < q \leq 1/2$ fixed. Set $\tilde t = \log \log x - \log 2$. There exists an absolute constant $A>0$ such that, for sufficiently large $U \le \log \log x$, we have
    \begin{align*}
        \mathbb{E}\left[\left|\frac{1}{\sqrt{x}}\sum_{n \leq x}X^{e}(n)\right|^{2q}\right] &\gg \mathbb{E}\left[\left(\frac{1}{\log x}\int_{-1/2}^{1/2}\left|\frac{F_{\tilde t}(4U/\log x+ih)}{{4U}/{\log x}+ih}\right|^2\mathrm{d}h\right)^q\right]\\
        &- A \mathbb{E}\left[\left(\frac{1}{e^{U}\log x}\int_{-\infty}^{\infty}\left|\frac{F_{\tilde t}(2U/\log x+ih)}{{2U}/{\log x}+ih}\right|^2 \log^2 (|h|+2)\mathrm{d}h\right)^q\right] -\frac{A}{(\log x)^q}.
    \end{align*}
\end{prop}
\begin{prop} 
    Let $x$ be sufficiently large and $0 < q \leq 1$ fixed. Set $\tilde t = \log \log x - \log 2$. There exists an absolute constant $A>0$ such that, for sufficiently large $U \le \log \log x$, we have
    \begin{align*}
        \mathbb{E}\left[\left|\sum_{n \leq x}\frac{X(n)}{\sqrt{n}}\right|^{2q}\right] &\gg \mathbb{E}\left[\left(\frac{1}{\log x}\int_{-1/2}^{1/2}\left|\frac{F_{\tilde t}(4U/\log x+ih)}{{4U}/{\log x}+ih}\right|^2\mathrm{d}h\right)^q\right]\\
        &- A \mathbb{E}\left[\left(\frac{1}{e^{U}\log x}\int_{-\infty}^{\infty}\left|\frac{F_{\tilde t}(2U/\log x+ih)}{{2U}/{\log x}+ih}\right|^2\mathrm{d}h\right)^q\right] - \frac{A}{(\log x)^q}.
    \end{align*} \label{lower_weighted_reduction}
\end{prop}

We note that the right-hand sides in Propositions \ref{lower_extended_reduction} and \ref{lower_weighted_reduction} are potentially negative, in which case we interpret the above $\gg$ statements as void. By taking $x$ sufficiently large, we will be able to ensure that this does not happen in the final estimates.

We omit the proof of Proposition \ref{lower_weighted_reduction}, as it is essentially identical to that of \cite[Proposition 2.5]{G1}, as well as that of (an easier version of) Proposition \ref{lower_extended_reduction}. The range of $0 < q \le 1$ has no effect on the proof aside from changing the absolute constant.

Lastly, we note that on the range $q \le \frac{1}{2}$, it seems likely that one can replace the statement of Proposition \ref{upper_reduction} by a similar expression without summing over $k$, akin to the argument leading to \cite[Lemma 1.2]{GorodetskyWongHelson}. This strategy does not appear to be sufficient for larger $q$ since the resulting contribution from small values of $h$ would turn out to be too large. 

\begin{proof}[Proof of Proposition \ref{upper_reduction}]
    The proof proceeds in exactly the same way as that in the Steinhaus case \cite[Propositions 2.3 and 2.4]{G1} for $q \le 1/2$, with $\alpha = 1$. Combining these two parts, one obtains that
    \begin{align*}
        \mathbb{E}\bigg[\bigg| \sum_{n \leq x}\frac{X(n)}{\sqrt{n}}\bigg|^{2q}\bigg] 
        &\ll \sum_{0 \leq k \leq \mathcal{K}}\mathbb{E}\left[\left(\frac{1}{\log x}\int_{-\infty}^{\infty}\left|\frac{F_{t-k}(\sigma_k+ih)}{\sigma_k+ih}\right|^2 \mathrm{d}h\right)^q\right] \\
        &+ \sum_{0 \le k \le \mathcal{K}} e^{-e^k q} \mathbb{E}[|F_{t-k-1}(0)|^{2q}] + (\log x)^{q(2q-1)},
    \end{align*}
    for $\sigma_k=-2(k+1)/\log x$.
    Note that the last error term here was only bounded by $1$ in \cite[Proposition 2.3]{G1}, but the proof readily provides a much stronger bound (by increasing the constant in the application of Rankin's trick).
    Note also that the bounds become uniform in $q$ once one assumes that $q$ is bounded away from $0$.
    The claim then follows readily from an application of Lemma \ref{lem:MomentsEulerProducts} to bound the Euler product moments.
    For the range of $\frac{1}{2} \le q \le 1$, one replaces all applications of subadditivity of $x \mapsto x^{2q}$ by Minkowski's inequality.
\end{proof}

The proof of Proposition \ref{lower_extended_reduction} will rely on the following second moment bound for sums of $X^e$ in short intervals.

\begin{lem}\label{lem:short_intervals} Uniformly in $x \ge y > 0$, we have
\[
    \mathbb{E}\Big[\Big|\sum_{x-y<n\leq x}X^e(n)\Big|^2\Big]\ll \min\Big\{(y^2+1), \frac{y^2}{x}\Big(\log_+\frac{y^2}{x}+1\Big)+y\Big(\log\frac{x}{y}+1\Big)^2+\Big(\frac{x}{y}\Big)^2\Big\},
\]
with $\log_+x=\max(\log x,0)$.
Moreover, the left-hand side is increasing as a function of $y$ for a given $x$.
\end{lem}
\begin{proof}
    The bound $\ll(y^2+1)$ follows trivially, noting that $|X^e(n)|^2 = 1$. For the second bound, we express each integer $n$ as a product $d^2m$ and use the fact that $X^e(n)=X(m)$ to write
    \begin{align*}
                 \mathbb{E}\Big[\Big|\sum_{x-y<n\leq x}X^e(n)\Big|^2\Big]&= \mathbb{E}\Big[\Big|\sum_{1\leq m\leq x}X(m)\sum_{\sqrt{\frac{x-y}{m}}<d\leq \sqrt{\frac{x}{m}}}1\Big|^2\Big]\\
                 &=\sum_{1\leq m\leq x}\mu(m)^2\bigg(\sum_{\sqrt{\frac{x-y}{m}}<d\leq \sqrt{\frac{x}{m}}}1\bigg)^2.
    \end{align*}
    From this, the monotonicity in $y$ is immediately apparent.
    
    For the first claim, we split this quantity according to whether $m\leq y^2/{x}$ or $m> y^2/{x}$, roughly determining whether the inner sum is (typically) empty or not. 

    In the first case, we have that
    \[
        \sqrt{\frac{x}{m}}-\sqrt{\frac{x-y}{m}}\gg 1,
    \]
    and it follows that 
    \[
        \sum_{1\leq m\leq \frac{y^2}{x}}\bigg(\sum_{\sqrt{\frac{x-y}{m}}<d\leq \sqrt{\frac{x}{m}}}1\bigg)^2 \ll \sum_{1\leq m\leq \frac{y^2}{x}} \Big(\frac{y^2}{mx}+1\Big)\ll \frac{y^2}{x}\bigg(\log_+\Big(\frac{y^2}{x}\Big)+1\bigg).
    \]
    For $m>y^2/{x}$, we exchange the order of summation to get 
    \begin{align*}
        \sum_{\frac{y^2}{x}<m\leq x} \sum_{\sqrt{\frac{x-y}{m}}<d_1,d_2\leq \sqrt{\frac{x}{m}}} 1 
            &\leq \sum_{d_1,d_2\leq \frac{x}{y}} \sum_{\frac{x-y}{d_1d_2}<m<\frac{x}{d_1d_2}}  1\\
        &\leq \sum_{d_1,d_2\leq \frac{x}{y}} \Big(\frac{y}{d_1d_2}+1\Big)\ll y\Big(\log\Big(\frac{x}{y}\Big)+1\Big)^2+\Big(\frac{x}{y}\Big)^2,
    \end{align*}
    with the range of $m$ stemming from the fact that $\sqrt{x-y}/d_j<\sqrt{m}\leq \sqrt{x}/d_j$ for $j\in \{1,2\}$.

\end{proof}

\begin{proof}[Proof of Proposition \ref{lower_extended_reduction}]
    We adapt ideas from the proof of Proposition 2 in \cite{hnr}, and begin by showing that 
    \begin{equation}\label{eq:reducetolargeprime}
               \mathbb{E}\bigg[\Big|\sum_{n \leq x}X^{e}(n)\Big|^{2q}\bigg] \gg \mathbb{E}\bigg[\Big|\sum_{\substack{n \leq x \\ P(n) > \sqrt{x}}}X^{e}(n)\Big|^{2q}\bigg].        
    \end{equation}
    After re-expressing the right-hand side as
    \begin{equation*}
        \frac{1}{2^{2q}}\mathbb{E}\bigg[\Big|\sum_{\substack{n \leq x \\ P(n) > \sqrt{x}}}X^{e}(n) + \sum_{\substack{n \leq x \\ P(n) \leq \sqrt{x}}}X^{e}(n) + \sum_{\substack{n \leq x \\ P(n) > \sqrt{x}}}X^{e}(n)-\sum_{\substack{n \leq x \\ P(n) \le \sqrt{x}}}X^{e}(n)\Big|^{2q}\bigg],
    \end{equation*}
    we can bound it above by
    \begin{align*}
        &\mathbb{E}\bigg[\Big|\sum_{\sqrt{x} < p \leq x}X(p)\sum_{n \leq x/p}X^{e}(n) + \sum_{\substack{n \leq x \\ P(n) \leq \sqrt{x}}}X^{e}(n)\Big|^{2q}\bigg] \\&+ \mathbb{E}\bigg[\Big|\sum_{\sqrt{x} < p \leq x}X(p)\sum_{n \leq x/p}X^{e}(n)-\sum_{\substack{n \leq x \\ P(n) \le \sqrt{x}}}X^{e}(n)\Big|^{2q}\bigg].
    \end{align*}
    Both of these expectations are then equal to 
    \[
    \mathbb{E}\biggl[\bigg|\epsilon\sum_{\sqrt{x} < p \leq x}X(p)\sum_{n \leq x/p}X^{e}(n) + \sum_{\substack{n \leq x \\ P(n) \leq \sqrt{x}}}X^{e}(n)\Big|^{2q}\bigg]
    \]
    for $\epsilon$ a Rademacher random variable independent from everything else. Noting that conditionally on $(X(p))_{p \le \sqrt{x}}$, the families $(X(p))_{\sqrt{x} < p \le x}$ and $(\epsilon X(p))_{\sqrt{x} < p \le x}$ have the same law, we conclude that
    \[
            \mathbb{E}\biggl[\Big|\epsilon\sum_{\sqrt{x} < p \leq x}X(p)\sum_{n \leq x/p}X^{e}(n) + \sum_{\substack{n \leq x \\ P(n) \leq \sqrt{x}}}X^{e}(n)\Big|^{2q}\bigg]
=\mathbb{E}\bigg[\Big|\sum_{n \leq x}X^{e}(n)\Big|^{2q}\bigg],
    \]
    and \eqref{eq:reducetolargeprime} follows. 
    
    Next, we apply a conditional version of Khintchine's inequality (see, e.g., \cite[Lemma 3.8.1]{gut06}), which we note holds for any $q >0$, and compute the expectation conditionally on $(X(p))_{p\leq \sqrt{x}}$ to get
    \[ \mathbb{E}\biggl[\bigg|\sum_{\substack{n \le x \\ P(n) > \sqrt{x}}}X^e(n)\bigg|^{2q}\biggl] \gg \mathbb{E}\biggl[\biggl(\sum_{\sqrt{x} < p \le x}\bigg|\sum_{n \le x/p}X^e(n)\biggl|^2\biggl)^q\biggl] \gg \frac{1}{(\log x)^q}\mathbb{E}\biggl[\biggl(\sum_{\sqrt{x} < p \le x}\log p\bigg|\sum_{n \le x/p}X^e(n)\biggl|^2\biggl)^q\biggl].\]
    Focusing on the term inside the $q$-power, we follow the proof of \cite[Proposition 3]{Harper1} to smoothen the sum over primes into an integral. For a parameter $Z = \exp(\sqrt{\log x})$, this gives 
    \begin{equation}\label{eq:integral_smooth}\begin{split}
        \mathbb{E}\Bigg[\bigg(\sum_{\sqrt{x} < p \le x}\log p \,&\bigg|\sum_{\substack{n \le x/p}}X^e(n)\bigg|^{2}\bigg)^q\Bigg] \gg \mathbb{E}\Bigg[\bigg(\sum_{\sqrt{x} < p \le x}\frac{1}{4}\frac{Z\log p}{p}\int_p^{p(1+1/Z)}\bigg|\sum_{n \le x/s}X^e(n)\bigg|^2\mathrm{d}s\bigg)^q\Bigg] \\
        &-\mathbb{E}\Bigg[\bigg(\sum_{\sqrt{x} < p \le x}\frac{Z\log p}{p}\int_{p}^{p(1+1/Z)}\bigg|\sum_{x/s < n \le x/p}X^e(n)\bigg|^{2}\mathrm{d}s\bigg)^q\Bigg].
    \end{split}
    \end{equation}
    To bound the subtracted term, note that by H\"older's inequality and Fubini's theorem,
    \begin{align*} 
    &\mathbb{E}\biggl[\bigg(\sum_{\sqrt{x} < p \le x}\frac{Z\log p}{p} \int_{p}^{p(1+1/Z)}\bigg|\sum_{x/s < n \le x/p}X^e(n)\bigg|^{2}\mathrm{d}s\bigg)^q\biggl] \\ 
    &\ll \bigg( \sum_{\sqrt{x} < p \le x} \frac{Z \log p}{p} \int_{p}^{p(1+1/Z)} \mathbb{E} \Big[ \Big| \sum_{x/s < n \le x/p} X^e(n) \Big|^2 \bigg]\mathrm{d}s \bigg)^q.
    \end{align*}
    For $x/Z^2<p\leq x$, the bound in Lemma \ref{lem:short_intervals} (for $y=x/p-x/s\ll x/pZ$) yields
    \[
        \sum_{x/Z^2< p \le x} \frac{Z \log p}{p} \int_{p}^{p(1+1/Z)} \mathbb{E} \Big[ \Big| \sum_{x/s < n \le x/p} X^e(n) \Big|^2 \bigg]\mathrm{d}s \ll \sum_{x/Z^2< p \le x}\log p\bigg(\Big(\frac{x}{pZ}\Big)^2+1\bigg) \ll x,
    \]
    using standard estimates on prime sums (see, e.g., \cite[Corollary 3.4]{MV07}). For the remaining $p$, we first use the monotonicity  
    claim in Lemma \ref{lem:short_intervals} to replace $s$ by $p(1+1/Z)$ in the integrand (for an upper bound), followed by the bound in said lemma to get
    \[
        \mathbb{E} \Big[ \Big| \sum_{\frac{x}{p(1+1/Z)} < n \le x/p} X^e(n) \Big|^2 \bigg] \ll \Big(\frac{x}{pZ}(\log x)^2+Z^2\Big).
    \]
    It follows that
    \begin{align*}
                \sum_{\sqrt{x}<p\leq x/Z^2}  \frac{Z \log p}{p} \int_{p}^{p(1+1/Z)} \mathbb{E} \Big[ \Big| \sum_{x/s < n \le x/p} X^e(n) \Big|^2 \bigg]\mathrm{d}s &\ll \sum_{\sqrt{x}<p\leq x/Z^2} \log p \Big(\frac{x}{pZ}(\log x)^2+Z^2\Big)\\
                &\ll \frac{x(\log x)^3}{Z}+x,
    \end{align*}
    and in turn that the subtracted term in \eqref{eq:integral_smooth} is $\ll x^q$.
    
    As for the first term in \eqref{eq:integral_smooth}, we exchange the order of summation and integration on the outer sum and perform the change of variable $z = x/s$ to find
    \[ \mathbb{E}\Bigg[\bigg|\sum_{\substack{n \le x \\ P(n) > \sqrt{x}}}X^e(n)\bigg|^{2q}\Bigg] \gg \frac{x^q}{\log^q x} \mathbb{E}\biggl[\bigg(\int_{1}^{\frac{\sqrt{x}}{(1+1/Z)}}\bigg|\sum_{n \le z}X^e(n)\bigg|^2\frac{\mathrm{d}z}{z^2}\bigg)^q\biggl]-\frac{Ax^q}{(\log x)^q},
    \]
    for some absolute constant $A>0$.
    Finally, we extend the range of integration in the above in order to apply Parseval's theorem (in the form of \cite[Equation (5.26)]{MV07}). For $U>0$ a large parameter, we have
    \begin{equation*}\begin{split}
        &\mathbb{E}\biggl[\bigg(\int_{1}^{\frac{\sqrt{x}}{(1+1/Z)}}\bigg|\sum_{n \le z}X^e(n)\bigg|^2\frac{\mathrm{d}z}{z^2}\bigg)^q\biggl] \ge \mathbb{E}\biggl[\bigg(\int_{1}^{\frac{\sqrt{x}}{(1+1/Z)}}\bigg|\sum_{\substack{n \le z \\ P(n) \le \sqrt{x}}}X^e(n)\bigg|^2\frac{\mathrm{d}z}{z^{2+8U/\log x}}\bigg)^q\biggl] \\
        \ge& \frac{1}{2^q}\mathbb{E}\biggl[\bigg(\int_{1}^\infty \bigg|\sum_{\substack{n \le z \\ P(n) \le \sqrt{x}}}X^e(n)\bigg|^2\frac{\mathrm{d}z}{z^{2+8U/\log x}}\bigg)^q\biggl] - \mathbb{E}\biggl[\bigg(\int_{\frac{\sqrt{x}}{(1+1/Z)}}^\infty \bigg|\sum_{\substack{n \le z \\ P(n) \le \sqrt{x}}}X^e(n)\bigg|^2\frac{\mathrm{d}z}{z^{2+8U/\log x}}\bigg)^q\biggl] \\
        \gg& \mathbb{E}\biggl[\bigg(\int_{1}^\infty \bigg|\sum_{\substack{n \le z \\ P(n) \le \sqrt{x}}}X^e(n)\bigg|^2\frac{\mathrm{d}z}{z^{2+8U/\log x}}\bigg)^q\biggl] - \frac{1}{e^{2Uq}}\mathbb{E}\biggl[\bigg(\int_{1}^\infty \bigg|\sum_{\substack{n \le z \\ P(n) \le \sqrt{x}}}X^e(n)\bigg|^2\frac{\mathrm{d}z}{z^{2+4U/\log x}}\bigg)^q\biggl],
    \end{split}
    \end{equation*}
    and can now apply \cite[Equation (5.26)]{MV07} to each of these integrals to get
    \begin{equation} \label{extended_euler}
    \begin{split}
        \mathbb{E} \bigg[ \bigg|\sum_{n \le x}X^e(n)\bigg|^{2q} \bigg] \gg\,\,& \mathbb{E}\left[\left(\frac{x}{\log x}\int_{-1/2}^{1/2}\left|\frac{F_{\tilde t}^e(4U/\log x+ih)}{1/2+ih}\right|^2 \mathrm{d}h \right)^q\right] \\
        - A &\mathbb{E}\left[\left(\frac{x}{e^{U}\log x}\int_{-\infty}^{\infty}\left|\frac{F_{\tilde t}^e(2U/\log x+ih)}{1/2+ih}\right|^2 \mathrm{d}h \right)^q\right] -\frac{Ax^q}{(\log x)^q} \end{split}
    \end{equation}
    for $A>0$ a large absolute constant.
    All that is now left to do is to relate these integrals of $F^e$ to integrals of $F$. For any $\sigma \asymp \frac{U}{\log x}$, we have 
    \begin{equation*} \begin{split}
    \prod_{p \leq \sqrt{x}}\bigg|1-\frac{X(p)}{p^{1/2+\sigma+ih}}\bigg|^{-2}\bigg|1+\frac{X(p)}{p^{1/2+\sigma+ih}}\bigg|^{-2} &= \prod_{p \leq \sqrt{x}}\left|1-\frac{1}{p^{1+2\sigma+2ih}}\right|^{-2} \\ &\asymp \exp \bigg(2\sum_{p \leq \sqrt{x}}\frac{\cos(2h \log p)}{p^{1+2\sigma}}\bigg). \end{split}
\end{equation*}
Applying Lemma \ref{lem:CosineSum}, it follows that the right-hand side in (\ref{extended_euler}) is
\begin{equation*}
    \begin{split}
        &\gg\mathbb{E}\left[\left(\frac{x}{\log x}\int_{-1/2}^{1/2}\left|\frac{F_{\tilde t}(4U/\log x+ih)}{\frac{4U}{\log x}+ih}\right|^2 \mathrm{d}h \right)^q\right] \\
        &\quad- A \mathbb{E}\left[\left(\frac{x}{e^{U}\log x}\int_{-\infty}^{\infty}\left|\frac{F_{\tilde t}(2U/\log x+ih)}{\frac{2U}{\log x}+ih}\right|^2 \log^2(|h|+2) \mathrm{d}h\right)^q\right] -\frac{A x^q}{\log^q x }, \end{split}
    \end{equation*}
    after potentially replacing $A$ by a larger absolute constant. 
\end{proof}

\section{Upper Bounds for Random Euler Products} \label{sec:upper}
The aim of this section is to prove the upper bounds in \eqref{q_uniformity} in the weighted Rademacher setting. In view of Propositions \ref{ExtToWeighted} and \ref{BoundFractPart}, this will imply the same bounds for the extended Rademacher setting. 

Throughout this section, and the rest of this work, $x$ will be a sufficiently large parameter, $\sigma$ a (small) real number and we set $t_\sigma = [\log \min \{ \log x, |\sigma|^{-1}\}]$, with the understanding that this minimum is given by the first term when $\sigma = 0$.

By Proposition \ref{upper_reduction}, it suffices to obtain upper bounds for integrals of the shape
\[ \mathbb{E}\left[\left(\frac{1}{\log x}\int_{-\infty}^{\infty}\left|\frac{F_{t-k}(\sigma_k+ih)}{\sigma_k+ih}\right|^2\mathrm{d}h\right)^q\right]. \]
We will instead prove bounds for the slightly larger expressions
\[ \mathbb{E}\left[\left(\frac{1}{\log x}\int_{-\infty}^{\infty}\left|\frac{F_{t-k}(\sigma_k+ih) \log (|h|+2)}{\sigma_k+ih}\right|^2\mathrm{d}h\right)^q\right] \]
with the aim of reusing them to obtain upper bounds for the subtracted terms in Propositions \ref{lower_extended_reduction} and \ref{lower_weighted_reduction}. The proofs will make use of a random Euler product moment bound due to Harper, which we record here for future reference\footnote{The statement there has very slightly different parameter choices, but the proof readily gives the result below without changes to the argument.}.

\begin{lem}\label{lem:constant_chaos}
    (Multiplicative Chaos Result 3, \cite{Harper4}) Suppose that $x$ is sufficiently large, let $\sigma  \in [-10/\log x,1/(\log x)^{0.01}]$ and set $t_\sigma = [\log \min \{ \log x, |\sigma|^{-1}\}]$. Uniformly for $|N| \leqslant (\log x)^{1000}$ and $q \in [0,1]$, we have
    \begin{equation*}
        \mathbb{E}\bigg[\Big(\int_{N-1/2}^{N+1/2}|F_t(\sigma+ih)|^2 \mathrm{d}h\Big)^q\bigg] \ll \big(\log\log(|N|+10)\big)^q\Big(\frac{e^{t_\sigma}}{1+(1-q)\sqrt{t_\sigma}}\Big)^q.
    \end{equation*}
\end{lem}

Throughout this section, we will continue to assume that $q \ge \frac{1}{3}$, since the upper bound for smaller $q$ follows from a simple application of H\"older's inequality together with the bound at $q = \frac{1}{3}$.

Another key input will be the following adaptation of the main result in \cite{Harper1} to \textit{mesoscopic intervals} $I$, meaning those $I\subseteq [0,1]$ satisfying $(\log x)^{-1} \le |I| \le 1$ as $x\to\infty$. The proof is deferred to Section \ref{sec:harper_mesoscopic}. 

\begin{prop}\label{prop:HarperMesoscopic} Let $x$ be sufficiently large, and $-\frac{10}{\log x} \le \sigma \le \frac{\log \log x}{\log x}$. Uniformly for $e^{-t_\sigma} \le H \le 1$ and $0 \le q \le 1$, we have
\[
    \mathbb{E}\bigg[\bigg( \frac{1}{H}\int_{H}^{eH} |F_{\log 1/H, t_\sigma}(\sigma+ih)|^2 \mathrm{d}h\bigg)^{q}\bigg]  \ll \Big( \frac{He^{t_\sigma}}{1+(1-q)\sqrt{t_\sigma + \log H}} \Big)^q.
\]
\end{prop}

\subsection{Reduction to a mesoscopic integral moment}

\begin{prop}\label{prop:upper_uniform_REP}
    Let $x$ be sufficiently large and $\frac{1}{3} \le q \le 1$. Suppose that $-10/\log x\le \sigma \le \log\log x/\log x$ and set $t = \log \log x$ as well as $Y = Y(q,x) = \min \{t, |2q-1|^{-1} \}$. Uniformly in $x, q$ and $\sigma$, we have
    \begin{equation}\label{typical_q}\begin{split}
        &\mathbb{E}\bigg[\Big(\frac{1}{\log x}\int_{\mathbb{R}} \bigg| \frac{F_t(\sigma+ih)}{\sigma+ih}\bigg|^2 \log^2(|h|+2) \mathrm{d}h\Big)^{q}\bigg] \ll (\log x)^{q(2q-1)} \left \{ (|\sigma|\log x)^{-q} +  Y^{3/4} \right\} + t^{-q/2} Y.
    \end{split}\end{equation}
\end{prop}

\begin{proof}
We begin by noting that since $|F_t(\sigma +ih)|=|F_t(\sigma-ih)|$ for every $h\in \mathbb{R}$, we may restrict the range of integration to $h \ge 0$.
We continue by dividing the range of integration into several subranges according to the size of $h$ as follows. Set
\[ I_1 = [0,e^{-t_\sigma}], \quad I_2 = [e^{-t_\sigma}, 1], \quad \quad I_3 = [1,t^2], \quad I_4 = [t^2,\infty) \]
and for $1 \le i \le 4$, define 
\[ S_i = \mathbb{E}\bigg[\Big(\frac{1}{\log x}\int_{I_i} \frac{|F_t(\sigma+ih)|^2\log^2(|h|+2)}{|\sigma+ih|^2}\mathrm{d}h\Big)^{q}\bigg]. \]
Using the subadditivity of $x\mapsto x^{q}$, we have
\begin{equation}\label{REPIntegral_upper}
    \mathbb{E}\bigg[\Big(\frac{1}{\log x}\int_{\mathbb{R}} \frac{|F_t(\sigma+ih)|^2\log^2(|h|+2)}{|\sigma+ih|^2}\mathrm{d}h\Big)^{q}\bigg] \ll \sum_{i=1}^4 S_i.
\end{equation}

To bound $S_4$, we use H\"older's inequality and Lemmas \ref{lem:MomentsEulerProducts} and \ref{lem:PNT} to get
\begin{align*}
    S_4 &= \mathbb{E}\bigg[\Big(\frac{1}{\log x}\int_{ t^2}^{\infty} \frac{|F_t(\sigma + ih)|^2\log^2(h+2)}{|\sigma + ih|^2}\mathrm{d}h\Big)^{q}\bigg] \\ &\le \Big(\frac{1}{\log x} \int_{ t^2}^{\infty}\frac{\mathbb{E}[|F_t(\sigma + ih)|^2]\log^2(h+2)}{h^2}\mathrm{d}h\Big)^q \\
    &\ll \Big( \int_{ t^2}^{\infty} \frac{\log^2 h}{h^2}\mathrm{d}h\Big)^q
    \ll t^{-q},
\end{align*}
which can be absorbed into the third term on the right-hand side of \eqref{typical_q} with room to spare. 

Next, we will show that
    \[S_3 \ll \big(1+(1-q)\sqrt{t}\big)^{-q} \]
uniformly for $\frac{1}{3} \le q \le 1$ and $x$. This can be absorbed into the third term of \eqref{typical_q} for $q \le \frac{2}{3}$, and comfortably into the second term for $\frac{2}{3} \le q \le 1$.

We may thus assume that $q \ge \frac{2}{3}$, since the lower range follows from an application of H\"older's inequality and the statement at $q = \frac{2}{3}$.
By Lemma \ref{lem:constant_chaos} and noting that $t \asymp t_\sigma$, we then have
\begin{align*}
    S_3 &= \mathbb{E}\bigg[\Big(\frac{1}{\log x}\int_1^{t^2} \frac{|F_t(\sigma + ih)|^2\log^2(h+2)}{|\sigma+ih|^2}\mathrm{d}h\Big)^{q}\bigg] \\ 
    &\ll \sum_{N \le t^2} \frac{\log^{2q}(N+2)}{N^{2q}} \mathbb{E}\bigg[\Big( \frac{1}{\log x}\int_{N}^{N+1} |F_t(\sigma+ih)|^2\mathrm{d}h\Big)^q\bigg]\\
    &\ll \sum_{N \le (\log\log x)^2}\frac{\log^{2q} (N + 2) (\log \log (N+10))^q }{N^{2q}} \frac{\min \{ \log x, |\sigma|^{-1} \}^q}{(1+(1-q)\sqrt{t})^q(\log x)^q} \ll \big(1+(1-q)\sqrt{t}\big)^{-q}
\end{align*}
as required.

We continue by bounding $S_1$. On this range, the numerator is straight-forward to handle since $p^{-ih}$ does not vary much on the range in consideration, for any $p \le x$. We thus have
\begin{align*}
    S_1 &= \mathbb{E}\bigg[\Big(\frac{1}{\log x}\int_0^{e^{-t_\sigma}} \frac{|F_t(\sigma + ih)|^2\log^2(h+2)}{|\sigma+ih|^2}\mathrm{d}h\Big)^{q}\bigg] \\
    &\ll \frac{1}{(\log x)^q} \mathbb{E}\bigg[|F_t(\sigma)|^{2q(1-q)} \Big( \int_0^{e^{-t_\sigma}} \frac{|F_t(\sigma + ih)|^2 |F_t(\sigma)|^{-2(1-q)}}{|\sigma+ih|^2} \mathrm{d}h\Big)^{q}\bigg] \\
    &\le \frac{\mathbb{E}[|F_t(\sigma)|^{2q}]^{1-q}}{( \log x)^q} \bigg( \int_0^{e^{-t_\sigma}} \frac{\mathbb{E}[|F_t(\sigma + ih)|^2 |F_t(\sigma)|^{-2(1-q)}]}{|\sigma+ih|^2} \mathrm{d}h \bigg)^{q}.
\end{align*}
Applying Lemma \ref{lem:MomentsEulerProducts}, we have $\mathbb{E}[|F_t(\sigma)|^{2q}]^{1-q}\asymp e^{t_\sigma q(2q-1)(1-q)}$, and
\[ \mathbb{E}[|F_t(\sigma + ih)|^2 |F_t(\sigma)|^{-2(1-q)}] \asymp \exp \left( \sum_{p \le x} \frac{4-5q+2q^2-4(1-q) \cos(h \log p)}{p^{1+2\sigma}} \right) \asymp e^{t_\sigma q(2q-1)} \]
by splitting the sum at $p = \min \{ x, e^{1/|\sigma|} \}$ and estimating $\cos(x) = 1 + O(x^2)$. Hence, we deduce that
\begin{align*}
    S_1 &\ll \frac{e^{t_\sigma q(2q-1)(1-q)}}{(\log x)^q} \Big( \int_{0}^{e^{-t_\sigma}} \frac{e^{t_\sigma q(2q-1)}}{|\sigma+ih|^2} \mathrm{d}h \Big)^q \ll e^{t_\sigma q(2q-1)} (|\sigma| \log x)^{-q}.
\end{align*} 
When $|\sigma| \le \frac{1}{\log x}$, we have $e^{t_\sigma} \asymp \log x$ and this can be absorbed into the first term on the right-hand side of \eqref{typical_q}. If on the other hand $|\sigma| > \frac{1}{\log x}$ then $e^{t_\sigma} \asymp |\sigma|^{-1}$, so that this whole expression is
\[ \le |\sigma|^{-2q^2} (\log x)^{-q} \le (\log x)^{q(2q-1)}, \]
which may be absorbed into the second term.

It remains to bound $S_2$, which is the main difficulty and the purpose of Proposition \ref{prop:HarperMesoscopic}. We cover this range by $e$-adic intervals $[H_j,H_{j+1})$ with endpoints $(H_j)_{0 \leq j\leq J}$, where $H_j=e^{j-t_\sigma}$ for each $j\geq 0$, and $J:=\min\{j: H_j\geq 1\} = t_\sigma $. This ensures in particular that on each interval $(H_j,H_{j+1}]$, the denominator $|\sigma+ih|^{-2}$ in the integrand in $S_2$ is of order $H_j^{-2}$. It follows that
\[
    S_2 \ll \mathbb{E}\bigg[\Big(\frac{1}{\log x}\int_{e^{-t_\sigma}}^1 \frac{|F_t(\sigma+ih)|^2}{|\sigma +ih|^2}\mathrm{d}h\Big)^{q}\bigg] \ll \sum_{j=0}^J \frac{1}{ (H_j \log x)^q} \mathbb{E}\Big[\Big( \frac{1}{H_j}\int_{H_j}^{H_{j+1}}|F_t(\sigma+ih)|^2\mathrm{d}h\Big)^{q}\Big].
\] 
On each block $(H_j,H_{j+1}]$, it is also a consequence of Lemma \ref{lem:MomentsEulerProducts} that the partial Euler products $F_{\log 1/H_j}(\sigma+ih)=\prod_{p\leq e^{1/H_j}}(1+X(p)p^{-1/2-\sigma-ih})$ are highly correlated across $h \in [0,H_{j+1}]$.  One should then morally be able to factor the contribution of $|F_{\log 1/H_j}(\sigma)|^2$, say, out of the integral. Moreover, we should be able to treat $F_{t_\sigma,t}(\sigma+ih)$ in an elementary manner, since all of its moments are bounded uniformly in $h$ and $q$ by Lemma \ref{lem:MomentsEulerProducts}.

Thus, let $\mathcal{F}_j$ denote the $\sigma$-algebra generated by the random variables $X(p)$ for primes $p$ satisfying $\log 1/H_j < \log \log p \le t_\sigma$. Using the tower property and H\"older's inequality conditionally on $\mathcal{F}_j$, we obtain
\begin{align*}
    &\mathbb{E}\Big[\Big( \frac{1}{H_j}\int_{H_j}^{H_{j+1}}|F_t(\sigma+ih)|^2\mathrm{d}h\Big)^{q}\Big] \\
    = \,  &\mathbb{E} \bigg[ \mathbb{E}\Big[ |F_{\log 1/H_j}(\sigma)|^{2q(1-q)} \Big( \frac{1}{H_j}\int_{H_j}^{H_{j+1}}|F_t(\sigma+ih)|^2 |F_{\log 1/H_j}(\sigma)|^{-2(1-q)} \mathrm{d}h\Big)^{q} \, \Big\vert \mathcal{F}_j \Big] \bigg] \\
    \le \, &\mathbb{E}[|F_{\log 1/H_j}(\sigma)|^{2q}]^{1-q} \mathbb{E} \bigg[ \Big( \frac{1}{H_j} \int_{H_j}^{H_{j+1}}\mathbb{E}[|F_{t_\sigma,t}(\sigma+ih)|^2] \times \\&\hspace{100pt} \mathbb{E}[|F_{\log 1/H_j}(\sigma+ih)|^2 |F_{\log 1/H_j} (\sigma)|^{-2(1-q)}] |F_{\log 1/H_j,t_\sigma}(\sigma+ih)|^2 \mathrm{d}h \Big)^q \bigg].
\end{align*}
Evaluating all of these expectations once more using Lemmas \ref{lem:MomentsEulerProducts} and \ref{lem:PNT} yields that this is in turn
\[ \ll H_j^{-q (2q-1)} \mathbb{E}\bigg[\bigg( \frac{1}{H_j}\int_{H_j}^{eH_j} |F_{\log 1/H_j, t_\sigma}(\sigma+ih)|^2 \mathrm{d}h\bigg)^{q}\bigg].  \]
It therefore follows that
\[ S_2 \ll (\log x)^{-q} \sum_{j=0}^J H_j^{-2q^2} \mathbb{E}\bigg[\bigg( \frac{1}{H_j}\int_{H_j}^{eH_j} |F_{\log 1/H_j, t_\sigma}(\sigma+ih)|^2 \mathrm{d}h\bigg)^{q}\bigg],  \]
which brings us into a position to apply Proposition \ref{prop:HarperMesoscopic}. In doing so, we obtain
\[ S_2 \ll (\log x)^{-q} \sum_{j = 0}^J H_j^{-2q^2} \bigg( \frac{H_j e^{t_\sigma}}{1+(1-q)\sqrt{t_\sigma+\log H_j}} \bigg)^q = \frac{e^{2t_\sigma q^2}}{(\log x)^q} \sum_{j=0}^J \left( \frac{e^{-j(2q-1)}}{1+(1-q)\sqrt{j}} \right)^q. \]
To finalise the bound on $S_2$, we split the range of $q$ into $3$ subranges.

We begin with the simplest range $\frac{2}{3} \le q \le 1$. In this case, we discard the denominator to obtain that
\[ S_2 \ll \frac{e^{2 t_\sigma q^2}}{(\log x)^q} \ll (\log x)^{q(2q-1)} (|\sigma| \log x)^{-q} \]
as required, using that the geometric sum converges uniformly in $q$ on that range. In the remaining cases, we can now simplify the denominator to $\sqrt{j+1}$.

Next, let us consider $\frac{1}{2} \le q \le \frac{2}{3}$ and set $J_0 = \min \{J, (2q-1)^{-1}\} \ll Y$. This is the range of $j$ up to which we can not expect substantial saving from the numerator, so we simply bound it by $1$. Bounding the exponent of the denominator from below by $1/4$, this gives
\begin{align*} S_2 \ll \frac{e^{2 t_\sigma q^2}}{(\log x)^q} \left( \sum_{j = 0}^{J_0} \frac{1}{(j+1)^{1/4}} + J_0^{-1/4} \sum_{j = J_0+1}^J e^{-jq(2q-1)} \right) &\ll Y^{3/4} \frac{e^{2 t_\sigma q^2}}{(\log x)^q} \ll (\log x)^{q(2q-1)} Y^{3/4}.
\end{align*}

Finally, suppose that $\frac{1}{3} \le q \le \frac{1}{2}$. In this case, we reverse the order of summation over $j$ to get
\[ S_2 \ll \left( \frac{e^{t_\sigma}}{\log x} \right)^q \sum_{j = 0}^J \frac{e^{-jq(1-2q)}}{(J-j+1)^{q/2}}. \]
We split the sum at $J/2$. For the first half of the sum itself, we have
$$\sum_{j \le J/2} (J-j+1)^{-q/2}e^{-q(1-2q)j} \ll J^{-q/2}\sum_{j \le J/2} e^{-q(1-2q)j} \ll t^{-q/2} Y. $$
For the second half of the sum, we also obtain
$$\sum_{J/2 < j \le J}(J-j+1)^{-q/2}e^{-q(1-2q)j} \ll e^{-q(1-2q)J/2}J^{1-q/2} \ll t^{-q/2}Y$$
using that $e^{-u} \ll \min \{1, u^{-1}\}$ for $u \ge 0$.
As a consequence, we conclude that on this range of $q$, we have $S_2 \ll t^{-q/2}Y$ and the claim follows.
\end{proof}

\subsection{Proof of the upper bounds in Theorems \ref{MainThm1_statement} and \ref{MainThm2_statement}}

The aim of this subsection is to prove the uniform upper bound \eqref{q_uniformity} in both the extended Rademacher and weighted Rademacher case. In particular, these immediately imply the upper bounds of Theorems \ref{MainThm1_statement} and \ref{MainThm2_statement}. Once again, by Propositions  \ref{ExtToWeighted} and \ref{BoundFractPart}, it suffices to prove the weighted Rademacher bound
\begin{align}\label{q_uniformity_weighted}
    \mathbb{E}\bigg[\Big| \sum_{n \le x}\frac{X(n)}{\sqrt{n}}\Big|^{2q}\bigg] \ll
	\begin{cases}
		Y t^{-q/2},  & \mbox{if } 0\leq q \le \tfrac{1}{2},\\
		Y^{3/4} (\log x)^{q(2q-1)}, & \mbox{if } \tfrac{1}{2}\le q \le 1,
	\end{cases}
\end{align}
with $t = \log \log x$ and $Y = \min \{t, |2q-1|^{-1} \}$.

\begin{proof}[Proof of \eqref{q_uniformity}, \eqref{q_uniformity_weighted} and the upper bounds of Theorems \ref{MainThm1_statement} and \ref{MainThm2_statement}]
We may assume, as before, that $q \ge \frac{1}{3}$.

Suppose first that $\frac{1}{3} \le q \le \frac{1}{2}$. Combining \eqref{upper_reduction_small} and Proposition \ref{prop:upper_uniform_REP} applied at $x_k = x^{e^{-k}}$ and $\sigma_k = \frac{-2(k+1)}{\log x} = \frac{-2(k+1)}{e^k \log x_k}$, we have
\begin{align*}
    \mathbb{E}\bigg[\bigg|\sum_{n \leq x}\frac{X(n)}{\sqrt{n}}\bigg|^{2q}\bigg]
    &\ll \sum_{0 \leq k \leq \mathcal{K}}e^{-kq} \mathbb{E}\left[\left(\frac{1}{\log x_k}\int_{-\infty}^{\infty}\left|\frac{F_{t-k}(\sigma_k+ih)}{\sigma_k+ih}\right|^2\mathrm{d}h\right)^q\right] + (\log x)^{q (2q-1)} \\
    &\ll \sum_{0 \leq k \leq \mathcal{K}} e^{-kq} (t-k)^{-q/2} \min \{ t-k, |1-2q|^{-1} \} + (\log x)^{q(2q-1)} \\ &+ \sum_{0 \leq k \leq \mathcal{K}} e^{-kq} (\log x_k)^{q(2q-1)} \Big\{ (|\sigma_k| \log x_k)^{-q} + \min \{ t-k, |1-2q|^{-1} \}^{3/4} \Big\}  \\
    &\ll Y t^{-q/2}  + (\log x)^{q(2q-1)} \sum_{0 \leq k \leq \mathcal{K}} \left(e^{-kq(2q-1)} + e^{-2kq^2}Y^{3/4}  \right) \\
    &\ll Y t^{-q/2}  + \mathcal{K} e^{q (1-2q) (\mathcal{K}-t)} + e^{t q(2q-1)} Y^{3/4}.
\end{align*}
It is an elementary computation to verify that all of these terms can be absorbed into the first one uniformly for $\frac{1}{3} \le q \le \frac{1}{2}$. 

To see this for the second term, note that $\mathcal{K} = \log t + O(1)$. Suppose first that $\frac{1}{2} - \frac{1}{2t} \le q \le \frac{1}{2}$, so that $Y = t$. We then have
\[ \mathcal{K} e^{q (1-2q) (\mathcal{K}-t)} \ll (\log t) e^{q (1-2q) (\log t - t)} \ll \log t \ll t^{1-q/2} = Y t^{-q/2}  \]
as required. On the complementary range $\frac{1}{3} \le q \le \frac{1}{2} - \frac{1}{2t}$ we have $Y = (1-2q)^{-1} $. In this case, we use the inequality $e^{-x} \le x^{-1} $ for $x \ge 0$ to deduce that 
\[ \mathcal{K} e^{q (1-2q) (\mathcal{K}-t)} \ll \frac{\log t}{q(1-2q) (t-\log t)} \asymp \frac{\log t}{(1-2q)t} \ll \frac{1}{t^{q/2} (1-2q) } \]
as needed.

Finally, to show that 
\[ e^{t q(2q-1)} Y^{3/4} \ll Y t^{-q/2} \]
uniformly for $\frac{1}{3} \le q \le \frac{1}{2}$, we can employ the same case distinction and inequalities to conclude the proof on this range.

Now suppose that $\frac{1}{2} \le q \le 1$. This time using \eqref{upper_reduction_large}, we obtain
\begin{equation}\begin{split}\label{sep25.1}
    &\mathbb{E}\bigg[\bigg|\sum_{n \leq x}\frac{X(n)}{\sqrt{n}}\bigg|^{2q}\bigg]^{1/2q} \\
    &\ll \sum_{0 \leq k \leq \mathcal{K}} e^{-k/2} \left[ (\log x_k)^{q(2q-1)} \left( (|\sigma_k| \log x_k)^{-q} + Y^{3/4} \right) + t^{-q/2} Y \right]^{1/2q} \\
    &\ll (\log x)^{q-1/2} \sum_{0 \leq k \leq \mathcal{K}} \frac{e^{-k(q-1/2)}}{\sqrt{k+1}} + (\log x)^{q-1/2} Y^{3/8q} \sum_{0 \leq k \leq \mathcal{K}} e^{-kq} + t^{-1/4} Y^{1/2q} \sum_{0 \leq k \leq \mathcal{K}} e^{-k/2} \\
    &\ll (\log x)^{q-1/2} Y^{3/8q} + t^{-1/4}Y^{1/2q} + (\log x)^{q-1/2} \sum_{0 \leq k \leq \mathcal{K}} \frac{e^{-k(q-1/2)}}{\sqrt{k+1}}.
\end{split}\end{equation}
Moreover, we have
\[ \sum_{0 \leq k \leq \mathcal{K}} \frac{e^{-k(q-1/2)}}{\sqrt{k+1}} \ll \min \{ \mathcal{K}, |1-2q|^{-1} \}^{1/2}, \]
for example by using an integral comparison and the identity
\[ \int_0^{\infty} \frac{e^{-au}}{\sqrt{u}} \mathrm{d}u = \sqrt{\pi/a} \]
for $a > 0$. This is in turn $\ll Y^{3/8q}$, as is easily seen by case distinction according to whether $q \le 3/4$. Hence, the third term on the last line of \eqref{sep25.1} may be absorbed into the first term. 

To see that the second term can also be absorbed into the first one, suppose first that $\frac{1}{2} \le q \le \frac{1}{2} + \frac{1}{2t}$, so that $Y = t$. This gives
\[ t^{-1/4} Y^{1/(2q)} = t^{1/(8q)-1/4} Y^{3/8q} \le Y^{3/(8q)} \le Y^{3/(8q)} (\log x)^{q-1/2}. \]
For $\frac{1}{2} + \frac{1}{2t} \le q \le 1$ we get $Y = (2q-1)^{-1}$ and, using $x^{-1/(8q)} \ll e^{x/2}$ for $x \ge 1$, that
\[ t^{-1/4}Y^{1/(2q)} \le (t (2q-1))^{-1/(8q)} Y^{3/(8q)} \ll e^{t(q-1/2)}Y^{3/(8q)} = (\log x)^{q-1/2} Y^{3/8q}. \]
The claim follows.

\end{proof}

\subsection{Upper bounds for random Euler products in mesoscopic intervals}\label{sec:harper_mesoscopic}
This section outlines the proof of Proposition \ref{prop:HarperMesoscopic}. We proceed in a very similar manner to \cite[Section 4]{Harper1}. Let $L = [ t_\sigma + \log H ]$ be the number of (doubly logarithmic) scales of primes involved in the Euler product, and take $B$ to be a sufficiently large integer constant. For $h \in [H,eH]$, set $h(L+1) = h$ and recursively define the approximation sequence (in reverse) as
$$h(\ell) :=\max\Big\{u \le h(\ell+1): u=\frac{n H}{\ell e^\ell}\text{ for some } n \in \mathbb{Z} \Big \}$$
for $B +2 \le \ell \le L$.
Then we define $\mathcal{G}_H(h)$ as the event that for all $B+2 \le m \le L$, we have
$$\Big|\sum_{\ell=B+2}^{m}\sum_{\frac{e^{\ell-1}}{H} < \log p \le \frac{e^{\ell}}{H}}\frac{X(p)}{p^{1/2+\sigma+ih(\ell)}}\Big| \le m + g(m)$$
where $g(m) = g(m,x,H,\sigma) =  C\min\{\sqrt{L}, (1-q)^{-1}\}+2 \log m$ with $C>0$ some large constant to be chosen later. Using this, we define $\mathcal{G}_H$ as the event that for all $h \in [H,eH]$, the event $\mathcal{G}_H(h)$ holds. One can then proceed in an identical fashion to that in \cite[Section 4]{Harper1} and obtain the following propositions.
\begin{prop}\label{prop:barriered_moment}
    Let $x$ be sufficiently large and $-\frac{10}{\log x} \le \sigma \le \frac{\log\log x}{\log x}$. Set $t_\sigma = \log \min \{ \log x, |\sigma|^{-1} \}$, assume that $e^{-t_\sigma} \le H \le 1$ and set $L = [ t_\sigma+\log H]$. Uniformly for $\frac{2}{3} \le q \le 1$ and sufficiently large $C$, we have
    $$\mathbb{E}\bigg[ \Big(\frac{\mathbf{1}(\mathcal{G}_H)}{H}\int_{H}^{eH}|F_{\log 1/H,t_\sigma}(\sigma+ih)|^2\mathrm{d}h\Big)^q \bigg] \ll \Big(\frac{C e^L}{1+(1-q)\sqrt{L}}\Big)^q.$$
\end{prop}
\begin{prop}\label{prop:good_event_fail}
    With the assumptions and notation as before, we have
    $$\mathbb{P}(\mathcal{G}_H^c) \ll e^{-2C\min\{\sqrt{L},(1-q)^{-1}\}}.$$
\end{prop}
With these two propositions in mind, one can follow the section ``Proof of the upper bound in Theorem 1, assuming Key Propositions 1 and 2'' to deduce the full bound in Proposition \ref{prop:HarperMesoscopic} for $q\in [2/3,1]$, with the modifications for the Rademacher case as outlined in Section 4.4 there. The bound for $q\in [0,2/3)$ follows by Hölder's inequality.

\section{Lower bounds for Random Euler Products}\label{sec:lower}
\subsection{The lower bound for fixed \texorpdfstring{$q \ne \frac{1}{2}$}{q \ne 1/2}}

\begin{proof}[Proof of the lower bounds in Theorems \ref{MainThm1_statement} and \ref{MainThm2_statement} for $q < \frac{1}{2}$]
    In view of Propositions \ref{lower_extended_reduction} and \ref{lower_weighted_reduction}, it suffices to show that
    \begin{equation}
        \mathbb{E}\bigg[\Big(\frac{1}{\log x}\int_{1/3}^{1/2}\left|F_{\tilde t}(4U/\log x+ih) \right|^2\mathrm{d}h\Big)^q\bigg] \gg  (U \sqrt{t} )^{-q}
    \end{equation}
    and
    \begin{equation}
        \mathbb{E}\left[\left(\frac{1}{e^U \log x}\int_{-\infty}^{\infty}\left|\frac{F_{\tilde t}(2U/\log x+ih)}{{2U}/{\log x}+ih}\right|^2 \log^2 (|h|+2)\mathrm{d}h\right)^q\right] \ll (e^U \sqrt{t})^{-q}.
    \end{equation}
    Taking $U$ to be a sufficiently large constant, the claim then follows.
    
    The first estimate here is \cite[p. 69, second display]{Harper1} in the Steinhaus case, with the Rademacher modifications outlined in Section 5.3 there.
    The second bound follows from Proposition \ref{prop:upper_uniform_REP}, noting that $Y \asymp 1$ for fixed $q \ne 1/2$.
\end{proof}

\begin{proof}[Proof of the lower bounds in Theorems \ref{MainThm1_statement} and \ref{MainThm2_statement} for $q > \frac{1}{2}$]
    Firstly, note that it suffices to show the weighted Rademacher case, Theorem \ref{MainThm2_statement} in view of Propositions \ref{ExtToWeighted} and \ref{BoundFractPart}.
    
    This follows with the same strategy as the previous proof, the only difference being that we bound the first term on the right-hand side of Proposition \ref{lower_weighted_reduction} from below by the contribution of small $h$, i.e. by
    \[ \left( \frac{\log x}{U}\right)^{2q}  \mathbb{E}\left[\left(\frac{1}{\log x}\int_{0}^{U/\log x}\left|F_{\tilde t}(4U/\log x+ih) \right|^2\mathrm{d}h\right)^q\right]. \]
    Using H\"older's inequality as a lower bound to pull the $q$-th power inside the integral, this is in turn
    \[ \gg U^{-q} \frac{\log x}{U} \int_0^{U/\log x} \mathbb{E}\Big[\left|F_{\tilde t}(4U/\log x+ih) \right|^{2q}\Big] \mathrm{d}h  \gg \frac{(\log x)^{q (2q-1)}}{U^{2q^2}}, \]
    where the last step follows from Lemma \ref{lem:MomentsEulerProducts}. Again applying Proposition \ref{prop:upper_uniform_REP} for the subtracted term in Proposition \ref{lower_weighted_reduction}, the claims follow.
\end{proof}

\subsection{The lower bound at \texorpdfstring{$q=1/2$}{q=1/2}}

By Propositions \ref{lower_extended_reduction} and \ref{lower_weighted_reduction}, both claims would once again follow if we can show that uniformly for sufficiently large $U$, we have      
\begin{equation}\label{jun9.2}
        \mathbb{E}\left[\left(\frac{1}{\log x}\int_{4U/\log x}^{1/2}\left|\frac{F_{\tilde t}(4U/\log x+ih)}{{4U}/{\log x}+ih}\right|^2\mathrm{d}h\right)^{1/2} \right] \gg \frac{(\log \log x)^{3/4}}{\sqrt{U}}
    \end{equation}
    and 
    \begin{equation}\label{jun9.3}
        \mathbb{E}\left[\left(\frac{1}{e^{U}\log x}\int_{-\infty}^{\infty}\left|\frac{F_{\tilde t}(2U/\log x+ih)}{{2U}/{\log x}+ih}\right|^2 \log^2(|h|+2) \mathrm{d}h\right)^{1/2} \right] \ll \frac{(\log \log x)^{3/4}}{ e^{U/2}}.
    \end{equation}
    
    The second claim follows immediately from Proposition \ref{prop:upper_uniform_REP}. 
    To verify \eqref{jun9.2}, we begin by partitioning the range of integration into $e$-adic intervals as was done for the upper bound. With $\sigma = \frac{4U}{\log x}$ and $t_\sigma = [\log \sigma^{-1}]$, we have
\begin{equation}\label{eq:lowerbound_1/2sum}
        \mathbb{E}\left[\left(\frac{1}{\log x}\int_{4U/\log x}^{1/2}\left|\frac{F_{\tilde t}(\sigma+ih)}{\sigma + ih}\right|^2 \mathrm{d}h \right)^{1/2}\right] \gg \mathbb{E}\left[\left( \sum_{k = 1}^{t_\sigma } \frac{e^{2k}}{\log x} \int_{e^{-(k+1)}}^{e^{-k}}|F_{\tilde t}(\sigma+ih)|^2 \mathrm{d}h \right)^{1/2}\right].
\end{equation}
As explained in the introduction, our strategy will then consist of showing that this sum is essentially dominated by its (random) largest term. To control the location and size of this term, we introduce the event 
\begin{equation}\label{eq:E(k,V)def}
    \mathcal{E}(k^*,V) := \mathcal{E}_{\le}(k^*,V) \cap \mathcal{E}_{>}(k^*)
\end{equation}
 with
\begin{align*}
     \mathcal{E}_{\le}(k^*,V) &:= \left \{ |F_{k_0,k}(\sigma)| < e^{V} \; \forall \; k_0 < k < k^* \right \} \cap \left\{ |F_{k_0,k^*}(\sigma)| \in \big[e^{V},e^{V+1}\big)  \right\}, \\
     \mathcal{E}_{>}(k^*) &:= \left \{ |F_{k^*,k}(\sigma)| < 1 \; \forall \; k^* < k \le t_\sigma \right\},
\end{align*}
where $k_0$ denotes a sufficiently large absolute constant and $k^*>k_0$. On $\mathcal{E}(k^*,V)$, the process $(F_{k_0,k}(\sigma))_{k_0<k\leq t_\sigma}$ attains its maximal value at ``time'' $k^*$, where $|F_{k_0,k^*}(\sigma)|\approx e^V$.

The following three propositions will allow us to work with this event, and their proofs are postponed to Section \ref{sec:E_results}. 

\begin{prop}\label{ProbLargeMax} 
    Let $x$ be sufficiently large and let $k_0$ be a sufficiently large fixed integer. Suppose that $\frac{-10}{\log x} \le \sigma \le \frac{\log \log x}{\log x}$. Then uniformly for $\frac{k_0}{3} \le V \le \frac{\log \log x}{10}$, and $|k^* - 2 V | \le \sqrt{V}$, we have
\[ \mathbb{P}(\mathcal{E}(k^*,V)) \asymp \mathbb{P}(\mathcal{E}_{\leq }(k^*,V))\asymp \frac{1}{\sqrt{V}}e^{-V}.
\]

\end{prop}

\begin{prop}\label{prop:condexpecttwists}
Let $x$ be sufficiently large, and let $k_0$ be a sufficiently large fixed integer. Suppose that $\frac{-10}{\log x} \le \sigma \le \frac{\log \log x}{\log x}$. Uniformly for $\frac{k_0}{3} \le V \le \frac{\log \log x}{10}$, for $k^* > k_0$ satisfying $|k^* - 2 V| \le \sqrt{V}$, and $|h_1|, |h_2| \le e^{-k^*}$, we have
    \begin{equation}\label{CondExpectSingleTwist}
        \mathbb{E} \left[ \left| \frac{F_{k_0,k^*}(\sigma+i h_1)}{F_{k_0,k^*}(\sigma)} \right|^2 \; \Bigg|  \; \mathcal{E}_{\le}(k^*,V) \right] \gg 1
    \end{equation}
    and
    \begin{equation}\label{eq:CondExpectDoubleTwist}
        \mathbb{E} \left[ \left| \frac{F_{k_0,k^*}(\sigma+i h_1)}{F_{k_0,k^*}(\sigma)} \right|^2 \left| \frac{F_{k_0,k^*}(\sigma+i h_2)}{F_{k_0,k^*}(\sigma)} \right|^2 \; \Bigg| \; \mathcal{E}_{\le}(k^*,V) \right] \ll 1.
    \end{equation}
\end{prop}
\begin{prop}\label{ModifiedBarrierLB} Let $x$ and $U \le \log \log x$ be sufficiently large. Set $\sigma = \tfrac{4U}{\log x}$ and let $k_0<k^* < t_\sigma$. Let $C_1$ be a sufficiently large constant and define $L_{k^*}$ to be the set of points $h\in [e^{-(k^*+1)}, e^{-k^*}]$ for which
\[  
    -C_1 (k-k^*) < \log |F_{k^*,k}(\sigma+ih)|  \le k - k^* - 2 \log (k-k^*) + C_1, \; \forall k^*< k \le t_\sigma,
\]
 and let $\mathcal{L}_{k^*}(h)$ denote the event that $h\in {L}_{k^*}$. Then uniformly for $e^{-(k^*+1)} \le h \le e^{-k^*}$ and in $x$, 
    \begin{equation}
        \mathbb{E} \left[ \mathbf{1}(\mathcal{L}_{k^*}(h)) \mathbf{1}(\mathcal{E}_{>}(k^*)) |F_{k^*,\tilde t}(\sigma+ih)|^2 \right] \gg  \frac{e^{t_\sigma-k^*}}{\sqrt{t_\sigma- k^*}}.
    \end{equation}
\end{prop}

A heuristic explanation of these results is in order prior to proving \eqref{jun9.2}. The first proposition compares the evolution of the path $k\mapsto \log|F_{k_0,k}(\sigma)|$ to that of a corresponding Gaussian random walk, for which the analogue of $\mathcal{E}(k^*,V)$ satisfies the claimed bounds. The second proposition lends rigour to the fact that with high probability, one expects that $ |F_{k_0,k^*}(\sigma+ih)|\approx |F_{k_0,k^*}(\sigma+ih')|$ when $|h-h'|\ll e^{-k^*}$, as the Euler product typically varies at a scale of $e^{-k^*}$. It also shows that this behaviour persists upon conditioning on $\mathcal{E}_{\leq }(k^*,V)$. In a similar vein, Proposition \ref{ModifiedBarrierLB} shows that for a suitable barrier event $\mathcal{L}_{k^*}(h)$, the bound one expects for $\mathbb{E}[\mathbf{1}(\mathcal{L}_{k^*}(h))|F_{k^*, \tilde t}(\sigma+ih)|^2]$ using a random-walk heuristic persists when restricting to samples in $\mathcal{E}_{>}(k^*)$.

Lastly, we will also need the following second moment bound for the integral of $|F_{k^*,\tilde t}(\sigma+ih)|^2$ over points for which $\mathcal{L}_{k^*}(h)$ holds. This follows straightforwardly from existing results in \cite{Harper1}, and the proof can be found in Appendix \ref{appendix} for completeness. 
\begin{prop}\label{ShortIntSecMomBarr}
    Let $x$ and $U \le \log \log x$ be sufficiently large. Set $\sigma = \frac{4 U}{\log x}$ and $t_\sigma = \left[ \log \sigma^{-1}\right]$, and let $k^* < t_\sigma$ be an integer. Let ${L}_{k^*}$ be as in the statement of Proposition \ref{ModifiedBarrierLB}. Then we have
    \begin{equation}
        \mathbb{E}\left[ \left( e^{2 k^*-t_\sigma} \int_{{L}_{k^*}} |F_{k^*,\tilde t}(\sigma+ih)|^2 \mathrm{d}h \right)^{2} \right] \ll \frac{1}{t_\sigma - k^*}.
    \end{equation}
\end{prop}

We are now equipped to complete the proof of \eqref{jun9.2}.

\begin{proof}[Proof of the lower bound in Theorems \ref{MainThm1_statement} and \ref{MainThm2_statement} for $q= \frac{1}{2}$]
   
We note that under $\mathcal{E}(k^*, V)$, we have
\[\max_{k_0 < k \le t_\sigma} |F_{k_0,k}(\sigma)| = |F_{k_0,k^*}(\sigma)| \asymp e^V.\]
Moreover, we remark that
\[ |F_{k_0}(\sigma+ih)| \ge \prod_{\log \log p \le k_0} \left| 1 - \frac{1}{p^{1/2+\sigma}} \right| \gg 1 \]
is deterministically bounded from below uniformly in $h$ and may thus be discarded.

By our reductions at the beginning of this subsection, it suffices to obtain lower bounds for the right-hand side of \eqref{eq:lowerbound_1/2sum}. Partitioning according to the occurrences of the events $\mathcal{E}(k^*, V)$, that expression is
\begin{align}\label{eq:lowerbound_conditional}
    & \ge\mathbb{E}\left[ \max_{k_0 < k \le t_\sigma } |F_{k_0,k}(\sigma)| \Big( \frac{e^{2k}}{\log x} \frac{1}{|F_{k_0,k}(\sigma)|^2} \int_{e^{-(k+1)}}^{e^{-k}}|F_{k_0,\tilde t}(\sigma+ih)|^2 \mathrm{d}h \Big)^{1/2}\right] \nonumber \\
    &\ge \sum_{\frac{k_0}{3} \le V \le \frac{t}{10}} \sum_{\substack{|k^* - 2 V| \le \sqrt{V} \\ k^* > k_0}} \mathbb{E}\bigg[ \max_{k_0 < k \le t_\sigma } \bigg\{ |F_{k_0,k}(\sigma)| \nonumber \\ &\hspace{100pt} \Big( \frac{e^{2k}}{\log x} \frac{1}{|F_{k_0,k}(\sigma)|^2} \int_{e^{-(k+1)}}^{e^{-k}}|F_{k_0,\tilde t}(\sigma+ih)|^2 \mathrm{d}h \Big)^{1/2} \bigg\} \Big \vert \mathcal{E}(k^*,V) \bigg] \mathbb{P}(\mathcal{E}(k^*,V)) \nonumber \\
    &\ge \sum_{\frac{k_0}{3} \le V \le \frac{t}{10}} e^V \!\!\!\!\!\! \sum_{\substack{|k^* - 2 V| \le \sqrt{V} \\ k^* > k_0}} \mathbb{E}\left[ \Big( \frac{e^{2k^*}}{\log x} \frac{1}{|F_{k_0,k^*}(\sigma)|^2} \int_{e^{-(k^*+1)}}^{e^{-k^*}}|F_{k_0,\tilde t}(\sigma+ih)|^2 \mathrm{d}h \Big)^{1/2} \Big \vert \mathcal{E}(k^*,V) \right] \mathbb{P}(\mathcal{E}(k^*,V)).
\end{align}

Inserting Proposition \ref{ProbLargeMax} into \eqref{eq:lowerbound_conditional} gives that this is in turn
\begin{equation}\label{jun9.1}
\gg \sum_{\frac{k_0}{3} \le V \le \frac{t}{10}} V^{-1/2} \sum_{\substack{|k^* - 2 V| \le \sqrt{V} \\ k^* > k_0}} \mathbb{E}\left[ \Big( \frac{e^{2k^*}}{\log x} \frac{1}{|F_{k_0,k^*}(\sigma)|^2} \int_{e^{-(k^*+1)}}^{e^{-k^*}}|F_{k_0,\tilde t}(\sigma+ih)|^2 \mathrm{d}h \Big)^{1/2} \Big \vert \mathcal{E}(k^*,V) \right].
\end{equation}
For a lower bound, we are free to restrict the integrals over $h$ in each summand to a subset of our choosing, by positivity of the integrand. We choose points for which the evolution of the paths $k\mapsto \log|F_{k^*,k}(\sigma+ih)|$ is constrained to lie under a barrier. Concretely, let $C_1$ and $L_{k^*}$ be as in the statement of Proposition \ref{ModifiedBarrierLB}.
Then by the reverse H\"older inequality, conditional on $\mathcal{E}(k^*,V)$, we have 
\begin{align}
    &\mathbb{E}\left[ \left( \frac{e^{2k^*}}{\log x} \frac{1}{|F_{k_0,k^*}(\sigma)|^2} \int_{e^{-(k^*+1)}}^{e^{-k^*}}|F_{k_0,\tilde t}(\sigma+ih)|^2 \mathrm{d}h \right)^{1/2} \Big \vert \mathcal{E}(k^*,V) \right] \nonumber \\
    &\ge \mathbb{E}\left[ \left( \frac{e^{2k^*}}{\log x} \frac{1}{|F_{k_0,k^*}(\sigma)|^2} \int_{L_{k^*}} |F_{k_0,\tilde t}(\sigma+ih)|^2 \mathrm{d}h \right)^{1/2} \Big \vert \mathcal{E}(k^*,V) \right] \nonumber \\
    &\ge \frac{\mathbb{E}\left[ \frac{e^{2k^*}}{\log x} \frac{1}{|F_{k_0,k^*}(\sigma)|^2} \int_{L_{k^*}} |F_{k_0,\tilde t}(\sigma+ih)|^2 \mathrm{d}h \Big \vert \mathcal{E}(k^*,V) \right]^{3/2}}{\mathbb{E}\left[ \left( \frac{e^{2k^*}}{\log x} \frac{1}{|F_{k_0,k^*}(\sigma)|^2} \int_{L_{k^*}} |F_{k_0,\tilde t}(\sigma+ih)|^2 \mathrm{d}h \right)^{2} \Big \vert \mathcal{E}(k^*,V) \right]^{1/2}}. \label{RevCondHolderBarrier}
\end{align}

The next step will be to obtain an upper bound for the denominator. In order to do so, we use that $\mathcal{E}(k^*,V) \subseteq \mathcal{E}_{\le}(k^*,V)$, which only depends on $X(p)$ for $p \le e^{e^{k^*}}$. Applying Proposition \ref{ProbLargeMax}, we deduce that

\begin{align*}
    &\mathbb{E}\left[ \left( \frac{e^{2k^*}}{\log x} \frac{1}{|F_{k_0,k^*}(\sigma)|^2} \int_{L_{k^*}} |F_{k_0,\tilde t}(\sigma+ih)|^2 \mathrm{d}h \right)^{2} \Big \vert \mathcal{E}(k^*,V) \right] \\
    &\ll \mathbb{E}\left[ \left( \frac{e^{2k^*}}{\log x} \frac{1}{|F_{k_0,k^*}(\sigma)|^2} \int_{L_{k^*}} |F_{k_0,\tilde t}
    (\sigma+ih)|^2 \mathrm{d}h \right)^{2} \Big \vert \mathcal{E}_{\le}(k^*,V) \right] \\
    &= \frac{ e^{4k^*}}{(\log x)^2} \iint \mathbb{E} \bigg[ \mathbf{1} (\mathcal{L}_{k^*}(h_1)) \mathbf{1} (\mathcal{L}_{k^*}(h_2)) \frac{|F_{k_0,\tilde t}(\sigma+ih_1)|^2 |F_{k_0,\tilde t}(\sigma+ih_2)|^2}{|F_{k_0,k^*}(\sigma)|^4} \bigg \vert \mathcal{E}_{\le}(k^*,V) \bigg] \mathrm{d}h_1 \mathrm{d}h_2,
\end{align*}
where the last double integral ranges over $[e^{-(k^*+1)}, e^{-k^*}]^2$. For the inner expectation, we can split everything into quantities that only depend on primes $p \le e^{e^{k^*}}$ and those that only depend on larger primes. This gives

\begin{align*}
    &\mathbb{E} \bigg[ \mathbf{1} (\mathcal{L}_{k^*}(h_1)) \mathbf{1} (\mathcal{L}_{k^*}(h_2)) \frac{|F_{k_0,\tilde t}(\sigma+ih_1)|^2 |F_{k_0,\tilde t}(\sigma+ih_2)|^2}{|F_{k_0,k^*}(\sigma)|^4} \bigg \vert \mathcal{E}_{\le}(k^*,V) \bigg]\\
    &=\mathbb{E} \big[ \mathbf{1} (\mathcal{L}_{k^*}(h_1)) \mathbf{1} (\mathcal{L}_{k^*}(h_2)) |F_{k^*,\tilde t}(\sigma+ih_1)|^2 |F_{k^*,\tilde t}(\sigma+ih_2)|^2 \big] \times \\ &\hspace{150pt} \mathbb{E} \Big[ \frac{|F_{k_0,k^*}(\sigma+ih_1)|^2}{|F_{k_0,k^*}(\sigma)|^{2}} \frac{|F_{k_0,k^*}(\sigma+ih_2)|^2}{|F_{k_0,k^*}(\sigma)|^{2}} \big \vert \mathcal{E}_{\le}(k^*,V) \Big].
\end{align*}
The conditional expectation here is $\ll 1$ by \eqref{eq:CondExpectDoubleTwist}, and thus we obtain that
\begin{align*}
    &\mathbb{E}\left[ \left( \frac{e^{2k^*}}{\log x} \frac{1}{|F_{k_0,k^*}(\sigma)|^2} \int_{{L}_{k^*}} |F_{k_0,\tilde t}(\sigma+ih)|^2 \mathrm{d}h \right)^{2} \Big \vert \mathcal{E}(k^*,V) \right] \\
    \ll \; &\mathbb{E}\left[ \left( \frac{e^{2k^*}}{\log x} \int_{{L}_{k^*}} |F_{k^*,\tilde t}(\sigma+ih)|^2 \mathrm{d}h \right)^{2} \right] \ll \frac{1}{U^2(t_\sigma-k^*)}
\end{align*}
by Proposition \ref{ShortIntSecMomBarr} after noting that $e^{-t_\sigma} \asymp \frac{U}{\log x}$. In view of \eqref{RevCondHolderBarrier}, our next aim is to show that
\begin{equation}
    \frac{e^{2k^*}}{\log x} \mathbb{E}\left[  \frac{1}{|F_{k_0,k^*}(\sigma)|^2} \int_{{L}_{k^*}} |F_{k_0,\tilde t}(\sigma+ih)|^2 \mathrm{d}h \,\Big \vert\,  \mathcal{E}(k^*,V) \right] \gg \frac{1}{U \sqrt{t_\sigma-k^*}}.
\end{equation}
To see this, we follow a similar strategy of separating small and large primes to deduce that
\begin{align*}
    &\mathbb{E}\left[  \frac{1}{|F_{k_0,k^*}(\sigma)|^2} \int_{{L}_{k^*}} |F_{k_0,\tilde t}(\sigma+ih)|^2 \mathrm{d}h \; \Big \vert \; \mathcal{E}(k^*,V) \right] \\
    &= \int_{e^{-(k^*+1)}}^{e^{-k^*}} \mathbb{E} \left[ \mathbf{1}\big(\mathcal{L}_{k^*}(h)\big) \frac{|F_{k_0,\tilde t}(\sigma+ih)|^2}{|F_{k_0,k^*}(\sigma)|^2} \; \Bigg\vert \; \mathcal{E}(k^*,V) \right] \mathrm{d}h \\
    &= \int_{e^{-(k^*+1)}}^{e^{-k^*}} \mathbb{E} \Big[ \mathbf{1}(\mathcal{L}_{k^*}(h) |F_{k^*,\tilde t}(\sigma+ih)|^2 \; \big\vert \; \mathcal{E}_{>}(k^*) \Big] \mathbb{E} \Big[ \frac{|F_{k_0,k^*}(\sigma+ih)|^2}{|F_{k_0,k^*}(\sigma)|^{2}} \; \big \vert \; \mathcal{E}_{\le}(k^*,V) \Big] \mathrm{d}h.
\end{align*}
In view of \eqref{CondExpectSingleTwist} and Proposition \ref{ModifiedBarrierLB}, this quantity is
\[ \gg e^{-k^*} \frac{e^{t_\sigma - k^*}}{\sqrt{t_\sigma-k^*}} \gg \frac{\log x}{e^{2k^*}} \frac{1}{U \sqrt{t_\sigma-k^*}}.  \]
Plugging this bound into \eqref{RevCondHolderBarrier}, we conclude that
\begin{equation}
    \mathbb{E}\left[ \left( \frac{e^{2k^*}}{\log x} \frac{1}{|F_{k_0,k^*}(\sigma)|^2} \int_{e^{-(k^*+1)}}^{e^{-k^*}}|F_{k_0,\tilde t}(\sigma+ih)|^2 \mathrm{d}h \right)^{1/2} \Big \vert \mathcal{E}(k^*,V) \right] \gg \frac{1}{\sqrt{U} (t_\sigma- k^*)^{1/4}}.
\end{equation}
Further inserting this into \eqref{jun9.1} shows that
\begin{align*}
    \mathbb{E}\left[\left(\frac{1}{\log x}\int_{4U/\log x}^{1/2}\left|\frac{F_{\tilde{t}}(\sigma+ih)}{\sigma + ih}\right|^2 \mathrm{d}h \right)^{1/2}\right] &\gg \sum_{\frac{k_0}{3} \le V \le \frac{t}{10}} V^{-1/2} \sum_{\substack{|k^* - 2 V| \le \sqrt{V}\\k^*>k_0}} \frac{1}{\sqrt{U}(t_\sigma-k^*)^{1/4}} \\ 
    &\gg \frac{1}{\sqrt{U} (\log \log x)^{1/4}} \sum_{\frac{k_0}{3} \le V \le \frac{t}{10}} 1 \gg \frac{(\log \log x)^{3/4}}{\sqrt{U}}.
\end{align*}
\end{proof}

\section{Proof of Propositions \ref{ProbLargeMax}-\ref{ModifiedBarrierLB}}\label{sec:E_results}

\subsection{Proof of Proposition \ref{ProbLargeMax}}

    For the first claim, it suffices to show that
    \[
        \mathbb{P}\big(|F_{k^*,k}(\sigma)| < 1 \forall k^* < k \le t_\sigma\big)
    \]
    is bounded below by a constant, uniformly in $k^*$, $t$ and $V$. For any $0<j<t_\sigma-k^*$, we use the following decomposition
    \begin{align*}
                \mathbb{P}\big(|F_{k^*,k}(\sigma)| < 1 \forall k^* < k \le t_\sigma\big)\geq \mathbb{P}\big(|F_{k^*+j,k}(\sigma)| < e^{-\sum _{\ell\leq j}M(\ell)}, \forall k>j+k^* \big) \\\prod_{\ell=1}^j\mathbb{P}\big(|F_{{k^*+\ell-1},k^*+\ell}(\sigma)| \le e^{M(\ell)}\big),
    \end{align*}
    where $M(\ell)=\mathbb{E}[\log|F_{k^*+\ell-1,k^*+\ell}(\sigma)|]$. Uniformly in our range of $k^*$ and $\sigma$, we have
    \begin{equation}\label{eq:lem6.1:mertens}
           -\sum_{\ell=1}^j M(\ell)=\frac{1}{2}\sum_{\exp(e^{k^*})<p\leq \exp(e^{k^*+j})}\!\!\!\!-\log|1-p^{-1-2\sigma}|\geq \frac{e^{-2\sigma e^{t_\sigma}}}{6}\sum_{\exp(e^{k^*})<p\leq \exp(e^{k^*+j})} p^{-1}
    \end{equation}
    and $e^{-2\sigma e^{t_\sigma}}/6 \gg 1$. It follows by Mertens' theorem that for any large constant $C>0$, there exists a $j_0=j_0(C)$ such that $e^{-\sum_{\ell\leq j_0}M(\ell)} > C$. Noting that  $\log|F_{k^*+\ell-1,k^*+\ell}(\sigma)|-M(\ell)$ is symmetric in law for each $1\leq \ell\leq j_0$, this implies that
    \begin{align}
                \mathbb{P}\big(|F_{k^*,k}(\sigma)| < 1 \forall k^* < k \le t_\sigma\big)&\geq 2^{-j_0}\cdot  \mathbb{P}\big(|F_{k^*+j_0,k}(\sigma)| < e^{-\sum _{\ell\leq j_0}M(\ell)}, \forall k>k^*+j_0\big)\nonumber\\
                &\gg_{C}   \mathbb{P}\big(|F_{k^*+j_0,k}(\sigma)| \le C, \forall k>k^*+j_0 \big)\nonumber \\
                &= 1- \mathbb{P}\big( \exists k>k^*+j_0, |F_{k^*+j_0,k}(\sigma)| > C \big).\label{eq:complementprob}
    \end{align}
    By a union bound and Markov's inequality, this probability is bounded above by
    \[
        \leq \sum_{k>k^*+j_0} \mathbb{P}\big(|F_{k^*+j_0,k}(\sigma)| > C\big) \le \sum_{k>k^*+j_0} C^{-1/2}\mathbb{E}[|F_{k^*+j_0,k}(\sigma)|^{1/2}],
    \]
    and by Lemmas \ref{lem:MomentsEulerProducts} and \ref{lem:PNT}, this is
    \begin{align*}
                \leq 2 C^{-1/2}\sum_{k>k^*+j_0} \exp\bigg(\!-\frac{1}{8}\sum_{\exp(e^{k^*+j_0})<p<\exp(e^{k})}\frac{1}{p^{1+2\sigma}}\bigg) \leq 4 C^{-1/2}\sum_{(k-k^*-j_0)>0} e^{-(k-k^*-j_0)/(8e^2)}
    \end{align*}
   for large enough $k_0$. This makes \eqref{eq:complementprob} bounded away from zero for an appropriate choice of $C$.

    \bigskip

    What's left is to show that 
    \begin{equation}\label{eq:WTS}
        \mathbb{P}\big(\mathcal{E}_{\leq}(k^*,V)\big)\asymp \frac{1}{\sqrt{V}}e^{-V}.
    \end{equation}
    We do this by introducing a tilted measure $\tilde{\mathbb{P}}$, defined through
    \begin{align}\label{eq:tilde_tilt}
                \frac{\mathrm{d}\tilde{\mathbb{P}}}{\mathrm{d}\mathbb{P}}=\frac{|F_{k_0,k^*}(\sigma)|^\lambda}{\mathbb{E}\big[|F_{k_0,k^*}(\sigma)|^{\lambda}\big]},
    \end{align}
    where $\lambda =1/2+V/(k^*-k_0)$. Note that $\lambda$ is bounded in our ranges of $V$ and $k^*$.
    Letting $\tilde{\mathbb{E}}=\mathbb{E}_{\tilde{\mathbb{P}}}$, we can write
    \begin{align*}
                \mathbb{P}\big(\mathcal{E}_{\leq}(k^*,V)\big)&=\mathbb{E}\big[|F_{k_0,k^*}(\sigma)|^{\lambda}\big]\cdot\tilde{\mathbb{E}}\big[|F_{k_0,k^*}(\sigma)|^{-\lambda}\mathbf{1}\big(\mathcal{E}_{\leq }(k^*,V)\big)\big]\\
                &\asymp e^{-\lambda V}\mathbb{E}\big[|F_{k_0,k^*}(\sigma)|^{\lambda}\big]\cdot\tilde{\mathbb{P}}\big(\mathcal{E}_{\leq }(k^*,V)\big).
    \end{align*}
   Lemmas \ref{lem:MomentsEulerProducts} and \ref{lem:PNT} then imply that
    \begin{equation}\label{eq:tilted_LD_V}
        e^{-\lambda V}\mathbb{E}\big[|F_{k_0,k^*}(\sigma)|^{\lambda}\big] \asymp_{}  \exp\Big(\frac{\lambda(\lambda-1)(k^*-k_0)}{2}-\lambda V\Big)=\exp\Big(-\frac{V^2}{2(k^*-k_0)}-\frac{(k^*-k_0)}{8}-\frac{V}{2}\Big),        
    \end{equation}
    which is $\asymp e^{-V}$ for $|k^*-2V|\leq \sqrt{V}$. It thus suffices to show that $\tilde{\mathbb{P}}(\mathcal{E}_{\leq }(k^*,V))\asymp V^{-1/2}$.

    The upper bound will follow from a standard application of the Berry-Esseen theorem, noting that
    \begin{equation}\label{eq:Elequpper}
                \tilde{\mathbb{P}}\big(\mathcal{E}_{\leq }(k^*,V)\big)\leq \tilde{\mathbb{P}}\big(\log|F_{k_0,k^*}(\sigma)|\in [V,V+1) \big),
    \end{equation}
    and that under $\tilde{\mathbb{P}}$, $\log|F_{k_0,k^*}(\sigma)|$ remains a sum of independent random variables. The application of this theorem will require estimates for the mean and variance of $\log|F_{k_0,k^*}(\sigma)|$ under $\tilde{\mathbb{P}}$. We show how to compute the former, the latter following by a similar argument.

    Let $z=1/2+\sigma$ for simplicity. Then
    \[
        \tilde{\mathbb{E}}\big[\log|1+X(p)p^{-z}|\big]=\frac{\partial}{\partial\beta}\bigg|_{\beta=\lambda}\log\mathbb{E}\big[|1+X(p)p^{-z}|^\beta\big]=\frac{\partial}{\partial \beta}\bigg|_{\beta=\lambda}\log\big(|1-p^{-z}|^\beta+|1+p^{-z}|^{\beta}\big).
    \]
    The right-hand side equals
    \[
        \frac{|1+p^{-z}|^\lambda\log|1+p^{-z}|+|1-p^{-z}|^\lambda\log|1-p^{-z}|}{|1+p^{-z}|^\lambda+|1-p^{-z}|^\lambda},
    \]
    and we can expand both the numerator and denominator as power series in $p^{-z}$
    (using the generalised binomial theorem for the $|1\pm p^{-z}|^\lambda$). Doing so, and using $(u+1)^{-1}=1-u+O(u^2)$ for $|u|<\frac{1}{2}$ to invert the series in the denominator, we reach
\begin{equation}\label{eq:lem7.1_tiltedmean}
        \tilde{\mathbb{E}}\big[\log|1+X(p)p^{-z}|\big]=\frac{\lambda-\tfrac{1}{2}}{p^{2z}}+O\big(p^{-4z}\big)=\frac{V}{(k^*-k_0)}\frac{1}{p^{2z}}+O\big(p^{-4z}\big).
\end{equation}
    An analogous, albeit more tedious calculation shows that
    \begin{align}\label{eq:lem7.1_tiltedvariance}  \text{Var}_{\tilde{\mathbb{{P}}}}\big(\log|1+X(p)p^{-z}|\big)=\frac{\partial^2}{\partial\beta^2}\bigg|_{\beta=\lambda}\log \mathbb{E}\big[|1+X(p)p^{-z}|^\beta\big]=p^{-2z}+O\big(p^{-4z}\big).
    \end{align}
    By summing over $\log\log p\in (k_0,k^*]$, using Lemma \ref{lem:PNT} and applying the estimate $\mathrm{Ei}(-2\sigma e^{k^*})-\mathrm{Ei}(-2\sigma e^{k_0})=(k^*-k_0)+o(1)$ for $k^*\leq t/5$, it follows that 
    \begin{align}\label{eq:tilted_mean_uptok^*}
    &\tilde{\mathbb{E}}\big[\log|F_{k_0,k^*}(\sigma)|\big]=V+O\big(e^{-ce^{k_0/2}}\big),\\&\text{Var}_{\tilde{\mathbb{P}}}\big(\log|F_{k_0,k^*}(\sigma)|\big)= \sum_{\log\log p\in (k_0,k^*]}\text{Var}_{\tilde{\mathbb{P}}}\big(\log|1+X(p)p^{-z}|\big)=(k^*-k_0)+O(1).
    \end{align}
    We also have 
    \begin{align*}       \sum_{p}\tilde{\mathbb{E}}\Big[\big|\log|1+X(p)p^{-z}|-\tilde{\mathbb{E}}[\log|1+X(p)p^{-z}|]\big|^3\Big] \ll\sum_p p^{-3z} <\infty,
    \end{align*}
    which follows from the fact that $|\log|1+X(p)p^{-z}||\ll p^{-z} $ uniformly in $p$ for large enough $x$ (since $X(p)$ is still supported on $\{-1,1\}$ under $\tilde{\mathbb{P}}$, $|X(p)|\leq 1$ deterministically).
    Using Lemma \ref{lem:BerryEsseen}, we can therefore conclude that
    \[
        \tilde{\mathbb{P}}\big(\log|F_{k_0,k^*}(\sigma)|\in [V,V+1) \big) \ll \frac{1}{\sqrt{k^*-k_0}}\ll\frac{1}{\sqrt{V}}.
    \]

    For the lower bound, we compare the evolution of the sum $k\mapsto \log|F_{k_0,k}(\sigma)|$ (under the measure $\tilde{\mathbb{P}}$) to that of a Gaussian random walk. For every $k>0$, let $\Delta_k=(1/16)k^{-4}$, so that $\sum_k \Delta_k<1/8$. (We could have picked any $\varepsilon>0$ small enough instead of $1/8$, but choose the latter for concreteness.) Let $\mathcal{T}$ be the set of tuples $(u_k)_{k_0<k\leq k^*}\subseteq \mathbb{R}^{k^*-k_0}$ satisfying $u_k\in \Delta_k\mathbb{Z}$ for every $k$. We will decompose the event $\mathcal{E}_{\leq}(k^*,V)$ into a disjoint union, over admissible tuples $(u_k)_{k_0<k\leq k^*}$, of events of the form $\cap_{k_0<k\leq k^*}\{Y_k\in[u_k,u_k+\Delta_k)\}$, where $Y_k:=\log|F_{k-1,k}(\sigma)|-\tilde{\mathbb{E}}[\log|F_{k-1,k}(\sigma)|]$ is the $k$-th ``increment'' of $\log|F_{k_0,k^*}(\sigma)|-\tilde{\mathbb{E}}[\log|F_{k_0,k^*}(\sigma)|]$. The probability of these events will then be estimated using the Berry-Esseen theorem, noting that each $Y_k$ is itself a sum of a large number of independent random variables.

    Instead of dealing with $\mathcal{E}_{\leq}(k^*,V)$ directly, it will be helpful to work on a more restrictive event $\mathcal{E}_{\leq }(k^*,V)\cap \mathcal{G}$, where
    \[
        \mathcal{G}=\bigg\{\Big|\sum_{k_0<j\leq k}Y_j\Big|<10k, \forall k_0<k<k^*\bigg\}.
    \]
    For a given sequence $(u_k)_{k<k^*}\in \mathcal{T}$, note that
    \[
        \sum_{k_0< j\leq k}u_j\leq \sum_{k_0< j\leq k}Y_j\leq \sum_{k_0< j\leq k}u_j+\sum_{k_0< j\leq k}\Delta_j\leq \sum_{k_0< j\leq k}u_j+1/8,\quad \text{for every $k_0<k\leq k^*$,}
    \]
    on the event $\cap_{k_0<k\leq k^*}\{Y_k\in[u_k,u_k+\Delta_k)\}$. One can then check that if $\mathcal{T}'\subseteq \mathcal{T}$ is the subset of tuples $(u_k)_{k_0<k\leq k^*}$ for which
    \begin{align}
        \forall k_0<k< k^*, -10k<\sum_{k_0<j\leq k}u_j<\min\big(10k,V-M_k\big)-1/8, \\
        \text{ and }\sum_{k_0<j\leq k^*}u_j\in [V-M_{k^*},V-M_{k^*}+7/8],
    \end{align}
    where $M_k:=\sum_{k_0<j\leq k}\tilde{\mathbb{E}}[\log|F_{j-1,j}(\sigma)|]$,
    then
    \[
        \bigsqcup_{(u_k)_{k_0<k\leq k^*}\in \mathcal{T}'} \{Y_k\in [u_k,u_k+\Delta_k)\} \subseteq \mathcal{E}_{\leq}(k^*,V)\cap \mathcal{G}.
    \]
    This union being disjoint, it follows that
    \[
        \tilde{\mathbb{P}}\big(\mathcal{E}_{\leq}(k^*,V)\cap \mathcal{G}\big)\geq \sum_{(u_k)_{k_0<k\leq k^*}\in \mathcal{T}'} \prod_{k_0<k\leq k^*}\tilde{\mathbb{P}}\big(Y_k\in [u_k,u_k+\Delta_k)\big).
    \]
    We now approximate the probabilities on the right-hand side by a natural Gaussian analogue. For each $k$, let $\mathfrak{v}_k:= \tilde{\mathbb{E}}[Y_k^2]$, noting that $\mathfrak{v}_k\gg1 $ uniformly in $k$ by the estimate in \eqref{eq:lem7.1_tiltedvariance} and Lemma \ref{lem:PNT}. There exists an absolute constant $C'>0$ such that
    \begin{align*}
                \sum_{\exp(e^{k-1})<p\leq \exp(e^k)}\tilde{\mathbb{E}}\Big[\big|\log|1+X(p)p^{-z}|-\tilde{\mathbb{E}}[\log|1+X(p)p^{-z}|]\big|^3\Big]&\leq C'\sum_{\exp(e^{k-1})<p
                \leq \exp(e^k)} p^{-3z} \\&\leq 2 C' e^{-ce^{k/2}}
    \end{align*}
    for some $c>0$ and every $k_0<k\leq k^*$, where $z=1/2+\sigma$ as above.
    By Lemma \ref{lem:BerryEsseen}, it thus follows that uniformly in $k$,
    \[
        \tilde{\mathbb{P}}\big(Y_k\in [u_k,u_k+\Delta_k)\big)={\mathbb{P}}\big(N_k\in [u_k,u_k+\Delta_k)\big)+O\big(e^{-ce^{k/2}}\big),
    \]
    where $(N_k)_{k_0<k\leq k^*}$ is a sequence of centred Gaussian random variables with variances $\mathfrak{v}_k$, and the implicit constant does not depend on $k$.
    Using the fact that $|u_k|\leq 20k$ on $\mathcal{G}$, this can be rewritten as
    \[
        \mathbb{P}\big(N_k\in [u_k,u_k+\Delta_k)\big)\cdot\big(1+O(k^{-3})\big),
    \]
    for an implicit constant which is again absolute. Since $\prod_{k > k_0} (1+O(k^{-3}))\asymp 1$ by taking $k_0$ sufficiently large, it follows that
    \[
         \tilde{\mathbb{P}}\big(\mathcal{E}_{\leq}(k^*,V)\cap \mathcal{G}\big) \gg \sum_{(u_k)_{k_0<k\leq k^*}\in \mathcal{T}'} \prod_{k_0<k\leq k^*}\mathbb{P}\big(N_k\in [u_k,u_k+\Delta_k)\big).
    \]
    By taking $k_0$ to be sufficiently large, we can also assume that $|V-M_{k^*}|\leq 1/8$ (see Equation \eqref{eq:tilted_mean_uptok^*}). Letting $\bar{S}_k:=\sum_{k_0<j\leq k}N_j$, we then have the inclusion
    \begin{align*}
                &\big\{-9k<\bar{S}_k<\min(9k,V-M_k-1/8), \forall k_0<k<k^*, \bar{S}_{k^*}\in [1/4,3/4)\big\}\\
                &\subseteq \bigcup_{(u_k)_{k\leq k^*}\in \mathcal{T}'} \big\{N_k\in [u_k,u_k+\Delta_k)\big\},
    \end{align*}
    and we may conclude that $\tilde{\mathbb{P}}\big(\mathcal{E}_{\leq}(k^*,V)\cap \mathcal{G}\big)$ is 
    \begin{equation}\label{eq:GaussianE}
                 \gg \mathbb{P}\big(-9k<\bar{S}_k<\min(9k,V-M_k-1/8), \forall k_0< k<k^*, \bar{S}_{k^*}\in [1/4,3/4)\big).
    \end{equation}

    From here, the desired bound is a consequence of the fact that the endpoint of a Gaussian random walk is independent from its ``bridge''. Write $S_k=\bar{S}_k+M_k$ for every $k$, and denote the variance of $\bar{S}_k$ by $\Sigma_{k}=\sum_{k_0<j\leq k}\mathfrak{v}_j$. Letting $B_k:=\bar{S}_k-(\Sigma_k/\Sigma_{k^*})\bar{S}_{k^*}$, we then have that
    \[
        \mathbb{E}\big[B_k\bar{S}_{k^*}\big]=\mathbb{E}\big[(\bar{S}_k-(\Sigma_k/\Sigma_{k^*})\bar{S}_{k^*})\bar{S}_{k^*}\big]=0,\quad  \forall k_0<k<k^*,
    \]
    and a covariance calculation shows that $B_k\sim \mathcal{N}(0,\Sigma_{k}-\Sigma_k^2/\Sigma_{k^*})$.
    
    By conditioning on the value of $\bar{S}_{k^*}$, it follows that the right-hand side in \eqref{eq:GaussianE} is
    \begin{align*}
           &\gg \int_{1/4}^{3/4}\frac{e^{-u^2/(2\Sigma_{k^*})}}{\sqrt{\Sigma_{k^*}}} \mathbb{P}\big(-9k<B_k+u(\Sigma_k/\Sigma_{k^*})<\min(9k,V-M_k-1/8), \forall k_0<k<k^*\big)\mathrm{d}u
           \\
                   &\gg \frac{1}{\sqrt{\Sigma_{k^*}}}\mathbb{P}\big(-9k<B_k<\min(8k,V-M_k-7/8), \forall k_0<k<k^*\big).
    \end{align*}

    Since $\Sigma_{k^*}=(k^*-k_0)+O(1)\ll V$, what's left is to show that this probability is $\gg 1$. Noting that the variance of $B_{k^*-1}$ is $O(1)$, this probability is
\begin{align}\label{eq:finalprob_goodevent}
    \gg \int_{D}^{D+1}\mathbb{P}\big(-9k<B_k<\min(8k,V-M_k-7/8), \forall k_0<k<k^*\big|B_{k^*-1}=-C\big)\mathrm{d}C
\end{align}
for any fixed large constant $D>0$. We may in particular assume that $-D<V-M_{k^*-1}-7/8$ in what follows, so that the condition at $k=k^*-1$ is trivially satisfied in the above for all $C$. For smaller $k$, conditioning on $B_{k^*-1}=-C$ causes $B_k$ to acquire the following distribution
\[
    B_k\sim \mathcal{N}\bigg(-C\frac{\Sigma_k}{\Sigma_{k^*-1}},\Sigma_k\Big(1-\frac{\Sigma_k}{\Sigma_{k^*-1}}\Big)\bigg),\quad \forall k_0<k<k^*-1.
\]
With this law now explicit, we now show how to bound the probability in \eqref{eq:finalprob_goodevent} below uniformly in $C\in [D,D+1]$. We will do so by bounding its complement using a union bound over $k$, and estimating the resulting terms using standard Gaussian tail bounds. 

For large enough $k_0$, 
\begin{align}\label{eq:lem6.1_part1}
        \sum_{k_0<k< k^*-1}\mathbb{P}\big(-B_k>9k\,\big|\,B_{k^*-1}=-C\big) &\leq \sum_{k_0<k<k^*-1} \exp\bigg(-\frac{1}{2}\frac{(9k-C\Sigma_k/\Sigma_{k^*-1})^2}{\Sigma_k(1-\Sigma_k/\Sigma_{k^*-1})}\bigg)\\
        &\leq \sum_{k_0<k<k^*-1} \exp\bigg(-\frac{k}{2}\frac{k}{\Sigma_k}\bigg)
\end{align}
By \eqref{eq:lem7.1_tiltedvariance} and Lemma \ref{lem:PNT}, there exists an absolute constant $c>0$ such that $k/\Sigma_k>c$ uniformly in $k_0<k<k^*$. The right-hand side is therefore $\leq\sum_{k>k_0}e^{-(c/2)k}$, which is $<1/10$, say, for large enough $k_0$. Similarly, 
\begin{align}\label{eq:lem6.1_part2}
    \sum_{k_0<k<k^*-1} \mathbb{P}\big(B_k>8k\,\big|\,B_{k^*-1}=-C\big) \leq \sum_{k_0<k<k^*-1} \exp\bigg(-\frac{1}{2}\frac{k}{\Sigma_k}(64k)\bigg)<\frac{1}{10}.
\end{align}
Lastly, we bound the probability that $B_k$ exceeds $(V-M_k-7/8)$ at some $k$. This is where the assumption that $B_{k^*-1}=-C$ is used; without it, the terms arising from $k$ near $k^*$ could give a contribution exceeding $1$.

Note that by \eqref{eq:lem7.1_tiltedmean},
\[
    V-M_k=V-\frac{V}{k^*-k_0}\sum_{k_0<\log\log p\leq k}p^{-2z} +O\Big(\sum_{k_0<\log\log p\leq k}p^{-4z}\Big).
\]
By making $k_0$ large enough and using the restriction $|k^*-2V|\leq \sqrt{V}$, one can use Lemma \ref{lem:PNT} to show that $V-M_{k}-7/8 \geq A k^*(k^*-k)/(k^*-k_0)$, for some absolute constant $A>0$. It follows that
\begin{align}
    &\sum_{k_0<k<k^*-1}\mathbb{P}\bigg(B_k>(V-M_k-7/8)\,\bigg|\,B_{k^*-1}=-C\bigg)\\
    &\leq \sum_{k_0<k<k^*-1} \exp\bigg(-\frac{1}{2}\frac{\Sigma_{k^*-1}}{\Sigma_{k}(\Sigma_{k^*-1}-\Sigma_{k})}\Big(A\frac{k^*(k^*-k)}{k^*-k_0}+C\frac{\Sigma_k}{\Sigma_{k^*-1}
    }\Big)^2\bigg)\\
    &\leq \sum_{k_0<k<k^*-1}\exp\bigg(-AC \frac{(k^*-k)k^*}{(k^*-k_0)(\Sigma_{k^*-1}-\Sigma_k)}\bigg)\exp\bigg(-\frac{A^2}{2} \frac{\Sigma_{k^*-1}(k^*-k)^2{k^*}^2}{\Sigma_{k}(\Sigma_{k^*-1}-\Sigma_k)(k^*-k_0)^2}\bigg).\label{eq:6.1_final}
\end{align}
Lemma \ref{lem:PNT} also implies that
\begin{equation}\label{eq:CapitalSigmaAsymptotic}
    \Sigma_k=(k-k_0)+O\big(e^{-ce^{k_0/2}}\big)
\end{equation}
uniformly in $k_0<k<k^*-1$. By taking $x, k_0$ large enough, this yields the estimates
\[
    \frac{(k^*-k)k^*}{(k^*-k_0)(\Sigma_{k^*-1}-\Sigma_k)}\geq 1
\]
on $k_0<k<k^*-1$, and 
\begin{align*}
    \frac{\Sigma_{k^*-1}(k^*-k)^2{k^*}^2}{\Sigma_{k}(\Sigma_{k^*-1}-\Sigma_k)(k^*-k_0)^2}&\geq \frac{(k^*-k_0-1.1)(m+1)}{(k^*-k_0-0.9-m)} \geq \frac{1}{2}(m+1)
\end{align*}
where $m=k^*-1-k>0$, whenever $k^*-k_0\geq 3$. The sum in \eqref{eq:6.1_final} is therefore bounded above by $\leq e^{-AC} \sum_{m\geq 1} e^{-(A^2/4)(m+1)}$,
which is $\leq 1/10$ for large enough $D$ (and thus $C$). Combining this with the bound on \eqref{eq:lem6.1_part1} and \eqref{eq:lem6.1_part2}, we conclude that 
\[  
    \mathbb{P}\Big(-9k<B_k<\min(8k,V-M_k-7/8), \forall k_0<k<k^*\,\Big|\, B_{k^*-1}=-C\Big)\geq 1-\frac{3}{10}>0.
\]
uniformly over $C\in [D,D+1]$.
\qed

\subsection{Proof of Proposition \ref{prop:condexpecttwists}}

We begin with the second claim. Letting $\mu$ be the measure defined through
\begin{equation}\label{eq:mu_change}
       \frac{\mathrm{d}\mu}{\mathrm{d}\mathbb{P}}=\frac{|F_{k_0,k^*}(\sigma+ih_1)F_{k_0,k^*}(\sigma)^{-1}|^2|F_{k_0,k^*}(\sigma+ih_2)F_{k_0,k^*}(\sigma)^{-1}|^2}{\mathbb{E}[|F_{k_0,k^*}(\sigma+ih_1)F_{k_0,k^*}(\sigma)^{-1}|^2|F_{k_0,k^*}(\sigma+ih_2)F_{k_0,k^*}(\sigma)^{-1}|^2]},
\end{equation}
we deduce from Proposition \ref{ProbLargeMax} that the left-hand side in \eqref{eq:CondExpectDoubleTwist} is
\[
    \asymp e^V\sqrt{V}\cdot{ \mu\big(\mathcal{E}_{\leq }(k^*,V)\big)}\cdot  \mathbb{E} \left[ \left| \frac{F_{k_0,k^*}(\sigma+i h_1)}{F_{k_0,k^*}(\sigma)} \right|^2 \left| \frac{F_{k_0,k^*}(\sigma+i h_2)}{F_{k_0,k^*}(\sigma)} \right|^2\right].
\]
Setting $z=1/2+\sigma, z_1=1/2+\sigma+ih_1$, and $ z_2=1/2+\sigma+ih_2$, we have
\begin{align*}
&\mathbb{E} \left[ \left| \frac{F_{k_0,k^*}(\sigma+i h_1)}{F_{k_0,k^*}(\sigma)} \right|^2 \left| \frac{F_{k_0,k^*}(\sigma+i h_2)}{F_{k_0,k^*}(\sigma)} \right|^2\right]\\&=
    \prod_{\exp(e^{k_0})<p\leq \exp(e^{k^*})}\mathbb{E}\bigg[\bigg|\frac{(1+X(p)p^{-z_1})(1+X(p)p^{-z_2})}{(1+X(p)p^{-z})(1+X(p)p^{-z})}\bigg|^2\bigg]\\
    &=\prod_{\exp(e^{k_0})<p\leq \exp(e^{k^*})}\frac{1}{2}\bigg(\bigg|\frac{(1+p^{-z_1})(1+p^{-z_2})}{(1+p^{-z})^2}\bigg|^2+\bigg|\frac{(1-p^{-z_1})(1-p^{-z_2})}{(1-p^{-z})^2}\bigg|^2\bigg)\\
    &=:\prod_{\exp(e^{k_0})<p\leq \exp(e^{k^*})} A_{\sigma,h_1,h_2}(p).
\end{align*}
Using the identity  $|1\pm p^{-z_j}|^2=(1\pm p^{-z})^2\mp 4p^{-z}\sin^2(|h_j|(\log p)/2)$, we can write
\begin{align}\label{eq:tedious_identity}
   \bigg|\frac{(1\pm p^{-z_1})(1\pm p^{-z_2})}{(1\pm p^{-z})^2}\bigg|^2&= \Big(1\mp 4\frac{p^{-z}}{(1\pm p^{-z})^2}S_1\Big)\Big(1\mp 4\frac{p^{-z}}{(1\pm p^{-z})^2}S_2\Big)
\end{align}
for $S_j= \sin^2(|h_j|(\log p)/2)$, which implies that
\begin{align*}
        A_{\sigma,h_1,h_2}(p)&=1+2p^{-z}(S_1+S_2)\Big(\frac{1}{(1-p^{-z})^2}-\frac{1}{(1+p^{-z})^2}\Big)+(S_1S_2)\Big(\frac{8p^{-2z}}{(1+p^{-z})^4}+\frac{8p^{-2z}}{(1-p^{-z})^4}\Big)\\
        &\leq 1+\frac{8p^{-2z}}{(1-p^{-2z})^2}(S_1+S_2)+\frac{8p^{-2z}}{(1-p^{-z})^4}(S_1+S_2),
\end{align*}
    since $(S_1S_2)\leq \tfrac{1}{2}(S_1+S_2)$.  It follows that for an absolute constant $c>0$,
\[
    \log\big(A_{\sigma,h_1,h_2}(p)\big)\leq \log \big(1+cp^{-2z}(S_1+S_2)\big)\leq cp^{-2z}\big({|h_1|^2(\log p)^2+|h_2|^2(\log p)^2}\big).
\]
Since $p^{2\sigma}\asymp 1$ uniformly in $\exp(e^{k_0})<p\leq \exp(e^{k^*})$ on our range of $k^*$ and $\sigma$, this gives
\begin{align*}
        \log\bigg(\prod_{\exp(e^{k_0})<p\leq \exp(e^{k^*})} A_{\sigma,h_1,h_2}(p)\bigg) \ll &(|h_1|^2+|h_2|^2) \sum_{\exp(e^{k_0})<p\leq \exp(e^{k^*})} \frac{(\log p)^2}{p^{2z}} \\
        \asymp &(|h_1|^2+|h_2|^2) \sum_{\exp(e^{k_0})<p\leq \exp(e^{k^*})} \frac{(\log p)^2}{p}. 
\end{align*}
By Lemma \ref{lem:log_PNT}, we conclude that 
\begin{equation}\label{eq:2pointratio}
    \mathbb{E} \left[ \left| \frac{F_{k_0,k^*}(\sigma+i h_1)}{F_{k_0,k^*}(\sigma)} \right|^2 \left| \frac{F_{k_0,k^*}(\sigma+i h_2)}{F_{k_0,k^*}(\sigma)} \right|^2\right]\asymp  1
\end{equation}
for $|h_1|,|h_2|<e^{-k^*}$, and it therefore suffices to prove that
\begin{align}\label{eq:mu_LD}
        \mu\big(|F_{k_0,k^*}(\sigma)|\in\big[e^V,e^{V+1}\big)\big)\ll \frac{1}{\sqrt{V}} e^{-V}
\end{align}
in order to establish \eqref{eq:CondExpectDoubleTwist}.

We will do so by tilting $\mu$ in order to make $\log|F_{k_0,k^*}(\sigma)|$ typically $\approx V$, similarly to the proof of \eqref{eq:Elequpper}. This will require the following estimate, which we prove first:
\begin{align}\label{eq:huge_numerator_nu}
       \mathbb{E}\bigg[|F_{k_0,k^*}(\sigma)|^{\lambda}\bigg|\frac{F_{k_0,k^*}(\sigma+ih_1)F_{k_0,k^*}(\sigma+ih_2)}{F_{k_0,k^*}(\sigma)^2}\bigg|^2\bigg] \asymp \mathbb{E}[|F_{k_0,k^*}(\sigma)|^\lambda].
\end{align}
Using the identity in \eqref{eq:tedious_identity} and the same notation as above, the per-prime contribution to the left-hand side in \eqref{eq:huge_numerator_nu} equals
\begin{align*}
\frac{|1+p^{-z}|^{\lambda}}{2}\bigg(1-4\frac{p^{-z}}{(1+p^{-z})^2}(S_1+S_2)+16\frac{p^{-2z}}{(1+p^{-z})^4}(S_1S_2)\bigg)\\
    +\frac{|1-p^{-z}|^{\lambda}}{2}\bigg(1+4\frac{p^{-z}}{(1-p^{-z})^2}(S_1+S_2)+16\frac{p^{-2z}}{(1-p^{-z})^4}(S_1S_2)\bigg),
\end{align*}
which can be rearranged into
\begin{align}\label{eq:recover_untiltedexp}
    \frac{1}{2}\big(|1+p^{-z}|^{\lambda}+|1-p^{-z}|^\lambda\big)\big(1+E_p(S_1,S_2)\big)=\mathbb{E}[|1+X(p)p^{-z}|^\lambda]\big(1+E_p(S_1,S_2)\big)
\end{align}
for $ E_p(S_1,S_2)$ equalling
\begin{align*}
\frac{4p^{-z}\big((1-p^{-z})^{\lambda-2}-(1+p^{-z})^{\lambda-2}\big)(S_1+S_2)+16p^{-2z}\big((1+p^{-z})^{\lambda-4}+(1-p^{-z})^{\lambda-4}\big)S_1S_2}{(1+p^{-z})^{\lambda}+(1-p^{-z})^{\lambda}}.
\end{align*}
Note that by our assumptions on $V$ and $k^*$, $\lambda \ll k_0$, which we recall is fixed. Since $|p^{-z}|<1$, the generalised binomial theorem can thus be applied to $(1\pm p^{-z})^{\lambda-2}$ and $(1\pm p^{-z})^{\lambda-4}$ to get 
\begin{align*}
    &\big((1-p^{-z})^{\lambda-2}-(1+p^{-z})^{\lambda-2}\big)=-2(\lambda-2)p^{-z}+O(k_0^3p^{-3z}),\\
    &\big((1+p^{-z})^{\lambda-4}+(1-p^{-z})^{\lambda-4}\big)=2+O(k_0^2p^{-2z}),
\end{align*}
 Inserting these estimates into the definition of $E_p(S_1,S_2)$ thus gives
\begin{equation}\label{eq:error_pestimate}
    |E_p(S_1,S_2)|\leq c'k_0p^{-2z}(S_1+S_2) \leq c'k_0p^{-2z}\big(|h_1|^2(\log p)^2+|h_2|^2(\log p)^2\big),
\end{equation}
for $c'>0$ an absolute constant, which implies that $|E_p(S_1,S_2)|\leq 1/2$ for large enough $k_0$ by the double-exponential decay of $p^{-2z}$ in $k_0$. Lemma \ref{lem:log_PNT} and the estimate  $|\log(1+E_p(S_1,S_2))|\ll |E_p(S_1,S_2)|$ thus give
\[
        \Big|\sum_{\exp(e^{k_0})<p\leq \exp(e^{k^*})}\!\!\log\big(1+E_p(S_1,S_2)\big)\Big|\ll \sum_{\exp(e^{k_0})<p\leq \exp(e^{k^*})} \frac{(|h_1|^2+|h_2|^2)(\log p)^2}{p^{2z}} \ll 1,
\]
and combining this with \eqref{eq:recover_untiltedexp}, we obtain
\begin{align*}
        &\mathbb{E}\bigg[|F_{k_0,k^*}(\sigma)|^{\lambda}\bigg|\frac{F_{k_0,k^*}(\sigma+ih_1)F_{k_0,k^*}(\sigma+ih_2)}{F_{k_0,k^*}(\sigma)^2}\bigg|^2\bigg]  \\
        &= \mathbb{E}[|F_{k_0,k^*}(\sigma)|^\lambda]\exp\Big(\sum_{\exp(e^{k_0})<p\leq \exp(e^{k^*})}\!\!\log\big(1+E_p(S_1,S_2)\big)\Big) 
        \\&\asymp \mathbb{E}[|F_{k_0,k^*}(\sigma)|^\lambda].
\end{align*}

We may now proceed with the proof of \eqref{eq:mu_LD}. Let $\nu$ be the measure defined by
\[
    \frac{\mathrm{d}\nu}{\mathrm{d}\mu}=\frac{|F_{k_0,k^*}(\sigma)|^\lambda}{\mathbb{E}_\mu[|F_{k_0,k^*}(\sigma)|^\lambda]}
\]
for $\lambda=1/2+V/(k^*-k_0)$, noting that $\mathbb{E}_\mu [|F_{k_0,k^*}(\sigma)|^\lambda]\asymp \mathbb{E}[|F_{k_0,k^*}(\sigma)|^\lambda]$ by \eqref{eq:2pointratio} and \eqref{eq:huge_numerator_nu}. For this choice of $\lambda$, recall that we have shown in  \eqref{eq:tilted_LD_V} that $e^{-\lambda V} \mathbb{E}\big[|F_{k_0,k^*}(\sigma)|^\lambda\big]\asymp e^{-V}$. It follows that
\begin{align*}
    \mu\big(|F_{k_0,k^*}(\sigma)|\in\big[e^V,e^{V+1}\big)\big)&=\mathbb{E}_\mu\big[|F_{k_0,k^*}(\sigma)|^\lambda \big]\cdot \mathbb{E}_\nu\big[|F_{k_0,k^*}(\sigma)|^{-\lambda}\mathbf{1}\big(|F_{k_0,k^*}(\sigma)|\in\big[e^V,e^{V+1}\big)\big)\big]\\
    &\asymp e^{-\lambda V} \mathbb{E}\big[|F_{k_0,k^*}(\sigma)|^\lambda\big]\cdot    \nu\big(|F_{k_0,k^*}(\sigma)|\in\big[e^V,e^{V+1}\big)\big)
    \\ &\asymp e^{-V} \nu\big(|F_{k_0,k^*}(\sigma)|\in\big[e^V,e^{V+1}\big)\big),
\end{align*}
and must show that the remaining probability is $\ll V^{-1/2}$. By Lemma \ref{lem:BerryEsseen}, this would follow directly from the following estimates:
\begin{align}
    &\sum_{\exp(e^{k_0})<p\leq \exp(e^{k^*})} \mathbb{E}_\nu\Big[\big|\log|1+X(p)p^{-z}|-\mathbb{E}_\nu\big[\log|1+X(p)p^{-z}|\big]\big|^3\Big] \ll \sum_p p^{-3z} <\infty,\\
    &\sum_{\exp(e^{k_0})<p\leq \exp(e^{k^*})}\mathrm{Var}_{\nu}\big(\log|1+X(p)p^{-z}|\big) \asymp V \label{eq:nu_varianceestimate}.
\end{align}
The first follows from the fact that $\nu$ is absolutely continuous with respect to $\mathbb{P}$, as this implies that $|\log|1+X(p)p^{-z}||\leq Cp^{-z}$ for some absolute constant $C>0$. For the second, we compute the variances directly. By definition of $\nu$, 
\begin{align*}
    &\mathbb{E}_\nu\big[(\log|1+X(p)p^{-z}|)^2\big]= \nu\big(X(p)=1\big)\big(\log|1+p^{-z}|\big)^2+\nu\big(X(p)=-1\big)\big(\log|1-p^{-z}|\big)^2\\
    &\mathbb{E}_\nu\big[\log|1+X(p)p^{-z}|\big]^2=\Big(\nu\big(X(p)=1\big)\log|1+p^{-z}|+\nu\big(X(p)=-1\big)\log|1-p^{-z}|\Big)^2,
\end{align*}
which gives
\begin{align*}
        &\mathrm{Var}_{\nu}\big(\log|1+X(p)p^{-z}|\big)= \nu\big(X(p)=1\big)\nu\big(X(p)=-1\big)\Big(\log\Big|\frac{1+p^{-z}}{1-p^{-z}}\Big|\Big)^2.
\end{align*}
These probabilities can then be expressed in terms of the quantities $S_1,S_2$ and $E_p$ defined above to get
\begin{align*}
    \nu\big(X(p)=\pm 1\big)&=\mathbb{E}\Big[\mathbf{1}(X(p)=\pm 1)\frac{\mathrm{d}\mu}{\mathrm{d}\mathbb{P}}\frac{\mathrm{d}\nu}{\mathrm{d}\mu}\Big]\\
    &=\frac{1}{2}\frac{|1\pm p^{-z}|^\lambda}{\mathbb{E}\big[|1+X(p) p^{-z}|^\lambda\big]}\frac{1}{1+E_p(S_1,S_2)}\Big(1\mp 4\frac{p^{-z}}{(1\pm p^{-z})^2}S_1\Big)\Big(1\mp 4\frac{p^{-z}}{(1\pm p^{-z})^2}S_2\Big),
\end{align*}
and inserting this into the expression  for $\mathrm{Var}_{\nu}\big(\log|1+X(p)p^{-z}|\big)$ gives
\begin{equation}\label{eq:variance_comparison}
        \mathrm{Var}_{\nu}\big(\log|1+X(p)p^{-z}|\big)=\mathrm{Var}_{\tilde{\mathbb{P}}}\big(\log|1+X(p)p^{-z}|\big)\frac{\prod_{j=1,2}\Big(1+\frac{16p^{-2z}}{(1-p^{-2z})^2}S_j(1-S_j)\Big)}{(1+E_p(S_1,S_2))^2}
\end{equation}
for
\[
    \frac{\mathrm{d}\tilde{\mathbb{P}}}{\mathrm{d}\mathbb{P}}=\frac{|F_{k_0,k^*}(\sigma)|^\lambda}{\mathbb{E}\big[|F_{k_0,k^*}(\sigma)|^\lambda\big]}.
\]
In particular, this is the same measure $\tilde{\mathbb{P}}$ defined in \eqref{eq:tilde_tilt}, for which we have shown that 
\[
    \mathrm{Var}_{\tilde{\mathbb{P}}}\big(\log|1+X(p)p^{-z}|\big) = p^{-2z}+O\big(p^{-4z}\big)
\]
(cf.\ Equation \eqref{eq:lem7.1_tiltedvariance}). Combining this with the fact that 
\begin{align}
    \frac{1}{(1+E_p(S_1,S_2))^2}\prod_{j=1,2}\Big(1+\frac{16p^{-2z}}{(1-p^{-2z})^2}S_j(1-S_j)\Big)=1+O\big(p^{-2z}\big),
\end{align}
 uniformly over $p>\exp(e^{k_0})$ as a consequence of \eqref{eq:error_pestimate},
yields the estimate in \eqref{eq:nu_varianceestimate}.

 \bigskip

To prove the first claim, Equation \eqref{CondExpectSingleTwist}, one proceeds similarly as above, showing first that
    \begin{equation}
        \mathbb{E} \left[ \left| \frac{F_{k_0,k^*}(\sigma+i h)}{F_{k_0,k^*}(\sigma)} \right|^2\right] \asymp 1,
    \end{equation}
and then that
\begin{equation}
        \eta\Big(\mathcal{E}_{\leq}(k^*,V)\Big)\gg \frac{e^{-V}}{\sqrt{V}},\quad \text{ for }\quad  \frac{\mathrm{d}\eta}{\mathrm{d}\mathbb{P}}=\frac{| {F_{k_0,k^*}(\sigma+i h)}{F_{k_0,k^*}(\sigma)}^{-1} |^2}{\mathbb{E} [ | {F_{k_0,k^*}(\sigma+i h)}{F_{k_0,k^*}(\sigma)}^{-1} |^2]}.
\end{equation}
We only sketch the proofs of these claims, as they are straightforward modifications of previous arguments in this section. The first proceeds by using the identity $|1\pm p^{-z_j}|^2=(1\pm p^{-z})^2\mp 4p^{-z}\sin^2(|h_j|(\log p)/2)$ and Lemma \ref{lem:log_PNT}, as in the proof of \eqref{eq:CondExpectDoubleTwist}. The second follows from the proof of Proposition \ref{ProbLargeMax} upon substituting $\mathbb{P}$ by $\eta$ and making the necessary changes, noting that a simpler version of the argument used to prove \eqref{eq:nu_varianceestimate} shows that for $\mathrm{d}\tilde{\eta}/\mathrm{d}\eta=|F_{k_0,k^*}(\sigma)|^\lambda/\mathbb{E}_\eta[|F_{k_0,k^*}(\sigma)|^\lambda]$,
\begin{align}
        &\mathrm{Var}_{\tilde{\eta}}\big(\log|1+X(p)p^{-z}|\big)=\mathrm{Var}_{\tilde{\mathbb{P}}}\big(\log|1+X(p)p^{-z}|\big)\Big(1+O\big(p^{-2z}S\big)\Big),\\
        &\mathbb{E}_{\tilde{\eta}}\big[\log|1+X(p)p^{-z}|\big]=\mathbb{E}_{\tilde{\mathbb{P}}}\big[\log|1+X(p)p^{-z}|\big]+O\big(p^{-2z}S\big)
\end{align}
uniformly over $p>\exp(e^{k_0})$, for $S=\sin^2(|h|(\log p)/2)$.
\qed

\subsection{Proof of Proposition \ref{ModifiedBarrierLB}}
 To begin with, note that 
        \begin{align*}
        &\mathbb{E} \left[ \mathbf{1}(\mathcal{L}_{k^*}(h)) \mathbf{1}(\mathcal{E}_{>}(k^*)) |F_{k^*,\tilde t}(\sigma+ih)|^2 \right] \\&=\mathbb{E} \left[ \mathbf{1}(\mathcal{L}_{k^*}(h)) \mathbf{1}(\mathcal{E}_{>}(k^*)) |F_{k^*,t_\sigma}(\sigma+ih)|^2 \right]\mathbb{E}\big[|F_{t_\sigma,\tilde{t}}(\sigma+ih)|^2\big]\\
        &\asymp \mathbb{E} \left[ \mathbf{1}(\mathcal{L}_{k^*}(h)) \mathbf{1}(\mathcal{E}_{>}(k^*)) |F_{k^*,t_\sigma}(\sigma+ih)|^2 \right]
    \end{align*}
    by independence and Lemma \ref{lem:MomentsEulerProducts}, and it therefore suffices to prove that
    \begin{equation}
        \mathbb{E} \left[ \mathbf{1}(\mathcal{L}_{k^*}(h)) \mathbf{1}(\mathcal{E}_{>}(k^*)) |F_{k^*,t_\sigma}(\sigma+ih)|^2 \right] \gg \frac{e^{t_\sigma-k^*}}{\sqrt{t_\sigma-k^*}}.
    \end{equation}

    We will show that, roughly, every occurrence of $k^*$ in the left-hand side can be replaced by $k^*+j_0$, for $j_0$ a large absolute constant. This shift will ensure that the error terms arising in a Gaussian approximation step later in the proof are sufficiently small.
    
    To that end, let $j_0 = [C_1/4]$ be sufficiently large and let
    \[ \mathcal{L}_{k^*,k^*+j_0}(h) := \{ | \log |F_{k^*,k}(\sigma+ih)|| < C_1 \; \forall \; k^* < k \le k^* + j_0 \}, \]
    for $C_1$ as in the proposition's statement (if $t_\sigma - k^*\leq [C_1/4]$, take $j_0 = t_\sigma - k^*$). We also let
    \begin{equation*}\begin{split}
        \tilde{\mathcal{L}}_{k^*+j_0}(h) :=& \{ -C_1 (k- k^* -j_0) <  \log |F_{k^*+j_0,k}(\sigma+ih)| \\ &\le k-k^* - 2 \log(k-k^*) ; \forall k^*+j_0 < k \le t_\sigma \}. 
    \end{split}
    \end{equation*} 
    Then one can straightforwardly verify that
    \[ \mathcal{L}_{k^*}(h) \supseteq \mathcal{L}_{k^*,k^*+j_0}(h) \cap \tilde{\mathcal{L}}_{k^*+j_0}(h). \]
    This containment is immediate for $k \le k^*+j_0$, and otherwise follows from the fact that $\log|F_{k^*,k}(\sigma+ih)| = \log|F_{k^*,k^*+j_0}(\sigma+ih)| + \log|F_{k^*+j_0,k}(\sigma+ih)|$. Note in particular that under $\mathcal{L}_{k^*,k^*+j_0}(h)$, we have
    \begin{equation}\label{eq:bound_k^*_j_0_1}
        |F_{k^*,k^*+j_0}(\sigma+ih)|^2 \gg 1. 
    \end{equation} 
    Next, let
    \[ \tilde{\mathcal{E}}_>(k^*+j_0) := \{ \log | F_{k^*+j_0,k}(\sigma)| <0 \; \forall \; k^*+j_0 < k \le t_\sigma \} \]
    and
    \[ \mathcal{M}(k^*,k^*+j_0) := \big\{ \log |F_{k-1,k}(\sigma)| \le \mathbb{E}\big[\log |F_{k-1,k}(\sigma)|\big] \; \forall k^* < k \le k^*+j_0 \big\}, \]
    noting that by definition of $j_0$ and the negativity of $\mathbb{E}[\log|F_{k-1,k}(\sigma)|]$ for each $k$ (cf.\ Equation \eqref{eq:lem6.1:mertens}), we have
    \[\mathcal{E}_>(k^*) \supseteq \mathcal{M}(k^*,k^*+j_0) \cap \tilde{\mathcal{E}}_>(k^*+j_0). \]
    Using the independence of the $X(p)$ on different intervals of primes, and the estimate in \eqref{eq:bound_k^*_j_0_1}, we can then write
    \begin{equation}\begin{split}\label{jun24.2}
        &\mathbb{E} \left[ \mathbf{1}(\mathcal{L}_{k^*}(h)) \mathbf{1}(\mathcal{E}_{>}(k^*)) |F_{k^*,t_\sigma}(\sigma+ih)|^2 \right] \\
        \ge \, &\mathbb{E} \left[ \mathbf{1}(\mathcal{L}_{k^*,k^*+j_0}(h)) \mathbf{1}(\mathcal{M}(k^*,k^*+j_0)) |F_{k^*,k^*+j_0}(\sigma+ih)|^2 \right] \times \\
        &\hspace{100pt}
        \mathbb{E} \left[ \mathbf{1}(\tilde{\mathcal{L}}_{k^*+j_0}(h)) \mathbf{1}(\tilde{\mathcal{E}}_{>}(k^*+j_0)) |F_{k^*+j_0,t_\sigma}(\sigma+ih)|^2 \right] \\
        \gg \, &\mathbb{P} \left( \mathcal{L}_{k^*,k^*+j_0}(h) \cap \mathcal{M}(k^*,k^*+j_0) \right) 
        \mathbb{E} \left[ \mathbf{1}(\tilde{\mathcal{L}}_{k^*+j_0}(h)) \mathbf{1}( \tilde{\mathcal{E}}_{>}(k^*+j_0)) |F_{k^*+j_0,t_\sigma}(\sigma+ih)|^2 \right],
    \end{split}\end{equation}
    and must therefore show that
    \begin{align}\label{jun24.1}
        \mathbb{P} \left( \mathcal{L}_{k^*,k^*+j_0}(h) \cap \mathcal{M}(k^*,k^*+j_0) \right) \gg 1.
    \end{align}

   A union bound on the complement of the event in the left-hand side gives
    \[ \mathbb{P} \left( \mathcal{L}_{k^*,k^*+j_0}(h) \cap \mathcal{M}(k^*,k^*+j_0) \right) \ge 1 - \mathbb{P} \left( \mathcal{L}_{k^*,k^*+j_0}(h)^c \right) - \mathbb{P} \left( \mathcal{M}(k^*,k^*+j_0)^c \right). \]
    Noting that
    \[ \mathbb{P} \left( \mathcal{M}(k^*,k^*+j_0)^c \right) = 1 - \mathbb{P} \left( \mathcal{M}(k^*,k^*+j_0) \right) \le 1 - 2^{-j_0} \]
    by the independence and symmetry in law of the increments $\log|F_{k-1,k}(\sigma)|$, it thus suffices to show that
    \begin{equation}\label{jun26.1}
        \mathbb{P} \left( \mathcal{L}_{k^*,k^*+j_0}(h)^c \right) \le \frac{1}{2^{j_0+1}}
    \end{equation}
    in order to establish \eqref{jun24.1}.
    We do so by applying another union bound, followed by Markov's inequality, to deduce that 
    \begin{align*}
        \mathbb{P} \left( \mathcal{L}_{k^*,k^*+j_0}(h)^c \right) 
        &\le \sum_{k^* < k \le k^*+j_0} \mathbb{P}\big( | \log |F_{k^*,k}(\sigma+ih)| | > C_1 \big) \\
        &\le \sum_{k^* < k \le k^*+j_0} \frac{\mathbb{E}[|F_{k^*,k}(\sigma+ih)|] + \mathbb{E}[|F_{k^*,k}(\sigma+ih)|^{-1}]}{e^{C_1}}.
    \end{align*}
    Since Lemma \ref{lem:MomentsEulerProducts} implies that
    \begin{align*}
        \mathbb{E}[|F_{k^*,k}(\sigma+ih)|^{\pm 1}] &\le 2 \exp \left( \sum_{k^* < \log \log p \le k} \frac{1/4 + (1/4 \mp 1/2) \cos(2h \log p)}{p^{1+2 \sigma}} \right) \\
        &\le 2 \exp \left( \sum_{k^* < \log \log p \le k^*+j_0} \frac{1}{p} \right) \le 2 e^{2 j_0}
    \end{align*}
    for $k^* > k_0$ and $j_0$ both sufficiently large, it follows that
    \[ \mathbb{P} \left( \mathcal{L}_{k^*,k^*+j_0}(h)^c \right) \le \frac{4 j_0 e^{2j_0}}{e^{C_1}} \]
    which is $\leq 2^{-(j_0+1)}$ using the definition of $j_0$.

\bigskip 

    We are now left with the task of showing that
    \begin{equation}\label{eq:long_lower_ballot}
        \mathbb{E} \left[ \mathbf{1}(\tilde{\mathcal{L}}_{k^*+j_0}(h)) \mathbf{1}( \tilde{\mathcal{E}}_{>}(k^*+j_0)) |F_{k^*+j_0,t_\sigma}(\sigma+ih)|^2 \right] \gg \frac{e^{t_\sigma-k^*}}{\sqrt{t_\sigma- k^*}}.
    \end{equation}
    To do so, we introduce the tilted measure
    $$\frac{\mathrm{d}{\mathbb{Q}}}{\mathrm{d}\mathbb{P}} = \frac{|F_{t_\sigma}(\sigma+ih)|^2}{\mathbb{E}[|F_{t_\sigma}(\sigma+ih)|^2]},$$
    and rewrite the expectation in \eqref{eq:long_lower_ballot} as
    \[
        \mathbb{E}\big[|F_{k^*+j_0,t_\sigma}(\sigma+ih)|^2\big]\mathbb{Q}\Big(\tilde{\mathcal{L}}_{k^*+j_0}(h)\cap \tilde{\mathcal{E}}_{>}(k^*+j_0)\Big).
    \]
    The second moment of $|F_{k^*+j_0,t_\sigma}(\sigma+ih)|$ is $e^{t_\sigma-k^*-j_0} \asymp e^{t_\sigma-k^*}$ by Lemma \ref{lem:MomentsEulerProducts}. As for the remaining probability, we can bound it below by 
    \begin{equation}\label{eq:tilted_(2,0)_prob}
    \mathbb{Q}\Big(\tilde{\mathcal{L}}_{k^*+j_0}(h)\cap \tilde{\mathcal{E}}_{>}(k^*+j_0) \cap \big\{-C_1(k-k^*-j_0)\leq \log|F_{k^*+j_0,k}(\sigma)|\,\forall \,k^*+j_0<k\leq t_\sigma\big\}\Big)
    \end{equation}
    setting the stage for an application of Proposition \ref{GaussianComparisonGlobal}. Indeed, note that the probability of the event in \eqref{eq:tilted_(2,0)_prob} under $\mathbb{Q}$ is the same as under the measure $\tilde{\mathbb{P}}_{(h,0)}^{(2,0)}$ in this proposition, and that a tedious but nevertheless direct calculation shows that its conditions are satisfied for
    \begin{align*}
                &n_1=k^*+j_0, \quad n_2=t_\sigma,\\
                &L_k^{(1)}=L_k^{(2)}=-C_1(k-k^*-j_0)+ 1=-C_1(k-n_1) + 1,\\
                &U_k^{(1)}=(k-k^*)-2\log(k-k^*),\\
                &U_k^{(2)}=-1.
    \end{align*}
    We therefore conclude by Proposition \ref{GaussianComparisonGlobal} (more specifically, its corollary in \eqref{eq:remark_(2,0)}) that \eqref{eq:tilted_(2,0)_prob} is
    \begin{equation}
    \begin{split}
                 &\asymp \mathbb{P} \Big(L_k^{(1)} - (k-n_1)+O(1) \le \sum_{m = n_1+1}^k \overline{N}_m^{(1)} \le U_k^{(1)} - (k-n_1)+O(1) \; \forall \; n_1 < k \le n_2 \Big) \\&\times  \mathbb{P} \Big(L_k^{(2)}+\frac{(k-n_1)}{2}+O(1) \le \sum_{m = n_1+1}^k \overline{N}_m^{(2)} \le U_k^{(2)} +\frac{(k-n_1)}{2}+O(1) \; \forall \; n_1 < k \le n_2 \Big).
                 \end{split}
    \end{equation}
    Using standard prime sum estimates, one can show that the variances of $(\overline{N}_m^{(1)})_m$ and $(\overline{N}_m^{(2)})_m$ are in $[1/20,20]$ for each $m$, by taking $C_1$ (and thus $j_0$) sufficiently large. The ballot theorem in \cite[Probability Result 2]{Harper1} then implies that the first probability is
    \[
        \gg \frac{1}{\sqrt{t_\sigma-k^*}},
        \]
    whereas the second probability is easily shown to be $\gg 1$ by a union bound on its complement and standard Gaussian estimates, after conditioning on $\overline{N}_{n_1+1}^{(2)}\in [-C-1,-C]$ for a large constant $C>0$.
    \qed
\appendix

\section{Proof of Proposition \ref{ShortIntSecMomBarr}}\label{appendix}

We may and will assume throughout this section that $t_\sigma-k^*$ is sufficiently large. Otherwise, the claim follows by using Fubini's theorem to write
\begin{align*}
        &\mathbb{E}\left[ \left( e^{2 k^*-t_\sigma} \int_{L_{k^*}} |F_{k^*,\tilde{t}}(\sigma+ih)|^2 \mathrm{d}h \right)^{2} \right]\\ &\leq e^{4 k^*-2t_\sigma} \iint_{[e^{-(k^*+1)}, e^{-k^*}]^2}\mathbb{E}\big[|F_{k^*,\tilde{t}}(\sigma+ih_1)|^2|F_{k^*,\tilde{t}}(\sigma+ih_2)|^2\big] \mathrm{d}h_1\mathrm{d}h_2,
\end{align*}
and bounding the integrand using Lemma \ref{lem:MomentsEulerProducts}. We begin with some preliminary estimates.
\begin{lem}\label{EulerProductLemma2}
    Let $h_1, h_2, u_1, u_2, \lambda_1, \lambda_2$ be bounded real numbers and suppose that $y, z, \sigma$ are real numbers satisfying $100 (1+\lambda_1^2+\lambda_2^2+u_1^2+u_2^2) \le y \le z$ and $\sigma \ge -10/ \log z$. Then we have
    \begin{equation}\begin{split}
        \mathbb{E} & \bigg[\prod_{y < p \leq z}\Big|1+\frac{X(p)}{p^{1/2+\sigma+ih_1}}\Big|^{\lambda_1+iu_1}\Big|1+\frac{X(p)}{p^{1/2+\sigma+ih_2}}\Big|^{\lambda_2+iu_2}\bigg] \\
        = \exp\Bigg(&\sum_{y < p \leq z}\frac{(\lambda_1+iu_1)^2 + (\lambda_2+iu_2)^2}{4 p^{1+2\sigma}}+ \frac{(\lambda_1+iu_1)(\lambda_2+iu_2) \cos((h_1-h_2)\log p)}{2 p^{1+2\sigma}} \\
        &+ \sum_{y < p \leq z} \frac{(\lambda_1+iu_1)(\lambda_1 - 2 +i u_1) \cos(2 h_1 \log p)+(\lambda_2+iu_2)(\lambda_2 - 2 + i u_2) \cos(2 h_2 \log p)}{4 p^{1+2\sigma}} \\
        &+ \sum_{y < p \leq z} \frac{(\lambda_1+iu_1)(\lambda_2+iu_2) \cos((h_1+h_2) \log p)}{2 p^{1+2\sigma}}
        + T(u_1,u_2) \Bigg),
    \end{split}\end{equation}
    where $T(u,v)$ satisfies
    \begin{align*} |T(u,v)|&\ll \frac{1 + |u|^3+|v|^3}{\sqrt{y} \log y}, \qquad\frac{\partial T(u,v)}{\partial u}, \frac{\partial T(u,v)}{\partial v} \ll  \frac{1 + u^2+v^2}{\sqrt{y} \log y} \\
    \frac{\partial^2 T(u,v)}{\partial u \partial v} &\ll \frac{1 + |u|+|v|}{\sqrt{y} \log y}.
    \end{align*}
\end{lem}

\begin{proof}
    This is \cite[Lemma 6]{Harper1} with the necessary modifications for the Rademacher case as outlined in Section 5.3 there.
\end{proof}

In what follows, for any $h_1,h_2, \lambda_1,\lambda_2 \in \mathbb{R}$, let $\tilde{\mathbb{P}}_{h_1,h_2}^{(\lambda_1,\lambda_2)}$ denote the measure defined by
\[
    \frac{\mathrm{d}\tilde{\mathbb{P}}_{h_1,h_2}^{(\lambda_1,\lambda_2)}}{\mathrm{d}\mathbb{P}}=\frac{|F_{t_\sigma}(\sigma+ih_1)|^{\lambda_1
    }|F_{t_\sigma}(\sigma+ih_2)|^{\lambda_2
    }}{\mathbb{E}[|F_{t_\sigma}(\sigma+ih_1)|^{\lambda_1
    }|F_{t_\sigma}(\sigma+ih_2)|^{\lambda_2
    }]}.
\]
\begin{lem}\label{jun18.1}
    Let $n_1 < n_2 \in \mathbb{N}$ be sufficiently large and $B > 0$ a sufficiently large constant. Suppose that $0 \le h_1, h_2 \le 1$ are real numbers satisfying
    \[ |h_1 - h_2| \ge {Be^{-n_1}}. \]
    Suppose that $|\sigma| \le e^{-n_2}$ is a real number and let $(u_k)_{k=n_1+1}^{n_2}, (v_k)_{k={n_1+1}}^{n_2}$ be sequences of real numbers satisfying (say)
    \[ |u_k|, |v_k| \le \frac{e^{k/2}}{40} + 2. \]
    Lastly, let $Y_k(h_j)=\log|F_{k-1,k}(\sigma+ih_j)|$ for $j\in \{1,2\}$. Then we have
    \begin{equation}\begin{split}\label{SlicedGaussianComparison}
        &\tilde{\mathbb{P}}_{h_1,h_2} \Big( \substack{u_k \le Y_k(h_1) \le u_k + (k-n_1)^{-2} \\v_k \le Y_k(h_2) \le v_k + (k-n_1)^{-2} }\;, \forall \; n_1 < k \le n_2 \Big) \\ 
        &= \left(1 + O\left( e^{-e^{n_1}/100} \right) \right) \mathbb{P} \Big( \substack{u_k \le N_k^{(1)} \le u_k + (k-n_1)^{-2}\\ v_k \le N_k^{(2)} \le v_k + (k-n_1)^{-2}} \;, \forall \; n_1 < k \le n_2 \Big),
    \end{split}\end{equation}
    where $(N_k^{(1)},N_k^{(2)})_{k=n_1+1}^{n_2}$ is a sequence of independent bivariate Gaussians with mean vector
    \[ \boldsymbol{\mu}_k = \sum_{k-1 < \log \log p \le k} \frac{1}{2 p^{1+2 \sigma}} 
    \begin{pmatrix}
        \lambda_1 + \lambda_2 C_p(h_1-h_2) + (\lambda_1-1) C_p(2 h_1 ) + \lambda_2 C_p(h_1+h_2) \\
        \lambda_2 + \lambda_1 C_p(h_1-h_2)+ (\lambda_2-1) C_p(2 h_2) + \lambda_1 C_p(h_1+h_2)
    \end{pmatrix}
    \]
    and covariance matrix
    \[  
    \boldsymbol{\Sigma}_k = \sum_{k-1 < \log \log p \le k} \frac{1}{2p^{1+2\sigma}}
    \begin{pmatrix}
        1 + C_p(2 h_1) & C_p(h_1-h_2) + C_p(h_1+h_2) \\
        C_p(h_1-h_2)+ C_p(h_1+h_2) & 1 + C_p(2 h_2)
    \end{pmatrix},
    \]
    for $C_p(\cdot):=\cos(\cdot \log p)$.
    If we additionally have
    \[ |u_k|, |v_k| \le \frac{e^{k/2} \sqrt{|h_1-h_2|}}{(k-n_1)^2}-1, \]
    then
    \begin{equation}\begin{split}
        &\tilde{\mathbb{P}}_{h_1,h_2} \Big( \substack{u_k \le Y_k(h_1) \le u_k + (k-n_1)^{-2} \\v_k \le Y_k(h_2) \le v_k + (k-n_1)^{-2} }\;, \forall \; n_1 < k \le n_2 \Big) \\& \asymp \mathbb{P} \Big( u_k \le N_k^{(1)} \le u_k + (k-n_1)^{-2} \; \forall \; n_1 < k \le n_2 \Big)\mathbb{P} \Big( v_k \le N_k^{(2)} \le v_k + (k-n_1)^{-2} \; \forall \; n_1 < k \le n_2 \Big).
    \end{split}\end{equation}
\end{lem}

\begin{proof}
    One verifies by a straightforward calculation using Lemma \ref{EulerProductLemma2} that we picked $\boldsymbol{\mu}_k$ and $\boldsymbol{\Sigma}_k$ so that the characteristic function $e^{i u_1 N_k^{(1)} + i u_2 N_k^{(2)}}$ matches precisely the main term in the characteristic function 
    \[ \tilde{\mathbb{E}}_{h_1,h_2} \left[ e^{i u_1 Y_k(h_1) + i u_2Y_k(h_2)|} \right] = \frac{\mathbb{E}\left[|F_{k-1,k}(\sigma+ih_1)|^{\lambda_1+iu_1} |F_{k-1,k}(\sigma+ih_2)|^{\lambda_2+iu_2} \right]}{\mathbb{E} \left[ |F_{k-1,k}(\sigma+ih_1)|^{\lambda_1} |F_{k-1,k}(\sigma+ih_2)|^{\lambda_2} \right]}. \]
    By following the exact argument of the proof of \cite[Lemma 7]{Harper1}, one deduces that for a given $k$ with $n_1 < k \le n_2$,
    \begin{equation}\begin{split}
        &\tilde{\mathbb{P}}_{h_1,h_2} \Big( u_k \le Y_k(h_1) \le u_k + (k-n_1)^{-2}, v_k \le Y_k(h_2) \le v_k + (k-n_1)^{-2} \Big) \\
        = \, &\mathbb{P} \Big( u_k \le N_k^{(1)} \le u_k + (k-n_1)^{-2}, v_k \le N_k^{(2)} \le v_k + (k-n_1)^{-2} \Big) + O\left(e^{-e^k/50}\right).
    \end{split}\end{equation}
    Similarly, the rest of the proof there goes through without change after noting that by our assumptions, we have
    \[ \sum_{k-1 < \log \log p \le k } \frac{\cos((h_1 \pm h_2)\log p)}{p^{1+2 \sigma}} \ll \frac{1}{|h_1 \pm h_2| e^k} \ll \frac{1}{B e^{k-n_1}}. \]
    In particular, these terms are small when $B$ is sufficiently large. In the same way, one shows that say
    \[ \mathbb{P} \Big( u_k \le N_k^{(1)} \le u_k + (k-n_1)^{-2}, v_k \le N_k^{(2)} \le v_k + (k-n_1)^{-2} \Big) \gg e^{-e^k/100} \]
    using the bounds on $|u_k|$ and $|v_k|$, which allows one to turn the absolute error term into a relative one and deduce the first claim by taking the product over $k$.

    For the second claim, let $\rho_k := \frac{\Sigma_{1,2}}{\sqrt{\Sigma_{1,1} \Sigma_{2,2}}}$ be the correlation coefficient of $(N_k^{(1)},N_k^{(2)})$ and note that since $\Sigma_{i,i} \asymp 1$, we have from Lemma \ref{lem:PNT}
    \[ \rho_k \ll \frac{1}{|h_1-h_2| e^k}. \]
    Moreover, we can ensure that $|\rho_k| \le \tfrac{1}{2}$ for $n_1 < k \le n_2$ by making $B$ sufficiently large. It follows that
    \begin{align*}
        &\mathbb{P} \Big( u_k \le N_k^{(1)} \le u_k + (k-n_1)^{-2}, v_k \le N_k^{(2)} \le v_k + (k-n_1)^{-2} \Big) \\
        &= \int_{\frac{u_k-\mu_k^{(1)}}{\sigma_k^{(1)}}}^{\frac{u_k+(k-n_1)^{-2}-\mu_k^{(1)}}{\sigma_k^{(1)}}} \int_{\frac{v_k-\mu_k^{(2)}}{\sigma_k^{(2)}}}^{\frac{v_k+(k-n_1)^{-2}-\mu_k^{(2)}}{\sigma_k^{(2)}}} e^{-(x^2-2 \rho_k x y+y^2)/2(1-\rho_k^2)} \mathrm{d}x \mathrm{d}y \\
        &= \exp(O((k-n_1)^{-4})) \int_{\frac{u_k-\mu_k^{(1)}}{\sigma_k^{(1)}}}^{\frac{u_k+(k-n_1)^{-2}-\mu_k^{(1)}}{\sigma_k^{(1)}}} \int_{\frac{v_k-\mu_k^{(2)}}{\sigma_k^{(2)}}}^{\frac{v_k+(k-n_1)^{-2}-\mu_k^{(2)}}{\sigma_k^{(2)}}} e^{-(x^2+y^2)/2} \mathrm{d}x \mathrm{d}y.
    \end{align*}
    Here, we used in the last step that on the range of integration, we have
    \[ \rho_k x y \ll \frac{(|u_k|+1)( |v_k| + 1)}{|h_1 - h_2| e^k} \ll \frac{1}{(k-n_1)^4} \]
    and similar (but easier) estimates for $\rho_k^2 x^2$ and $\rho_k^2 y^2$ as well as for $(1-\rho_k^2)^{-1/2}$.
    But this is in turn
    \begin{align*} = &e^{O((k-n_1)^{-4})} \mathbb{P} \Big( u_k \le N_k^{(1)} \le u_k + (k-n_1)^{-2}  \Big) \mathbb{P} \Big( v_k \le N_k^{(2)} \le v_k + (k-n_1)^{-2} \Big). 
    \end{align*}
    Taking the product over all $n_1 < k \le n_2$ gives the second claim.
\end{proof}
\begin{prop}\label{GaussianComparisonGlobal}
    Suppose that $0\leq h_1, h_2\leq 1$.
    Let $C$ be a sufficiently large constant, and $k^* \le n_1 < k \le n_2$ be positive integers. Suppose further that
    \[ n_1 \ge 2 \log |h_1 - h_2|^{-1} - k^* + C \]
    and that $|\sigma| \le e^{-n_2}$ is a real number. Lastly, let  $L_k^{(1)} \le  U_k^{(1)}$ and $L_k^{(2)} \le U_k^{(2)}$ be real numbers satisfying
    \[ |L_k^{(1)}|, |L_k^{(2)}|, |U_k^{(1)}|, |U_k^{(2)}| \le \frac{e^{(k-k^*)/4}}{250(k-k^*)^2}, \]
    and $({N}_k^{(1)}, N_k^{(2)})_k$ be a sequence of independent bivariate Gaussians with means $\boldsymbol{\mu}_k$ (as in Lemma \ref{jun18.1}), and covariance matrices
    \[  
    \boldsymbol{\Sigma}_k = \sum_{k-1 < \log \log p \le k}\frac{1}{2p^{1+2\sigma}}
    \begin{pmatrix}
        1 + C_p(2 h_1) & 0\\
        0 & 1 + C_p(2 h_2)
    \end{pmatrix},
    \]
    for $C_p(\cdot):=\cos(\cdot \log p)$. Then we have 
    \begin{align*}\label{splitA3conclusion}
\tilde{\mathbb{P}}_{h_1,h_2}^{(\lambda_1,\lambda_2)} &\Big( L_k^{(1)} \le \log |F_{n_1,k}(\sigma+i h_1)| \le U_k^{(1)}, L_k^{(2)} \le \log |F_{n_1,k}(\sigma+i h_2)| \le U_k^{(2)} \; \forall \; n_1 < k \le n_2 \Big) \\
        \hspace{30pt} &\ll\mathbb{P} \Big(L_k^{(1)} -2 \le \sum_{m = n_1+1}^k {N}_m^{(1)} \le U_k^{(1)}+2  \; \forall \; n_1 < k \le n_2 \Big) \\&\quad\times\mathbb{P} \Big(L_k^{(2)}-2  \le \sum_{m = n_1+1}^k {N}_m^{(2)} \le U_k^{(2)}+2  \; \forall \; n_1 < k \le n_2 \Big),
    \end{align*}
    as well as the lower bound 
\begin{equation}\label{splitA3conclusion_lower}\begin{split}
        &\gg \mathbb{P} \Big(L_k^{(1)} +2 \le \sum_{m = n_1+1}^k {N}_m^{(1)} \le U_k^{(1)}-2  \; \forall \; n_1 < k \le n_2 \Big) \\&\quad\times\mathbb{P} \Big(L_k^{(2)}+2  \le \sum_{m = n_1+1}^k {N}_m^{(2)} \le U_k^{(2)}-2  \; \forall \; n_1 < k \le n_2 \Big).
    \end{split}\end{equation}
\end{prop}

\begin{proof}
    This follows from the ``slicing'' argument used to prove \cite[Lemma 4]{Harper1}. We only give the details for the upper bound, the lower bound following by a similar argument.
    
    Note first that if $L_k^{(1)} \le \log |F_{n_1,k}(\sigma+i h_1)| \le U_k^{(1)}$ then in particular we have
    \[ \left|\log |F_{k-1,k}(\sigma+ih_1)| \right| \le |U_{k}^{(1)} - L_{k-1}^{(1)}| \le \frac{e^{(k+1-k^*)/4}}{125(k-k^*)^2} \le \frac{e^{(k-k^*)/4}}{40(k-k^*)^2} \le \frac{e^{k/2}}{40}, \]
    and similarly for $\log |F_{k-1,k}(\sigma+ih_2)|$. For $n_1 < k \le n_2$, let 
    \[ \mathcal{R}_k := \left\{ u \in (k-n_1)^{-2} \mathbb{Z} : |u| \le \frac{e^{(k-k^*)/4}}{40 (k-k^*)^2} + 2 \right\}. \]
    If $L_k^{(1)} \le \log |F_{n_1,k}(\sigma+i h_1)| \le U_k^{(1)}$, then there must exist a tuple $(u_k^{(1)})_{k=n_1+1}^{n_2} \in \prod_{k=n_1+1}^{n_2} \mathcal{R}_k$ such that
    \[ u_k^{(1)} \le \log |F_{k-1,k}(\sigma+i h_1)| \le u_k^{(1)}+(k-n_1)^{-2} \]
    and
    \[ L_k^{(1)} - \sum_{m=n_1+1}^{n_2} \frac{1}{(m-n_1)^2} \le \sum_{m=n_1+1}^k u_m^{(1)} \le U_k^{(1)}, \]
    and the analogous estimate at $h_2$. Next, note that by our assumption on $n_1$, we have
    \[ |u_k^{(1)}|, |u_k^{(2)}| \le \frac{e^{(k-k^*)/4}}{(k-k^*)^2}+2 = e^{k/2} \frac{e^{-(k+k^*)/4}}{(k-k^*)^2} +2 \le \frac{e^{k/2}}{(k-k^*)^2} e^{-(n_1+k^*)/4} + 2 \le \frac{e^{k/2} \sqrt{|h_1-h_2|}}{(k-k^*)^2} - 1.  \]
    We are thus in a position to apply (the second part of) Lemma \ref{jun18.1} to deduce that
    \begin{align*}
        &\tilde{\mathbb{P}}_{h_1,h_2} \Big( L_k^{(1)} \le \log |F_{n_1,k}(\sigma+i h_1)| \le U_k^{(1)}, L_k^{(2)} \le \log |F_{n_1,k}(\sigma+i h_2)| \le U_k^{(2)} \; \forall \; n_1 < k \le n_2 \Big) \\
        &\ll \sum_{\substack{ u_k^{(1)}, u_k^{(2)} \in \mathcal{R}_k \forall n_1 < k \le n_2 \\  L_k^{(i)} - 2 \le \sum_{m=n_1+1}^k u_m^{(i)} \le U_k^{(i)}, i=1,2 }} \mathbb{P} \Big( u_k^{(1)} \le N_k^{(1)} \le u_k^{(1)} + (k-n_1)^{-2}\; \forall \; n_1 < k \le n_2 \Big) \times \\[-20pt] &\hspace{200pt}\mathbb{P} \Big( u_k^{(2)} \le N_k^{(2)} \le u_k^{(2)} + (k-n_1)^{-2} \; \forall \; n_1 < k \le n_2 \Big) \\
        &\ll \mathbb{P} \Big(L_k^{(1)} - 2 \le \sum_{m = n_1+1}^k N_m^{(1)} \le U_k^{(1)} + 2 \; \forall \; n_1 < k \le n_2 \Big) \times \\
        &\hspace{150pt} \mathbb{P} \Big(L_k^{(2)} - 2 \le \sum_{m = n_1+1}^k N_m^{(2)} \le U_k^{(2)} + 2 \; \forall \; n_1 < k \le n_2 \Big).
    \end{align*}
\end{proof}
 
\begin{rem}
    We record the following refinements of the conclusion of the previous proposition in two special cases, which arise in the proofs of Proposition \ref{ModifiedBarrierLB} and of Proposition \ref{ShortIntSecMomBarr} below.

    The first is when $\lambda_1=\lambda_2=2$ and $h_1,h_2\geq e^{-(k^*+1)}$, in which case 
 \[ \sum_{m = n_1+1}^k \boldsymbol{\mu}_m = 
    \begin{pmatrix}
        k - n_1 + O(1) \\
        k - n_1 + O(1)
    \end{pmatrix} \]
    and the bounds in the previous proposition can be written as 
    \begin{equation}\begin{split}\label{eq:centered_final}
                 \asymp\,\, &\mathbb{P} \Big(L_k^{(1)} - (k-n_1)+O(1) \le \sum_{m = n_1+1}^k \overline{N}_m^{(1)} \le U_k^{(1)} - (k-n_1)+O(1) \; \forall \; n_1 < k \le n_2 \Big) \\&\times  \mathbb{P} \Big(L_k^{(2)} - (k-n_1)+O(1) \le \sum_{m = n_1+1}^k \overline{N}_m^{(2)} \le U_k^{(2)} -(k-n_1)+O(1) \; \forall \; n_1 < k \le n_2 \Big),
                \end{split}
    \end{equation}
    where $(\overline{N}_{k}^{(1)},\overline{N}_{k}^{(2)})$ now has mean zero, and the same covariances as $({N}_{k}^{(1)},{N}_{k}^{(2)})$.
    
    The second is when $(\lambda_1,\lambda_2)=(2,0)$, $h_2=0$ and $h_1\geq Be^{-n_1}$, in which case one obtains
    \begin{equation}\label{eq:remark_(2,0)}
    \begin{split}
                 \asymp &\,\,\mathbb{P} \Big(L_k^{(1)} - (k-n_1)+O(1) \le \sum_{m = n_1+1}^k \overline{N}_m^{(1)} \le U_k^{(1)} - (k-n_1)+O(1) \; \forall \; n_1 < k \le n_2 \Big) \\&\times  \mathbb{P} \Big(L_k^{(2)}+\frac{(k-n_1)}{2}+O(1) \le \sum_{m = n_1+1}^k \overline{N}_m^{(2)} \le U_k^{(2)} +\frac{(k-n_1)}{2}+O(1) \; \forall \; n_1 < k \le n_2 \Big)
                 \end{split}
    \end{equation}
    by the same argument.
\end{rem}

\begin{proof}[Proof of Proposition \ref{ShortIntSecMomBarr}]
    We begin by expanding the square to see that
    \begin{equation}\begin{split}\label{jun10.1}
        &\mathbb{E}\left[ \left( e^{2k^*-t_\sigma} \int_{{L}_{k^*}} |F_{k^*,\tilde t}(\sigma+ih)|^2 \mathrm{d}h \right)^{2} \right] \\
        &= e^{4k^*-2t_\sigma} \iint \mathbb{E} \big[ \mathbf{1} (\mathcal{L}_{k^*}(h_1)) \mathbf{1} (\mathcal{L}_{k^*}(h_2)) |F_{k^*,\tilde t}(\sigma+ih_1)|^2 |F_{k^*,\tilde t}(\sigma+ih_2)|^2 \big] \mathrm{d}h_1 \mathrm{d}h_2,
    \end{split}\end{equation}
    where the range of integration is $[e^{-(k^*+1)},e^{-k^*}]^2$.
    Since $\mathcal{L}_{k^*}(h_1),\mathcal{L}_{k^*}(h_2)$ only depend on primes $\le e^{e^{t_\sigma}}$ and $t_\sigma < \tilde{t}$ for sufficiently large $U$, we can factor out the larger primes. By Lemma \ref{lem:MomentsEulerProducts}, we have
    \begin{align*}
        \mathbb{E}[|F_{t_\sigma,\tilde t}(\sigma+ih_1)|^2 |F_{t_\sigma,\tilde t}(\sigma+ih_2)|^2] \ll \exp \left( \sum_{t_\sigma < \log \log p \le \tilde t} \frac{6}{p^{1+2 \sigma}} \right) \ll 1
    \end{align*}
    uniformly in $h_1, h_2$. Thus the right-hand side of \eqref{jun10.1} is in turn
    \begin{align}\label{jun10.2}
        \ll e^{4k^*-2t_\sigma} \iint \mathbb{E} \big[ \mathbf{1} (\mathcal{L}_{k^*}(h_1)) \mathbf{1} (\mathcal{L}_{k^*}(h_2)) |F_{k^*,t_\sigma}(\sigma+ih_1)|^2 |F_{k^*,t_\sigma}(\sigma+ih_2)|^2 \big] \mathrm{d}h_1 \mathrm{d}h_2.
    \end{align}
    We will split the range of integration into $3$ subranges according to the size of $|h_1 - h_2|$ as follows.

    Suppose first that $|h_1 - h_2| \le e^{-t_\sigma} \asymp \frac{U}{\log x}$. Under the event $\mathcal{L}_{k^*}(h_2)$ applied at the final point $k = t_\sigma$, we have
    \[ |F_{k^*,t_\sigma}(\sigma+i h_2)|^2 \le \frac{e^{2(t_\sigma-k^*)}}{(t_\sigma-k^*)^4}. \]
    On this range, we can now get rid of the barrier events to deduce that the contribution from $|h_1 - h_2 | \le e^{-t_\sigma}$ to \eqref{jun10.2} is
    \begin{align*}
        \ll \frac{e^{2 k^* - t_\sigma}}{(t_\sigma-k^*)^4} \int_{e^{-(k^*+1)}}^{e^{-k^*}} \mathbb{E}[|F_{k^*,t_\sigma}(\sigma+ih_1)|^2] \mathrm{d}h_1 \ll \frac{1}{(t_\sigma-k^*)^4}
    \end{align*}
    using Lemma \ref{lem:MomentsEulerProducts}, which is more than acceptable.

    Next, suppose that $e^{-t_\sigma} < |h_1 - h_2| \le e^{-t_\sigma} (e^{t_\sigma-k^*})^{2/3}$. We follow the same strategy as in the previous range, this time applying the barrier event $\mathcal{L}_{k^*}(h_2)$ with the choice $k = M := [\log |h_1-h_2|^{-1}]$. Note that this choice satisfies $k^* < M \le t_\sigma$. We deduce that under $\mathcal{L}_{k^*}(h_2)$ we have
    \begin{align*}
        |F_{k^*,M}(\sigma+i h_2)|^2 \ll \frac{1}{|h_1-h_2|^2 e^{2 k^*} (\log |h_1-h_2|^{-1} - k^*)^4 } \ll \frac{1}{|h_1-h_2|^2 e^{2 k^*} (t_\sigma-k^*)^4}.
    \end{align*}
    Again discarding the barrier events after applying this bound and using Lemma \ref{lem:MomentsEulerProducts} (together with standard prime sum estimates), the contribution to \eqref{jun10.2} from $e^{-t_\sigma} < |h_1 - h_2| \le e^{-t_\sigma} (e^{t_\sigma-k^*})^{2/3}$ with $e^{-(k^*+1)} < h_1, h_2 \le e^{-k^*}$ is
    \begin{align*}
        &\ll \frac{e^{2k^*-2t_\sigma}}{(t_\sigma-k^*)^4} \iint \mathbb{E} \big[  |F_{M,t_\sigma}(\sigma+ih_1)|^2 |F_{M,t_\sigma}(\sigma+ih_2)|^2 \big] \mathbb{E} \big[  |F_{k^*,M}(\sigma+ih_1)|^2 \big] \frac{\mathrm{d}h_1 \mathrm{d}h_2}{|h_1 - h_2|^2} \\
        &\ll \frac{e^{2k^*-2t_\sigma}}{(t_\sigma-k^*)^4} \iint (e^{t_\sigma- M})^2 e^{M-k^*} \frac{\mathrm{d}h_1 \mathrm{d}h_2}{|h_1 - h_2|^2} \ll \frac{e^{k^*}}{(t_\sigma-k^*)^4} \iint \frac{\mathrm{d}h_1 \mathrm{d}h_2}{|h_1-h_2|} \ll \frac{1}{(t_\sigma-k^*)^3}.
    \end{align*}

    It thus remains to bound the contribution to \eqref{jun10.2} of $h_1, h_2$ satisfying 
    \begin{equation}\label{jun19.2}
    e^{-t_\sigma} (e^{t_\sigma-k^*})^{2/3} \le |h_1 - h_2| \le e^{-k^*}(1-e^{-1}). 
    \end{equation}
    Similar to before, we start by employing the barrier event $\mathcal{L}_{k^*}(h_2)$ with the choice $k = M := [\log |h_1-h_2|^{-1}]+1$, which satisfies $k^* < M \le t_\sigma$. This gives that
    \[ |F_{k^*,M}(\sigma+i h_2)|^2 \ll \frac{1}{|h_1-h_2|^2 e^{2 k^*} (\log |h_1-h_2|^{-1} - k^*)^4 }. \]
    It follows that
    \begin{equation}\begin{split}\label{jun19.1}
        &\mathbb{E} \big[ \mathbf{1} (\mathcal{L}_{k^*}(h_1)) \mathbf{1} (\mathcal{L}_{k^*}(h_2)) |F_{k^*,t_\sigma}(\sigma+ih_1)|^2 |F_{k^*,t_\sigma}(\sigma+ih_2)|^2 \big] \\
        &\ll \frac{\mathbb{E} \big[ \mathbf{1} (\mathcal{L}_{k^*}(h_1)) \mathbf{1} (\mathcal{L}_{k^*}(h_2)) |F_{k^*,t_\sigma}(\sigma+ih_1)|^2 |F_{M,t_\sigma}(\sigma+ih_2)|^2 \big]}{|h_1-h_2|^2 e^{2 k^*} (\log |h_1-h_2|^{-1} - k^*)^4 }.
    \end{split}\end{equation}
    Next, let 
    \[ n_1 := \left[2 \log |h_1-h_2|^{-1}\right] -k^* + C_3  \]
    for some sufficiently large integer constant $C_3$, and set $n_2 := t_\sigma$. One verifies that this choice satisfies $M \le n_1 \le n_2$ on this range of $|h_1 - h_2|$.
    Subtracting the barrier condition in $\mathcal{L}_{k^*}(h_1)$ at $k$ from the one at $n_1$, we obtain that under that $\mathcal{L}_{k^*}(h_1)$ we have
    \begin{align*}
        -C_1(k-k^*+1) - &(n_1 - k^*) + 2 \log (n_1 - k^*) \le \log |F_{n_1,k}(\sigma + i h_1)| \\
        &\le k-k^* - 2 \log (k-k^*) + C_1(n_1-k^*+1) \; \forall \; n_1 < k \le n_2,
    \end{align*}
    and similarly at $h_2$. Denoting by $\mathcal{L}_{n_1}(h_1,h_2)$ the event that both of these conditions hold, we deduce from \eqref{jun19.1} that
    \begin{equation}\begin{split}\label{jun19.3}
        &\mathbb{E} \big[ \mathbf{1} (\mathcal{L}_{k^*}(h_1)) \mathbf{1} (\mathcal{L}_{k^*}(h_2)) |F_{k^*,t_\sigma}(\sigma+ih_1)|^2 |F_{k^*,t_\sigma}(\sigma+ih_2)|^2 \big] \\
        &\ll \frac{\mathbb{E} \big[ \mathbf{1} (\mathcal{L}_{n_1}(h_1,h_2)) |F_{M,t_\sigma}(\sigma+ih_1)|^2 |F_{M,t_\sigma}(\sigma+ih_2)|^2 \big] \mathbb{E}[|F_{k^*,M}(\sigma+ih_1)|^2 ]}{|h_1-h_2|^2 e^{2 k^*} (\log |h_1-h_2|^{-1} - k^*)^4 } \\
        &= \frac{\tilde{\mathbb{P}}_{h_1,h_2}^{(2,2)} ( \mathcal{L}_{n_1}(h_1,h_2))  \mathbb{E}\big[|F_{M,t_\sigma}(\sigma+ih_1)|^2 |F_{M,t_\sigma}(\sigma+ih_2)|^2 \big] \mathbb{E}[|F_{k^*,M}(\sigma+ih_1)|^2 ]}{|h_1-h_2|^2 e^{2 k^*} (\log |h_1-h_2|^{-1} - k^*)^4 }.
    \end{split}\end{equation}
    Two applications of Lemma \ref{lem:MomentsEulerProducts} yield that this is in turn
    \begin{align*}
        &\ll \frac{\tilde{\mathbb{P}}_{h_1,h_2}^{(2,2)} ( \mathcal{L}_{n_1}(h_1,h_2)) |h_1-h_2|^2 e^{2t_\sigma} ((|h_1-h_2| e^{k^*})^{-1} }{|h_1-h_2|^2 e^{2 k^*} (\log |h_1-h_2|^{-1} - k^*)^4 } = \frac{e^{2 t_\sigma - 3 k^*} \tilde{\mathbb{P}}_{h_1,h_2}^{(2,2)} ( \mathcal{L}_{n_1}(h_1,h_2))}{|h_1 - h_2| (\log |h_1-h_2|^{-1} - k^*)^4}.
    \end{align*}
    Next, we want to apply Proposition \ref{GaussianComparisonGlobal} to bound this tilted probability. Note that the upper limit in the barrier condition satisfies
    \[ k-k^* - 2 \log (k-k^*) + C_1(n_1-k^*+1) \le \frac{e^{(k-k^*)/4}}{250(k-k^*)^2} \]
    for all $n_1 < k \le n_2$ by taking $C_3$ in the definition of $n_1$ sufficiently large, and similarly for the lower limit. It follows from \eqref{eq:centered_final} that
    \begin{align*}
         \tilde{\mathbb{P}}_{h_1,h_2}^{(2,2)} ( \mathcal{L}_{n_1}(h_1,h_2))&\ll \mathbb{P} \Big(\! - B(k-n_1) - O(n_1-k^*) \le \!\sum_{m = n_1+1}^k \overline{N}_m^{(1)} \le O(n_1-k^*) \; \forall \; n_1 < k \le n_2 \Big)  \\
        &\times\mathbb{P} \Big(\! - B(k-n_1) - O(n_1-k^*) \le \!\sum_{m = n_1+1}^k \overline{N}_m^{(2)} \le O(n_1-k^*) \; \forall \; n_1 < k \le n_2 \Big).
    \end{align*}
    Now by standard prime number theorem estimates, we also have $\frac{1}{20} \le \text{Var}(\overline{N}_{m}^{(1)}),\text{Var}(\overline{N}_{m}^{(2)}) \le 20$, say, so that \cite[Probability Result 2]{Harper1} finally implies that
    \[ \tilde{\mathbb{P}}_{h_1,h_2}^{(2,2)} ( \mathcal{L}_{n_1}(h_1,h_2)) \ll \frac{(n_1 - k^*)^2}{n_2 - n_1} \ll \frac{(\log |h_1-h_2|^{-1} - k^*)^2}{t_\sigma - k^*}. \]
    Plugging all of this into \eqref{jun10.2}, we obtain that the contribution from $h_1, h_2$ satisfying \eqref{jun19.2} is
    \begin{align*}
        &\ll \frac{e^{k^*}}{t_\sigma-k^*} \iint \frac{1}{|h_1 - h_2| (\log |h_1-h_2|^{-1} - k^*)^2} \mathrm{d}h_1 \mathrm{d}h_2 \ll \frac{e^{k^*}}{t_\sigma-k^*} \int_{e^{-(k^*+1)}}^{e^{-k^*}} \mathrm{d}h_1 \ll \frac{1}{t_\sigma-k^*}
    \end{align*}
    as required.
\end{proof}

\section{Quantitative central limit theorem}
We record the following variant of the Berry-Esseen theorem which is used throughout the paper (see, e.g., Corollary 17.2 in \cite{BhaRao}).

\begin{lem}\label{lem:BerryEsseen}
    Let $(W_j)_j$ be a sequence of independent random variables on $(\mathbb{R}, \mathcal{B}(\mathbb{R}),\mathbb{P})$, and $\mathbb{Q}_m$ be the distribution of ${W}_1+\dots+{W}_m$. For any $\mu\in \mathbb{R}$ and $\sigma \in \mathbb{R}$, let $\eta_{\mu,\sigma}$ denote the Gaussian measure with mean $\mu$ and variance $\sigma$. Lastly, for each $m\geq 1$, let
$$
\mu_m=\sum_{j=1}^m\mathbb{E}[{W}_j], \ \ \ {\sigma}_m=\sum_{j=1}^m \mathrm{Var}[{W}_j].
$$
Then if $\mathcal{A}$ is the set of intervals in $\mathbb{R}$, there exists an absolute constant $c>0$ such that 
$$
\sup_{A\in\mathcal A}|\mathbb{Q}_m(A)-\eta_{\mu_m,{\sigma}_m}(A)|\leq \frac{c}{\sigma_m^{3/2}}\sum_{j=1}^m\mathbb{E}\big[|{W}_j-\mathbb{E}[{W}_j]|^3\big].
$$
\end{lem}


\section{Proof of Proposition \ref{weighted_extended_Rademacher}}\label{sec:proof_weighted_extended_Rademacher}

\begin{proof}[Proof of Proposition \ref{weighted_extended_Rademacher}]
    Decomposing $n = s^2 m$ with $m$ squarefree, we have
    \[ \sum_{n \le x}\frac{X^e(n)}{\sqrt{n}} = \sum_{m \le x} \frac{X(m)}{\sqrt{m}} \sum_{s \le \sqrt{x/m}} \frac{1}{s}. \]
    We write
    \[ \sum_{s \le \sqrt{x/m}} \frac{1}{s} = \log \sqrt{x/m} + c_m \]
    with $c_m \ll 1$. Using subadditivity, this gives
    \[ \mathbb{E} \bigg[ \bigg| \frac{1}{\log x}\sum_{n \le x} \frac{X^e(n)}{\sqrt{n}} \bigg|^{2 q} \bigg] \ll \mathbb{E} \bigg[ \bigg| \sum_{m \le x} \frac{X(m)}{\sqrt{m}} \Big( 1 - \frac{\log m}{\log x} \Big) \bigg|^{2 q} \bigg] + \mathbb{E} \bigg[ \bigg| \frac{1}{\log x} \sum_{m \le x} \frac{ c_m X(m)}{\sqrt{m}} \bigg|^{2 q} \bigg]. \]
    For the second summand, note simply that, by H\"older's inequality and orthogonality, we have
    \[ \mathbb{E} \bigg[ \bigg| \frac{1}{\log x} \sum_{m \le x} \frac{ c_m X(m)}{\sqrt{m}} \bigg|^{2 q} \bigg] \ll \frac{1}{(\log x)^{2q}} \bigg( \sum_{m \le x} \frac{|c_m|^2}{m} \bigg)^q \ll (\log x)^{-q}.  \]
    This contribution is more than acceptable. It thus remains to show that
    \[ \mathbb{E} \bigg[ \bigg| \sum_{m \le x} \frac{X(m)}{\sqrt{m}} \Big( 1 - \frac{\log m}{\log x} \Big) \bigg|^{2 q} \bigg] \ll (\log x)^{q(2q-1)}. \]
    To this end, we point out that by Perron's formula, for any $\sigma > 0$ we have
    \[ \sum_{m \le x} \frac{X(m)}{\sqrt{m}} \Big( 1 - \frac{\log m}{\log x} \Big) = \frac{1}{2 \pi \log x} \int_{\mathbb{R}} F_t(\sigma+ih) x^{\sigma+ih} \frac{dh}{(\sigma+ih)^2}. \]
    Taking $\sigma = \frac{1}{\log x}$ and bounding $x^{\sigma+ih} \ll 1$, this gives
    \[ \mathbb{E} \bigg[ \bigg| \sum_{m \le x} \frac{X(m)}{\sqrt{m}} \Big( 1 - \frac{\log m}{\log x} \Big) \bigg|^{2 q} \bigg] \ll \mathbb{E} \bigg[ \bigg( \frac{1}{\log x} \int_{\mathbb{R}} \frac{|F_t(\sigma+ih)|}{|\sigma+ih|^2} \mathrm{d}h \bigg)^{2q} \bigg]. \]
    Subdividing the range of integration $e$-adically as before, this is in turn
    \[ \ll (\log x)^{2q} \mathbb{E} \bigg[ \bigg( \int_{0}^{\frac{1}{\log x}} |F_t(\sigma+ih)| \mathrm{d}h \bigg)^{2q} \bigg] +  \sum_{\substack{H \ge \frac{1}{\log x} \\ H \; e\text{-adic}}} \frac{1}{(\log x)^{2q} H^{4q}} \mathbb{E} \bigg[ \bigg( \int_{H}^{eH} |F_t(\sigma+ih)| \mathrm{d}h \bigg)^{2q} \bigg]. \]
    Next, we bound the contribution from $H \ge 1$. To do so, we simply apply H\"older's inequality and Lemmas \ref{lem:MomentsEulerProducts} and \ref{lem:PNT} to deduce that
    \[ \mathbb{E} \bigg[ \bigg( \int_{H}^{eH} |F_t(\sigma+ih)| \mathrm{d}h \bigg)^{2q} \bigg] \ll \bigg( H \int_{H}^{eH} \mathbb{E}[|F_t(\sigma+ih)|^2] \mathrm{d}h \bigg)^{q} \ll H^{2q} (\log x)^q.\]
    Summing over $H \ge 1$ gives a contribution of
    \[ \ll \frac{1}{(\log x)^{q}} \sum_{\substack{H \ge 1 \\ H \; e\text{-adic}}} H^{-2q} \ll (\log x)^{-q}. \]
    This is again more than acceptable for our final bound.

    For $\frac{1}{\log x} \le H \le 1$, we apply a balanced version of H\"older's inequality. We include the range $[0,(\log x)^{-1}]$ in this discussion, which can be bounded in the same way as the $H = (\log x)^{-1}$ case, with the same bound.

    More precisely, we have
    \begin{align*} 
    \mathbb{E} \bigg[ \bigg( \int_{H}^{eH} |F_t(\sigma+ih)| \mathrm{d}h \bigg)^{2q} \bigg] &= \mathbb{E} \bigg[ |F_{\log 1/H}(\sigma)|^{2q(1-2q)} \bigg( \int_{H}^{eH} |F_t(\sigma+ih)| |F_{\log 1/H}(\sigma)|^{-(1-2q)} \mathrm{d}h \bigg)^{2q} \bigg] \\
    &\le \mathbb{E}[|F_{\log 1/H}(\sigma)|^{2q}]^{1-2q} \bigg( \int_{H}^{eH} \mathbb{E}[|F_t(\sigma+ih)| |F_{\log 1/H}(\sigma)|^{-(1-2q)}] \mathrm{d}h \bigg)^{2q}
    \end{align*}
    We split the Euler product factors in the second expectation at $e^{1/H}$ and use independence to factor them into two expectations. Then applying Lemmas \ref{lem:MomentsEulerProducts} and \ref{lem:PNT} to estimate the resulting three expectations and integrating gives an upper bound of
    \[ \ll H^{-q(2q-1)(1-2q)} H^{2q} H^{-q(2q-1)(2q)} (H \log x)^{q/2} = H^{-q(2q-3)} (H \log x)^{q/2}. \]
    Summing the contribution over $(\log x)^{-1} \le H \le 1$ hence gives
    \[ \ll (\log x)^{-3q/2} \sum_{\substack{ \frac{1}{\log x} \le H \le 1 \\ H \; e\text{-adic}}} H^{-q(2q+1/2)} \ll (\log x)^{q(2q-1)} \]
    as required.
\end{proof}

\section{Formalisation in Lean 4}\label{lean}

All results were formalised in Lean 4 using Claude. This includes not only the results proved and written up in the present paper, but also all theorems the work depends on. This concerns

\begin{itemize}
    \item work of Harper \cite{Harper1} on low moments of random multiplicative functions. Here, we formalised both Theorems $1$ and $2$ as stated, with uniformity in $q$ in both the Steinhaus and Rademacher case.
    \item work of Caich \cite{Caich1}, proving an analogue of Harper's result for moments of partial sums on short intervals. In this case, we opted for a formalisation of the upper bound in the Rademacher setting.
\end{itemize}

The formalisation also draws inspiration from the works of Gorodetsky and Wong \cite{GorodetskyWongHelson} and the third author \cite{HamdanSupercritical} to simplify Harper's proof. Moreover, it uses ideas from \cite{HarperSoundXu} for minor simplifications in Caich's argument.

Aside from Mathlib, the formalisation also relies on the PNT+ project \cite{PNTPlus2024} for sieve estimates and a version of the Prime Number Theorem.

The versions of Theorems \ref{MainThm1_statement} and \ref{MainThm2_statement} formalised are slightly more precise than the ones proved in the present work. Namely, they also provide a lower bound with uniformity in $q$, showing that for $x \ge 3$ we in fact have
\[ \mathbb{E} \bigg[ \bigg| \frac{1}{\sqrt{x}}\sum_{n \le x} X^e(n) \bigg|^{2 q} \bigg]
        \asymp
        \begin{cases}
            \min \{ \log \log x, (1-2q)^{-1} \} (\log \log x)^{-q/2}, &\mbox{if } 0 \le q \le \frac{1}{2}, \\
            \min \{ \log \log x, (2q-1)^{-1} \}^{3/4} (\log x)^{q(2q-1)}, &\mbox{if } \frac{1}{2} \le q \le 1.
        \end{cases}, \]
        the same bounds having also been verified in the weighted Rademacher setting. The proof of the uniform lower bound proceeds along the same lines as those of the present paper, with only minor modifications. We have opted for a slightly weaker statement here, primarily for the sake of clarity and brevity.

The formalisation is available at
\begin{center}
    \url{https://github.com/MaximGerspach/low-moments-extended-rademacher}.
\end{center}
\noindent There, the file \texttt{RMF/Paper.lean} states the results above in Lean, and \texttt{RMF/PaperAxioms.lean} checks, as part of the build, that each of them depends only on the three standard axioms of Lean; in particular, no step of the formalisation is left unproved. The file \texttt{RMF/PaperExplicit.lean} contains versions with explicit constants: for instance, for every $x \ge 3$ and $0 < q \le 1$, the right-hand side above multiplied by $10^{-10}$ is a lower bound, and multiplied by $30$ an upper bound.

\bibliography{lit}
\bibliographystyle{alpha}

\end{document}